\documentclass[draft]{amsart}
\usepackage{amsmath,amsfonts,amssymb,amscd,verbatim,multicol}
\usepackage[arrow,matrix,cmtip]{xy}

\usepackage{fullpage}

\newcommand{\lann}{\ell\operatorname{-ann}}
\newcommand{\rann}{r\operatorname{-ann}}

\newcommand{\rel}{relevant\ }
\newcommand{\trdeg}{\operatorname{tr.deg}}
\newcommand{\End}{\operatorname{End}}

\newcommand{\CO}{\mathcal O}

\newcommand{\too}{\longrightarrow} 

\newcommand{\GKdim}{\operatorname{GKdim}}

\newcommand{\supp}{\operatorname{Supp}}

\newcommand{\bpr}{\begin{proof}}
\newcommand{\epr}{\end{proof}}

\newcommand{\ol}{\overline}

\newcommand{\spec}{\operatorname{Spec}}

\DeclareMathOperator{\calHom}{\mathcal{H}\mathit{om}}

\newcommand{\ann} {\operatorname{ann}}

\newcommand{\mc}{\mathcal}
\newcommand{\mf}{\mathfrak}

\newcommand{\mb}{\mathbb}

\newcommand{\dra}{\dashrightarrow}

\newcommand{\wt}{\widetilde}

\newcommand{\coH}{\operatorname{H}}

\newcommand{\im}{\operatorname{Im}}

\newcommand{\coh}{\operatorname{coh}}

\newcommand{\rGr}{\operatorname{Gr-}\hskip -2pt}
\newcommand{\rgr}{\operatorname{gr-}\hskip -2pt}
\newcommand{\rQgr}{\operatorname{Qgr-}\hskip -2pt}
\newcommand{\rqgr}{\operatorname{qgr-}\hskip -2pt}

\newcommand{\rtors}{\operatorname{tors-}\hskip -2pt}

\newcommand{\aut}{\operatorname{Aut}}

\DeclareMathOperator{\Hom}{Hom}

\DeclareMathOperator{\HB}{H}

\newcommand{\dirlim}{\underrightarrow{\lim}}

\newcommand{\bbc}{birationally commutative}

\newcommand{\sets}{\mathbf{Sets}}
\newcommand{\kalg}{\mathbf{Alg}}

\newcommand{\kschemes}{\mathbf{Schemes}}

\newcommand{\naive}{na{\"\i}ve}

\newcommand{\nsr}{na{\"\i}ve blowup algebra}
\newcommand{\Nsr}{Na{\"\i}ve blowup algebra}
\newcommand{\Splus}{S^+}
\newcommand{\Sminus}{S^-}
\newcommand{\Splusex}{S^+_{\infty}}
\newcommand{\Sminusex}{S^-_{\infty}}

\newcommand{\wtY}{\mathbb{X}}
\newcommand{\YY}{Y}
\newcommand{\XX}{W}
\newcommand{\bg}{birationally geometric}
\newcommand{\Bg}{Birationally geometric}
\newcommand{\btwog}{birationally $2$-geometric}
\newcommand{\XXX}{X}
\newcommand{\PPP}{P}
\newcommand{\wh}{\widehat}

\newcommand{\LL}{\mathcal{L}}
\newcommand{\OO}{\mathcal{O}}

\newcommand{\II}{\mathcal{I}}
\newcommand{\JJ}{\mathcal{J}}
\newcommand{\GG}{\mathcal{G}}
\newcommand{\PP}{\mathcal{P}}
\newcommand{\QQ}{\mathcal{R}}
\newcommand{\FF}{\mathcal{F}}
\newcommand{\HH}{\mathcal{H}}

\newcommand{\RR}{\mathcal{R}}
\renewcommand{\AA}{\mathcal{A}}
\newcommand{\NN}{\mathcal{N}}

\numberwithin{equation}{section}
 \theoremstyle{plain}
\newtheorem{theorem}[equation]{Theorem}
\newtheorem{maintheorem}[equation]{Main Theorem}

\newtheorem{sublemma}[equation]{Sublemma}
\newtheorem{lemma}[equation]{Lemma}

\newtheorem{corollary}[equation]{Corollary}
\newtheorem{proposition}[equation]{Proposition}

\newtheorem{assumptions}[equation]{Assumptions}

\theoremstyle{definition}
\newtheorem{definition}[equation]{Definition}
\newtheorem{notation}[equation]{Notation}

\newtheorem{remark}[equation]{Remark}
\newtheorem{remarks}[equation]{Remarks}

\newtheorem{example}[equation]{Example}

\begin{document}
\title[noncommutative projective surfaces]{A class of  noncommutative projective surfaces}
\author{D. Rogalski and J. T. Stafford}

\address{(Rogalski)
Department of Mathematics, UCSD, La Jolla, CA 92093-0112, USA. }
\thanks{The first author was  partially supported by the NSF  through grants
DMS-0202479  and  DMS-0600834   while the second author was
  partially supported by the NSF through grants DMS-0245320 and
  DMS-0555750  and also  by the  Leverhulme Research Interchange Grant
F/00158/X.
 Part of this work was written up while the second author was visiting and supported by the Newton Institute, Cambridge. We would like to thank all three institutions for their financial support.}
\email{rogalski@math.mit.edu}

\address{(Stafford) Department of Mathematics, University
of Michigan, Ann Arbor, MI 48109-1043, USA.}
\email{jts@umich.edu}
\begin{abstract}
Let $A=\bigoplus_{i\geq 0}A_i$ be a connected graded,
noetherian $k$-algebra that is generated in degree one  over an algebraically closed
 field $k$.
Suppose that  the graded quotient ring $Q(A)$ has the form
$Q(A)=k(X)[t,t^{-1},\sigma]$, where $\sigma$ is an automorphism of the
integral projective surface $X$.  Then we prove that  $A$ can be written
 as a na{\"\i}ve blowup algebra  of a projective surface $\mathbb{X}$  birational to
 $X$. This enables one to obtain a deep understanding of the structure of these algebras;
  for example, generically
  they are not strongly noetherian and their point modules are
  not parametrized   by a projective scheme. This is  despite the fact that
  the simple objects in $\rqgr A$  will always be in 
(1-1) correspondence with the closed
   points of the    scheme~$\mathbb{X}.$
\end{abstract}
\maketitle

\tableofcontents
\clearpage

 \section{Introduction}
One of the outstanding problems of noncommutative projective
geometry is to classify all noncommutative projective  surfaces.
There are many ways of phrasing this problem; from a ring-theoretic
standpoint it means classifying (hopefully in some geometric way)
those connected graded rings  that should be thought of as
corresponding to surfaces. In this paper, we will achieve just such
a classification for a large class of noncommutative rings:
noetherian  connected graded rings that are \btwog, as defined
below.

\subsection{The results}   Let us first give the statement of
 the main theorem and then explain its significance and applications.
Throughout, $k$ will be an algebraically closed base field. A
$k$-algebra $A$ is called  \emph{connected graded}
(\emph{cg})\label{cg-defn}   if $A=\bigoplus_{n\geq 0} A_n$, where
$A_0=k$ and  $\dim_kA_n<\infty$ for each $n$. For  a noetherian cg
domain $A$, one can invert the nonzero homogeneous elements to
obtain the \emph{graded quotient ring} $Q(A)\cong
D[t,t^{-1};\sigma]$, for some automorphism  $\sigma$ of the division
ring $D=Q(A)_0$ (see also  Section~\ref{sect-ideals}). We say that
\emph{$A$ is \bg} if $D\cong k(V)$, the field of rational functions
on an integral projective scheme $V$, and $\sigma$ is induced from
an automorphism of $V$.  If in addition $V$ is
a surface, we say that \emph{$A$ is \btwog}.

The aim of the paper is to classify \btwog\     algebras  as   \nsr s,  in the sense of  \cite{KRS}.
These are described as follows.    Let $\XXX$ be an integral projective
scheme, with  zero dimensional subscheme $Z=Z_{\mc{I}}=\spec
\mc{O}_\XXX/\II$,   automorphism $\sigma$ and $\sigma$-ample invertible sheaf
$\LL$, as defined in Section~\ref{naive-section}.
 In a manner reminiscent
of the  Rees ring   construction of the  blowup of a commutative scheme,
 we form the \emph{bimodule algebra}
$$\RR=\RR(\XXX, Z, \mc{L},\sigma) =\mc{O}_\XXX\oplus \mc{J}_1\oplus\mc{J}_2\oplus\cdots,\
\quad\mathrm{where}\  \mc{J}_n = \mc{L}_n\otimes_{\mc{O}_\XXX}\mc{I}_n \ \mathrm{for}\
\II_n = \mc{I} \cdot  \sigma^* \mc{I} \cdots (\sigma^*)^{n-1}\mc{I}.$$ As is described in detail in
 \cite{KRS} and explained briefly here in Section~\ref{naive-section}, this
  bimodule algebra has a natural
multiplication and the  {\it \nsr\ of $\XXX$ at $Z$} is then defined to be the algebra of sections
  $$R=R(\XXX,Z, \LL,\sigma)
  =\mathrm{H}^0(\XXX, \, \RR) = k\oplus \mathrm{H}^0(\XXX, \, \mc{J}_1)
  \oplus\mathrm{H}^0(\XXX, \, \mc{J}_2)\oplus\cdots $$

Although implicitly   $Z$ is nonempty,
 this definition  makes perfect sense when $Z=\emptyset$, in which case
 $B(\XXX,\mc{L},\sigma)=R(\XXX,\emptyset,\mc{L},\sigma)$ is just
  the \emph{twisted homogeneous coordinate ring}  that is so important in other
aspects of noncommutative projective geometry; see, for example,
\cite{AV} or  \cite{SV}.   For any zero dimensional
subscheme $Z\subset \XXX$ there is a natural embedding $R=R(\XXX,Z,
\LL,\sigma)\hookrightarrow B=B(\XXX,\LL,\sigma)$ and the multiplication
in $R$ is also that induced from $B$.

We can now state the main result of this paper.

\begin{maintheorem}\label{intro-mainthm} { \rm [Theorem~\ref{mainthm}]}
Let $A$ be a cg noetherian domain that is generated in degree one
and \btwog.  Then, up to a finite dimensional vector space, either
$A\cong B(\wtY,\LL,\sigma)$ or $A\cong R(\wtY,Z,\LL,\sigma)$, for
some projective surface $\wtY$ with automorphism $\sigma\in
\mathrm{Aut}(\wtY)$, $\sigma$-ample invertible sheaf $\mc{L}$ and  
zero dimensional subscheme $Z$.
 \end{maintheorem}

  Noetherian \nsr s only appear when the algebras are far from commutative; explicitly for each
 closed point  $z\in Z$ the orbit $\langle \sigma\rangle\cdot z $ will be  \emph{critically dense} in the
 sense that  the orbit is infinite and all of its infinite subsets are Zariski dense in $\wtY.$

   Combined with the results from \cite{KRS} and \cite{RS},
   the  alternative $A\cong R(\wtY,Z,\LL,\sigma)$ in Theorem~\ref{intro-mainthm}
  has  remarkable consequences for the structure of the algebra $A$ and
  its representations.
Before giving  these results we need some further definitions.
      The category of noetherian  $\mb{Z}$-graded right $A$-modules, with homomorphisms
      being graded homomorphisms of degree zero, will
be denoted by $\rgr A$.\label{grade-defn}
 Let $\rtors A$ denote the full subcategory  of $\rgr A$ consisting of
 modules of finite length,
  with quotient category $\rqgr A=\rgr A/\hskip -2pt\rtors A$.\label{qgrade-defn}
 (Thinking   of $\rqgr A$ as  the category of
coherent sheaves on the imaginary scheme $\mathrm{Proj}\, A$ can often provide
useful intuition.)
 \emph{A point module for $A$} is a cyclic graded
$A$-module $M=\bigoplus_{n\geq 0} M_n$  such that  $\dim_k M_n=1$
for all $n$. A \emph{point module in $\rqgr A$} is defined to be the
   image in $\rqgr A$ of a cyclic graded $A$-module $M=\bigoplus_{n\geq 0} M_n$,
generated in degree zero,  such that $\dim_k M_n=1$ for all $n\gg
0$.  For point modules defined over arbitrary base rings and the
formalities on parametrization, see
 Section~\ref{sect-pts}.

 \begin{corollary}\label{intro-maincor}
 Let $A$ be as in Theorem~\ref{intro-mainthm} and keep the notation from
 that result.
\begin{enumerate}
\item[{\rm (I) }]   If $A\cong B(\wtY,\LL,\sigma)$ then:
 \begin{enumerate}
   \item $A$ is  \emph{strongly noetherian};  i.e.,
   $A\otimes_kC$ is noetherian
for all commutative noetherian $k$-algebras~$C$;\label{strong-noeth}
\item  $\rqgr A \simeq \coh \wtY$, the category of coherent sheaves on $\wtY$;
 \item the set of point modules for $A$ (both in $\rgr A$
 and in $\rqgr A$) is  parametrized by the scheme $\wtY$;
  \item $A$ has a balanced dualizing complex in the sense of \cite{Ye}.
\end{enumerate}\end{enumerate}
\begin{enumerate}
\item[{\rm (II)}] Otherwise $A\cong R(\wtY,Z,\LL,\sigma)$, for
some nonempty   finite  subscheme $Z$ and:
\begin{enumerate}
\item   $A$ is not strongly noetherian;
\item neither the  point modules in $\rgr A$ nor those in   $\rqgr A$ are   parametrized
by a scheme of locally finite type;
 \item   the simple objects in $\rqgr A$ are the point modules and
 are naturally in  $(1$-$1)$ correspondence with the closed points of $\wtY;$
 \item  $\rqgr A$ has finite cohomological dimension;
  \item $A$ does not have a balanced dualizing complex;
 \item $H^1(A)=\mathrm{Ext}^1_{\rqgr A}(A,A)$ is infinite dimensional;
   \item $A$ does not satisfy generic flatness in the sense of \cite{ASZ}.
  \end{enumerate}
 \end{enumerate}
 \end{corollary}

{\bf Proof of Corollary~\ref{intro-maincor}:}  \label{maincor-proof}
This is really just a case of quoting the literature. To be precise,
 part (Ia) follows from \cite[Proposition~4.13]{ASZ}, while
 part (Ib) is \cite[Theorem~1.3]{AV} and
(Id) is proved in \cite[Theorem~7.3]{Ye}.  Part (Ic) follows from
\cite[Theorem~1.1, 1.2]{RZ} and \cite[Proposition~10.2]{KRS}.
   All the results stated in part (II) are contained in
\cite[Theorem~1.1 and Remark~1.2]{RS}.     \qed

   We  should emphasize the striking  fact that the properties  described
by part (II)  of the corollary are not exceptional:  as soon as    $(\wtY,\sigma)$
has at least one  critically dense orbit then the
 theorem and its corollary imply that,  generically,
each noetherian cg subalgebra of $k(\wtY)[t,t^{-1};\sigma]$ has these properties.
 Despite the curious features displayed in  Corollary~\ref{intro-maincor}(II), the category  $\rqgr A$
described there  is surprisingly similar to the category of coherent sheaves on a projective scheme.
   This is illustrated  by  the following result, which also shows that the properties of
the \naive\ blowup   $\rqgr A$  are actually closer to those of the original scheme $\wtY$
 than they are to those of the classical blowup~$\wt{\wtY}$    of $\wtY$ at $Z$.

\begin{proposition}\label{something} {\upshape
\cite[Theorem~4.10]{RS}}   Keep the hypotheses  of  Theorem~\ref{intro-mainthm}.
 Then the subcategory of  (Goldie) torsion objects in $\rqgr R(\wtY, Z, \mc{L},\sigma)$  is
equivalent to the category of
 torsion  coherent $\mc{O}_\wtY$-modules.\end{proposition}

\subsection{An example}\label{first-eg}
Here is an easy example that illustrates the sorts of
algebras that appear in Theorem~\ref{intro-mainthm} and its corollary
(see, also, Example~\ref{simplest-eg}).

 A typical example of a twisted coordinate ring is the ``quantum polynomial ring''
$S=S_{pq}$ with generators $x,y,z$ and relations $xy-pyx,\,\,
xz-qzx$ and $ yz-qp^{-1}zy$ for  $p$ and $q\in k^*$. This can also
be written as   $S \cong B(\mb{P}^2,\mc{O}_{\mb{P}}(1), \sigma)$,
where   $\sigma\in \mathrm{Aut}( \mb{P}^2)$ is defined by
$(\lambda_0:\lambda_1 : \lambda_2)\mapsto
(\lambda_0:p\lambda_1:q\lambda_2) .$ It follows   that
$Q(S)_0=k(yx^{-1}, zx^{-1})=k(\mb{P}^2)$ and so, in this special
case, Theorem~\ref{intro-mainthm} classifies the noetherian cg
algebras $A$ with $Q(A)=Q(S)$. To see an explicit  example of case
(II) of the corollary, assume that the scalars $p,q$ are
algebraically independent over the base field, take $Z$ to be the
single reduced point $(1:1:1)$ and consider the ring $R(\mb{P}^2, Z,
\mc{O}(1), \sigma)$; it is easy to show this ring is equal to $k
\langle x-y, x-z \rangle \subseteq S$.

   It is also easy to find cg subalgebras
  of $Q(S_{pq})$ that are not noetherian. For example,  \cite[Proposition~4.8]{KRS}
   implies that, for this choice of $\sigma$, the \nsr\
 $ R(\mb{P}^2, Z', \mc{O}(1),\sigma)$ is not noetherian when $Z'=\{(1:0:0)\}$.
  For a detailed and more algebraic examination of these rings, see \cite{Ro1}.

 \subsection{History}\label{subsection-history}
  We briefly explain the history behind these results and their wider relevance.

  As we remarked earlier, it is useful to think    of $\rqgr A$ for a noetherian cg ring
  $A$ as coherent sheaves on the nonexistent scheme $\mathrm{Proj}\,A$.
  Extending this analogy further
 one can think of  a noncommutative projective integral curve, respectively surface, as
  $\rqgr A$ for a noetherian cg domain  $A$ with Gelfand-Kirillov dimension    two,
   respectively three. For the purposes of this discussion we will also restrict ourselves to algebras
   generated in degree one (see Section~\ref{subsection-questions} for comments
   about the general case).   Under this analogy, noncommutative integral projective curves
have been classified in \cite{AS1} (they are nothing more than
$\rqgr B(\XXX,\mc{L},\sigma)\simeq \coh \XXX$ for integral projective curves $\XXX$). Similarly,
      noncommutative projective planes and their algebraic analogues, the Artin-Schelter regular rings
     with the Hilbert series of a polynomial ring in three variables,
      have been classified in   \cite{ATV1, BP}.  Such an Artin-Schelter regular ring
      either equals $B(\mb{P}^2,\mc{O}(1),\sigma)$ for
some $\sigma$ or   has a factor isomorphic to $B(E,\mc{L},\sigma)$ for  a projective curve
 $E$ in which   case $(E,\sigma)$ determines $A$. These results are surveyed
 in \cite{SV}, and the reader is referred to that paper  for more details.

Motivated by these results, in \cite{AR1} Artin posed the problem of classifying all
noncommutative integral projective surfaces   and this paper was motivated in part
by that question.   In one sense, however,   Theorem~\ref{intro-mainthm} is   orthogonal
to his question, since he was
 interested in finding the noncommutative algebras whose properties are close to those of
 commutative integral surfaces. For example, he
included the existence of a balanced dualizing complex as one of his hypotheses!  Of course,
 if $A$ satisfies the hypotheses of Theorem~\ref{intro-mainthm}  and also has a balanced
 dualizing complex,  then Corollary~\ref{intro-maincor}(e)  implies  that $A$ is a
 twisted homogeneous coordinate ring, in which case  the category  $\rqgr A$ of coherent sheaves on
  the corresponding noncommutative projective scheme   is commutative:  $\rqgr A\simeq  \coh
  \wtY$.

 Further motivation for this paper comes from recent work on strongly noetherian algebras.
 These algebras are studied in detail in \cite{ASZ, AZ2} and have many pleasant properties,
  some of which  are described in Section~\ref{subsection-proofs} below. One is therefore interested in
 understanding  when cg algebras satisfy this property and for \bg\ algebras
  this is answered by the following result.

\begin{theorem}\label{intro-RZtheorem}
 {\upshape(1)  \cite[Theorem 1.2]{RZ}}
Suppose that  $A$  is a strongly noetherian,  cg domain that is   generated in degree one and
 \bg.    Then,  up to a finite
dimensional vector space, $A\cong B(\wtY,\LL,\sigma)$, for a
projective variety $\wtY$, with automorphism $\sigma\in \aut(\wtY)$ and
$\sigma$-ample invertible sheaf $\LL$.

{\upshape (2) \cite[Proposition~4.13]{ASZ}}  Conversely, if $\mc{L}$
is a $\sigma$-ample invertible sheaf on a projective scheme $\XXX$,
then $B(\XXX,\mc{L},\sigma)$  is strongly noetherian.\end{theorem}

 Since noetherian cg domains that are \emph{not} strongly noetherian exist in profusion, even among
 domains of Gelfand-Kirillov dimension three, it
 is important to understand their structure and Theorem~\ref{intro-mainthm} gives
 one step in that direction.

 \subsection{The proofs}\label{subsection-proofs}   We next want to  give the strategy behind
  the proof of the main theorem, but we first outline the proof of
  Theorem~\ref{intro-RZtheorem}(1)
   as the contrast is illuminating.   Strongly noetherian rings are known
to have a number of very pleasant properties and in particular one has the following
 result from \cite[Theorem~E4.3]{AZ2}: \emph{If $A$ is a strongly noetherian
 $k$-algebra then the point
  modules for $A$ are parametrized by a projective scheme $\wtY.$}
  The key step in the proof of Theorem~\ref{intro-RZtheorem}
 is that, under the hypotheses of the theorem  and up to a finite
dimensional vector space as always, \emph{there exists a surjection
 $\chi: A\twoheadrightarrow B(\wtY,\mc{L},\sigma)$ for this scheme $\wtY$ together with
 a sheaf $\mc{L}$ and automorphism $\sigma$}.
By hypothesis  $Q(A)_0=k(V)$ for some scheme $V$ and the next step
is to  regard $Q(A)_{\geq 0}=\sum_{n\geq 0}k(V)t^n$ as a  ``generic
 point module'' with coefficients from $k(V)$ (see Section~\ref{sect-ideals} for a
  discussion of this type of module).  The existence of this generic point module  can  then
  be used to show that $\chi$ is injective.

Now suppose that $A$ is merely noetherian. As we have seen, the
point modules for $A$ are no longer parametrized by a projective
scheme  and so there is no easy candidate for the scheme $\wtY;$
indeed the construction  of $\wtY$  forms a major part of the proof
of the theorem.
  The idea is as follows.   {\it
  A truncated  point module
of length $n+1$} is  a cyclic graded $A$-module
$M=\bigoplus_{i=0}^n M_i$ where    $\dim_kM_i=1$ for each $i$. By
\cite[Prop\-osition~3.9]{ATV1} these modules are parametrized by a
scheme $\XX_n$ and there are natural maps $\Phi_{n-1},\Psi_{n-1} :
\XX_{n}\to \XX_{n-1}$ given by further truncation:
$\Phi_{n-1}(M)=M/M_n$, respectively $\Psi_{n-1}(M)=M_{\geq 1}[1]$,
where $[1]$ denotes the shift in gradation. When $A$ is strongly
noetherian, it follows from  \cite[Corollary~E4.5]{AZ2} that
$\Phi_n$ and $\Psi_n$ become isomorphisms for   $n\gg 0$.
 The scheme   $\wtY$ of Theorem~\ref{intro-RZtheorem}  is then just  $\wtY=\XX_n$ for
such an $n$, with automorphism $\sigma=\Psi_n\Phi_n^{-1}$.  In
contrast, when $A$ is only noetherian these maps will  never be
isomorphisms. Fortunately there is a reduced and irreducible
component $\YY_n$ of $\XX_n$ such that (a) $\Phi_n$ and $\Psi_n$
restrict to maps between the $\YY_n$ and (b) $k(\YY_n)\cong k(V)$
and so there is   an analogous generic  truncated point module
$\bigoplus_{i = 0}^n Q(A)_i$ (see Corollary~\ref{zero-ideal}).

The main step   the proof of Theorem~\ref{intro-mainthm} is to show
that, for $n\gg 0$, one can appropriately blow down   $\pi: \YY_n\to
\wtY$  to  give a scheme $\wtY$  for   which the induced action of
the birational map $\sigma=\Psi_n\Phi_n^{-1}$ is actually  an
automorphism. This is the scheme appearing in the theorem.
The scheme   $\wtY$ appears in a second way, which
also indicates how it is constructed: given a truncated point module
$M'$ corresponding to some $y\in \YY_n$  and $n\gg 0$, write  $M'$
as the image of a (non-unique) point module $M$ and  take the image
$\ol{M}$ of $M$ in $\rqgr A$. As  Corollary~\ref{ptmodule-cor}
shows, it is  \emph{this}  object that is uniquely determined by
$M'$
 and the map $\pi$ corresponds to mapping  $M'\mapsto \ol{M}$.
 The proof of these assertions
  takes up Sections~\ref{sect-pts}---\ref{construction-section}. Once  $\wtY$ has been
constructed, one can use the ideas from \cite{AS1, KRS, RZ}  to prove  the theorem
(see Sections~\ref{naive-section}---\ref{main-section}).

\subsection{Questions}\label{subsection-questions}
We end the introduction with some comments and questions. First,
Theorem~\ref{intro-RZtheorem}  is actually stronger than stated since it works for any
cg noetherian domain  $A$ that is \emph{birationally commutative} in the sense that
 $Q(A)_0$ is a finitely generated field; in other words
one does not need to assume that the automorphism $\sigma$ is induced from
that of a scheme.  We conjecture that, for surfaces,
  Theorem~\ref{intro-mainthm} also holds
under this more general hypothesis, but there is one case we are as yet  unable to handle.
Using work of Diller and Favre  \cite{DF} which, essentially,
  classifies birational maps of  surfaces one has the following result.

\begin{theorem}\cite{Ro3}
 Let $Q = K[t, t^{-1}; \sigma]$, where $K$
is a finitely generated field extension of $k$ of transcendence
degree two.
  Then every cg Ore domain $A$ with graded quotient ring
$Q$ has the same Gelfand-Kirillov dimension  $d \in \{3,4,5, \infty
\}$.  If $d < \infty$, then $d = 3$ or $5$ if and only if $\sigma$
is induced from an automorphism of some projective surface $\XXX$ with
$K\cong k(\XXX)$. If $d=\infty$ then $A$ is not
noetherian.\end{theorem}

This just leaves the case when $d=4$. This definitely can occur; for
example take the automorphism $u\mapsto uv, v\mapsto v$ of
$k(\mb{P}^2)=k(u,v)$. We believe, however, that when $d=4$ the
algebra $A$ can never be noetherian, which would solve  the
conjecture. See  Example~\ref{non-noeth}  for a typical    example
in this dimension.  Note, also, that the case $d=5$ is covered by
Theorem~\ref{intro-mainthm} and so noncommutative surfaces are a
little more general than algebras of Gelfand-Kirillov dimension $3$.

We end with two further questions. First, as is true of
Theorem~\ref{intro-RZtheorem}, is there a version of
Theorem~\ref{intro-mainthm} describing any \bg\ $A$? Second,   in the present paper we
 always assume that our algebras are generated in degree one.
  In the commutative case this is not so much of a restriction since one can
  always reduce to that case by taking an appropriate Veronese ring. In the
  noncommutative case, however, it is a more severe restriction, since it excludes
  the idealizer examples that appear in \cite{Ro2, AS1, SZ}. It would be interesting
  to know whether, as happened for curves in \cite{AS1},  one can classify all
  noetherian cg algebras that are \btwog\ in terms of
  variants of \nsr s and idealizers.

\section{Point modules}\label{sect-pts}

The method  we  use  to introduce geometry
into the study of noncommutative algebras is the theory of
(truncated) point modules and the  schemes that parametrize
them.  In this section we recall the relevant definitions  and prove some basic
facts about these objects that will be needed in the sequel.
 These concepts were first introduced by Artin, Tate and Van
den Bergh in \cite[Section 3]{ATV1}, and more details can be found both
there and in \cite[Section~4]{RZ}. For the most part our geometric notation
 is standard and, unless
we say otherwise,   we will use the definitions given in \cite{Ha}.

As in the introduction, fix a  cg algebra $A=\bigoplus_{n\geq 0}
A_n$ with categories $\rgr A$ and $\rqgr A$. The category of all
graded right $A$-modules
 is written $\rGr A$,  but we will not need the corresponding
 quotient category $\rQgr A$.
Assume for the rest of the section that $A$ is generated by $A_1$ as
a $k$-algebra. Given a commutative $k$-algebra $R$  we write $A_R=A\otimes_kR$, regarded
as a graded $R$-algebra by putting $R$ into degree zero. Following
\cite[Definition~3.8]{ATV1}, a \emph{truncated $R$-point module of
length $n+1$}\label{trunc-defn}, for $n\geq 0$,  is a graded cyclic
right $A_R$-module $M = \bigoplus_{i = 0}^n M_i$ where $M_0=R$ and
each $M_i$  is a locally free $R$-module of rank $1$. An
\emph{$R$-point module}\label{point-defn} is given by the same
definition, except that now $n=\infty$. Under the identification of
$R$ with $R^{\mathrm{op}}$, we can and usually  will  regard
(truncated) right $R$-point modules as $(R,\, A)$-bimodules since
this will frequently  make their module structures more transparent.
Note that if $M$ is a (truncated) $R$-point module, then the
identification $M_0=R$ defines a unique isomorphism $M\cong  A_R/I$
for some right ideal $I$, and so this rigidifies the given
(truncated) $R$-point module.
 A (truncated) $k$-point module is just called  a (truncated) point module.
  Write $\PP(R)$, respectively $\PP_n(R)$, for the
set of isomorphism classes of right $R$-point modules, respectively
truncated right $R$-point modules of length $n + 1$.
 We write
$\PP=\PP(k)$, $\PP_n=\PP_n(k)$ and let $\PP^\ell(R),\dots,
\PP^\ell_n$ denote the left hand analogues of each of these sets.

Recall that Yoneda's Lemma embeds the category $\kschemes$ of
noetherian $k$-schemes inside the category of functors  from the
category $\kalg$ of commutative noetherian $k$-algebras to the
category $\sets$ of  sets, where a scheme $\XXX$ corresponds to the
functor $h_{\XXX}: R \mapsto \Hom(\spec R, \XXX)$ (see
\cite[Section~VI]{EH}).  A functor $h: \kalg\to \sets$   is said to
be \emph{represented} by the scheme $\XXX$ if $h$ and $h_\XXX$ are
isomorphic functors.

We will apply this theory to the  functors arising from  truncated point
modules over our algebra $A$.  By \cite[Proposition~3.9]{ATV1}
the functor $F_n: \kalg\to \sets$
 defined by   $R\mapsto \PP_n(R)$
  is represented by a (commutative) projective
scheme $\XX_n$ called the \emph{$n^{\mathrm{th}}$ truncated point
scheme} for  $A$.\label{X-defn} Since it will be used many times, we
recall the functorial property  of  $F_n$.  Thus, let  $M \in
\mc{P}_n(R)$ be a truncated $R$-point module corresponding to a
morphism of schemes $\theta : \spec R\to \XX_n$
   and let $\alpha: R \to S$ be a homomorphism of commutative
  noetherian $k$-algebras. The  induced map $ \spec S \to \spec R$
  will be denoted $\alpha^\star$.\label{star-defn}
Using our convention that $R$-point modules are written as
$(R,\,A)$-bimodules, $S\otimes_R M$ is naturally a truncated point
module over $A_R\otimes_RS \cong A_S$ and the morphism of schemes
$\spec S\to \XX_n$ corresponding to
 $S \otimes_R M$ is  simply $\theta\alpha^\star$.
If we need to remember the map $\alpha$ in this construction, we
will write $S{{}_\alpha\otimes}\, M$  for
$S\otimes_RM$.\label{twisted-tensor-defn}

Following \cite{ATV1} the schemes $\XX_n$ can be constructed
explicitly. Write $A \cong T(A_1)/J$,  where $J=\bigoplus_{n\geq 0}
J_n $ is a homogeneous ideal of the tensor algebra $T(A_1)$.  Let
$\mb{P}$ \label{pee-defn}  denote  the projective space
$\mb{P}(A_1^*)$ of lines in $A_1^*$ and for each $n \geq 1$ write
$\mb{P}^{\times n}$ for   the product of $n$ copies of $\mb{P}$. For
$n\geq 0$, we identify $T_n$ with  the global sections of the sheaf
$\mc{O}(1,1, \dots, 1)$ on $\mb{P}^{\times n}$.  Under this
identification, $J_n$ defines a vector subspace $\wt{J}_n$ of
$\HB^0(\mb{P}^{\times n}, \mc{O}(1,1, \dots, 1) )$
and, by the proof of   \cite[Proposition~3.9]{ATV1},  $\XX_n$  can be  identified
   with the vanishing locus of $\wt{J}_n$  in $ \mb{P}^{\times n}$.
We typically write this embedding as    $\iota_n: \XX_n
\hookrightarrow \mb{P}^{\times n}$.
For each $n$ and $1 \leq a \leq b \leq n$
let $\pi_{a,b} : \mb{P}^{\times n} \to \mb{P}^{\times
(b-a+1)}$\label{pi-ab-defn} denote the projection onto the
coordinates between $a$ and $b$ inclusively. By \cite[Proposition
3.5]{ATV1}, $\pi_{1,n}(\XX_{n+1}) \subseteq \XX_n$ and
$\pi_{2,n+1}(\XX_{n+1}) \subseteq \XX_n$ and we define projective
morphisms $\Phi_n: \XX_{n+1} \to \XX_n$ and $\Psi_n: \XX_{n+1} \to
\XX_n$\label{phi-defn} by restricting $\pi_{1,n}$, respectively
$\pi_{2,n+1}$ to $\XX_{n+1}$. We call the collection of schemes and
morphisms $ (\XX_n,   \Phi_n, \Psi_n)$\label{X-data-defn} the
\emph{point scheme data} for the algebra $A$, where the embeddings
$\iota_n$  are to be understood.

 The
morphisms $\Phi_n$ and $\Psi_n$ are closely related to truncations
of modules. Let $R$ be a commutative noetherian $k$-algebra and, for
an $R$-module $V$, set $V^*=\Hom_R(V,R)$.
If $M=\bigoplus M_i\in \rGr A$ and  $n\in \mathbb{Z}$, then the  {\it
shift}\label{shift-defn} $M[n]$ of $M$ is defined by
$M[n]_r=M_{n+r}$.
  Given a truncated
$R$-point module of length $n + 2$, say $M=\bigoplus_{i=0}^{n+1} M_i
\in \mc{P}_{n+1}(R)$, we define the \emph{right truncation}   of $M$ to be
\begin{equation}\label{rtrunc}
 \Phi'_n(M)\ = \
M/M_{n+1} \ = \  M_0 \oplus\cdots\oplus M_{n} \ \in\  \mc{P}_n(R),
\end{equation}
and we define the \emph{left truncation}
of $M$ to be the right
$A_R$-module
\begin{equation}\label{ltrunc}
\Psi'_n(M)\ = \  M_1^*\otimes_RM[1]_{\geq 0}  = \
(M_1^*\otimes_RM_1)\oplus\cdots\oplus (M_1^*\otimes_RM_{n+1})\  \in\
\mc{P}_n(R).
\end{equation}
 The  canonical identification
$ M_1^* \otimes_R M_1 = R$ ensures that   $\Psi'_n(M)$ really is a
truncated $R$-point module.

By repeating the previous few paragraphs  with left modules in place  of
right modules, we can define the analogous  \emph{left point scheme
data}\label{left-defn} $(\XX^{\ell}_n,   \Phi^{\ell}_n,
\Psi^{\ell}_n)$. For example, $\XX^{\ell}_n$ is a projective scheme
which represents the functor $F^{\ell}_n : \kalg\to \sets$ defined
by  $R\mapsto \mc{P}^{\ell}_n(R)$. As is shown by the next lemma,
 $(\XX^{\ell}_n,   \Phi^{\ell}_n,  \Psi^{\ell}_n)$    is a mirror of $ (\XX_n,   \Phi_n,
\Psi_n)$   and so
we will not need to consider it separately.

For the rest of this section, we analyze the morphisms $\Phi_n$ and $\Psi_n$
and begin with some elementary properties that will be used frequently.

\begin{lemma} \label{right left duality}
{\rm (1) }  For $n\geq 0$, the maps $\mc{P}_{n+1}(R) \to
\mc{P}_n(R)$ given by the rules $M  \mapsto \Phi'_n(M)$, and $M
\mapsto \Psi'_n(M)$ induce  the respective morphisms $\Phi_n:
\XX_{n+1} \to \XX_n$ and $\Psi_n: \XX_{n+1} \to \XX_n$ between the
representing schemes.

{\rm (2)} There are isomorphisms $\gamma_n: \XX_n
\buildrel{\sim}\over\longrightarrow  \XX^{\ell}_n$ for $n \geq 0$
satisfying $\gamma_n^{-1} \Phi^{\ell}_n \gamma_{n+1} = \Psi_n$ and
$\gamma_n^{-1} \Psi^{\ell}_n \gamma_{n+1} = \Phi_n$.

{\rm (3)}  For all $n\geq 0$ we have $\Psi_{n}\Phi_{n+1}=\Phi_{n}\Psi_{n+1}$.
\end{lemma}

\begin{proof}
(1)  See the proof of  \cite[Proposition~3.9]{ATV1}.

(2)  Fix $n \geq 0$ and suppose that $M = \bigoplus_{i = 0}^n M_i\in
\PP_n(R)$, for some noetherian commutative $k$-algebra $R$.  Define
$M^{\vee} =
\Big(\bigoplus_{i = -n}^0 M_{-i}^* \otimes_R M_n\Bigr)[-n]$ as a variant of the Matlis dual
and note that $M^\vee$
is a graded left $A_R$-module, under the convention that $A$ acts
 from the left and $R$ from the right.
 By  \cite[Proposition~10.2(1)]{KRS}\footnote{The proof
   in \cite{KRS} omitted the term $\otimes M_n$ which is need to ensure that
   $M^\vee_1=R$, but otherwise the proofs are identical.}  $M^{\vee} \in \PP_n^\ell (R).$

Similarly, if $N\in \PP_n^\ell(R) $, then  $N^{\vee} =\Bigl(
\bigoplus_{i = -n}^0 N_n\otimes_RN_{-i}^*\Bigr)[-n]\in \PP_n(R)$.
The two operations are inverses, compatible with change of rings,
and so the map $M\mapsto M^\vee$ induces an
 isomorphism  $\gamma_n: \XX_n
\buildrel{\sim}\over{\longrightarrow}  \XX^{\ell}_n$ between the
representing schemes.

Given $M = \bigoplus_{i = 0}^{n+1} M_i\in \PP_{n+1}(R)$, then a
short calculation shows that
\[
((M^{\vee})_{\leq n})^{\vee} \cong \bigoplus_{i = 0}^n M_1^*\otimes_R M_{i+1}
=\Psi_n'(M)  \in \mc{P}_n(R),
\]
and this isomorphism is compatible with change of rings.  Using part
(1) and  its left-sided analog  we conclude that $\gamma_n^{-1}
\Phi^{\ell}_n \gamma_{n+1} = \Psi_n$ as maps from $\XX_{n+1} $ to $
\XX_n$. The proof that $\gamma_n^{-1} \Psi^{\ell}_n \gamma_{n+1} =
\Phi_n$ is analogous.

(3)  The  morphisms $\pi_{1,n}\, \pi_{2, n+2}$ and $\pi_{2, n+1}\,
\pi_{1, n+1}$ give the same projection  $\mb{P}^{\times (n+2)} \to
\mb{P}^{\times n}$. Restriction to the subscheme $\XX_n \subseteq
\mb{P}^n$ gives $\Phi_n \Psi_{n+1} = \Psi_n \Phi_{n+1}$.
\end{proof}

In order to  avoid any  possible ambiguity we will clarify our  conventions
about a morphism
of schemes  $f: Z\to W$.
If we say that $f^{-1}$ is defined at a point $w\in W$, then we
do mean that it is defined as a morphism from an open neigbourhood of $w$
to an open neighbourhood of the (necessarily singleton)  set  $f^{-1}(w)$.
As such, $f^{-1}$ automatically
defines a local  isomorphism  at $w$, meaning that it defines an isomorphism from an open
 neighbourhood of $w$ to an open neighbourhood of  $f^{-1}(w)$.
    Finally,  we adopt the usual  convention   that an  automorphism
     $\tau$   acts on functions by $(\tau(\theta))(z)=\theta(\tau(z))$ for $z\in Z$.

\begin{lemma}  \label{point fiber}
 Let $f: Z \to W$ be a  proper morphism of finite type schemes and
 suppose that  $w \in W$ is a closed
point with set-theoretic fibre  $f^{-1}(w)
= \{z \}$  a single
point. Then:
\begin{enumerate}
\item If the induced map $f_z^{\sharp} : \mc{O}_{W,w} \to \mc{O}_{Z,z}$ is   an isomorphism,
then $f^{-1}$ is  a local  isomorphism  at  $w$.
 \item If $f_*(\OO_Z)=\OO_W$, then  $f^{-1}$ is a local  isomorphism at  $w$.
 \end{enumerate}
\end{lemma}

\begin{proof}   (1) It is routine to show that, since
$f_z^{\sharp}$ is an isomorphism, there are open neighborhoods $w \in V
\subseteq W$ and $z \in U \subseteq Z$ such that $f$ restricts to an
isomorphism $U
\to V$.  Since the fibre   $f^{-1}(w)
=\{z\}$ is a singleton and $f$ is a closed morphism,
we can shrink $V$ to a smaller open set $V' \ni w$ whose
inverse image satisfies $f^{-1}(V') \subseteq U$.  Then $f$ is
injective over $V'$ and  $f^{-1}$ is defined at $w$.

(2) Since $f_*(\mc{O}_Z) = \mc{O}_W$, the map $f_z^{\sharp}$ can be
identified with the map
$$\widehat{f} : \underset{w \in V}{\dirlim} \mc{O}_Z(f^{-1}(V)) \to
 \underset{z \in U}{\dirlim} \mc{O}_Z(U),$$
where $U$ and $V$ run through open sets in their respective schemes.
As in part $(1)$, for any such $U$ containing $z$ there is an open
set $V$ containing $w$ with $f^{-1}(V) \subseteq U$.  It follows
that $\widehat{f}$ is an isomorphism and so the result follows from
part~(1).
\end{proof}

It follows from \cite[Theorem~E4.3]{AZ2} together with
Lemma~\ref{right left duality} that the point scheme data $(\XX_n,
\Phi_n, \Psi_n)$ stabilizes for strongly noetherian rings, in the
sense that $\Phi_n$ and $\Psi_n$ are isomorphisms for large $n$. In
contrast, as we  remarked in the introduction, this need not be the
case for arbitrary noetherian rings. Fortunately, as we show in
Proposition~\ref{eventually iso},
 the point scheme data
does at least stabilize locally, and this will form the starting point of many of
the results in this paper.

In the next few results,  we  let $\chi$\label{chi-defn} be a symbol
which denotes either $\Phi$ or $\Psi$; thus any property claimed for
a map $\chi_n$ is supposed to hold for both $\Phi_n$ and $\Psi_n$.

\begin{proposition}  \label{eventually iso}
Assume that  $A$ is a noetherian, cg $k$-algebra  generated in
degree $1$ and pick $n_0\geq 0$.  Let $\{w_n \in \XX_n : n \geq n_0
\}$ be a sequence of  (not necessarily closed)  points  such that
$\chi_n(w_{n+1}) = w_n$ for all $n \geq n_0$. Then, for all $n\gg
n_0$,
the fibre   $\chi_n^{-1}(w_n)
= \{w_{n+1}\}$ is a singleton  and  $\chi_n^{-1}
$ is a local  isomorphism  at $w_{n}$.
\end{proposition}

\begin{proof} Assume that $\chi = \Phi$; the proof for $\chi = \Psi$ follows
similarly using left modules and Lemma~\ref{right left duality}(2).

Suppose first that the points $\{ w_n \}$ are closed points.  Then
this sequence $\{w_n\}$  corresponds to  a right $k$-point module $M
= A/J$; thus  $J$ is a  right ideal  with $\dim_k J_n = \dim_k A_n -
1$ for all $n \geq 0$. Choose $m_0$
such that $J$ is generated in degrees $\leq m_0$.  If possible, pick
$m\geq m_0$  such that the fibre $\Phi_m^{-1} (w_m)$  contains (at
least) two distinct closed points, say  $w_{m+1}$  and $y$,
corresponding to the truncated point modules $M_{\leq (m+1)}$ and
$N$. Write $N = A/I$ for  some right ideal $I$ and notice that
$I_{\leq m} = J_{\leq m}$, but $I_{m+1} \neq J_{m+1}$. However,  as
$A$ is generated in degree one and $J$ is generated in degrees $\leq
m$, this implies that $I_{m+1} \supseteq I_m A_1 = J_m A_1 =
J_{m+1}$. Since $\dim_kA_{m+1}/I_{m+1}=1=\dim A_{m+1}/J_{m+1}$, this
forces $I_{m+1} = J_{m+1}$, a contradiction.

Thus $\Phi_n^{-1}(w_n)    = \{w_{n+1}\}$ for all $n \gg 0$.
  By   \cite[Proposition~3.6(ii)]{ATV1} and  the last line of
  the proof of that result,  the induced map of
noetherian local rings $\widehat{\Phi}_n : \mc{O}_{\XX_n, w_n} \to
\mc{O}_{\XX_{n+1}, w_{n+1}}$ is surjective for any such $n$.
   For sufficiently large $n$, the map $\widehat{\Phi}_n$
is therefore an isomorphism.  By  \cite[Corollary~II.4.8(e)]{Ha}
$\Phi_n$ is proper and  so   Lemma~\ref{point fiber}(1)  implies
 that $\Phi_n^{-1}$  is a local  isomorphism  at $w_{n}$.

Now suppose that the points $w_n \in \XX_n$ are not necessarily
closed and write $\overline{w}_n$ for the closure of $w_n$. By
properness, $\Phi_n(\overline{w}_{n+1})$ is closed, and thus equal
to $\overline{w}_n$, for any $n\geq n_0$. Thus, beginning with any
closed point $y_{n_0}\in \overline{w}_{n_0}$, we can choose a
sequence of closed points
 $y_n\in \overline{w}_n$ such that
$\Phi_n(y_{n+1}) = y_n$, for each $n \geq n_0$.   For $n\gg 0$     the first part of
the proof shows that   the fibre $\Phi_n^{-1}(y_n)
=\{y_{n+1}\}$ is a singleton and that
$\Phi_n^{-1}$ is a local  isomorphism    at $y_{n}$.
 But then $\Phi_n^{-1}$ is then also  local  isomorphism   at
  $w_{n}$  and so certainly  $\Phi_n^{-1}(w_n)
= w_{n+1}$ is a singleton  for all such $n$.
\end{proof}

  For $m>n\geq 0$ define
morphisms $\Phi_{m,n}$ and $\Psi_{m,n}$ by
\begin{equation}\label{phi-mn}
\Phi_{m,n}=\Phi_n\Phi_{n+1}\cdots\Phi_{m-1} : \XX_m\to \XX_n\qquad\
\mathrm{and}\qquad \Psi_{m,n}=\Psi_n\Psi_{n+1}\cdots\Psi_{m-1} :
\XX_m\to \XX_n.
\end{equation}

\begin{corollary} \label{finite sets}  Fix $n_0\geq 0$.
For  each $n\geq n_0$ let $G_n$ be a finite set of (not necessarily
closed) points of $\XX_n$ such that  $\chi_n(G_{n+1}) \subseteq
G_n$. Then for $n \gg n_0$, the cardinality $|G_n|$ is   constant
and $\chi_n$ gives a bijection between $G_{n+1}$ and $G_n$.
Moreover, if $g\in G_n$ for any such  $n$,  then $\chi_n^{-1}$ is a
local  isomorphism  at $g$.
\end{corollary}

\begin{proof} As before, it suffices to prove the result for $\chi = \Phi$.
For each $n\geq n_0$, let $H_n $ be the subset of points $w\in G_n$
for which there exists  $m> n$ and two distinct points  $y,z\in
G_{m}$ such that $\Phi_{m,n}(y)=w=\Phi_{m,n}(z)$.
 Obviously $\Phi_n(H_{n+1}) \subseteq H_n$, so if
$H_{n+1}\not=\emptyset$ for some $n>n_0$, then $H_n\not=\emptyset$.
Assume that  $H_n\not=\emptyset$ for all $n \geq n_0$. Then we claim
that there exists a sequence of points $\{w_n \in H_n : n\geq n_0\}$
such that $\Phi_n(w_{n+1}) = w_n$  for all $n \geq n_0$. To see this,
for $w\in H_n$, set $\ell(w) = \sup\{r : w=\Phi_{r+n,n}(y) $ for
some $y\in H_{n+r}\}$. As $H_{n_0}$ is finite we may pick $w_{n_0}\in H_{n_0}$
with $\ell(w_{n_0})$ maximal. If $\ell(w_{n_0})=t<\infty$, then $H_{j}=\emptyset$
for all $j>t+n_0$,  contradicting our assumption.  So
$\ell(w_{n_0})=\infty$.  Now we may construct the sequence $w_n$
inductively:   having chosen $w_n \in H_n$ with $\ell(w_n) = \infty$,
since $H_{n+1}$ is finite there exists $w_{n+1}\in H_{n+1}$ such
that $\Phi_n(w_{n+1})=w_n$ and $\ell(w_{n+1})=\infty$.

Thus if $H_n\not=\emptyset$ for all $n\geq n_0$,  there is a
sequence of points $w_n \in \XX_n$ such that (i)  $\Phi_n(w_{n+1}) =
w_n$ for all $n \geq n_0$, and (ii) for infinitely many choices of
$n$ the set theoretic  fibre  $\Phi_n^{-1}(w_n)$ is not a singleton.
 This contradicts Proposition~\ref{eventually
iso}. So $H_n$ is empty for $n \gg n_0$ and $\Phi_n$ induces an
injection from $G_{n+1}$ to $G_n$ for all such $n$. Since each $G_n$
is finite, this implies that $\Phi_n$ must induce a bijection  for $G_{n+1}$ to $G_n$
for all $n \gg n_0$.  The final statement of  the
corollary now follows from the final assertion of
Proposition~\ref{eventually iso}.
\end{proof}

\section{Birationally commutative domains}
\label{sect-ideals}

Fix  a noetherian cg $k$-algebra $A$  that is  generated in degree
one, with point scheme data $(\XX_n,\Phi_n,\Psi_n)$. In this section
we examine how the ideal structure of $A$ is reflected in the
varietal structure of the schemes $\XX_n$.
 For the algebras $A$  of interest,
this will allow us to isolate subschemes $\YY_n\subseteq
\XX_n$  that better reflect the geometry of $A$.

   If $A$ is also a domain then \cite[C.I.1.6 and A.I.4.3]{NV} shows that
    the set $\mc{C}$  of nonzero homogeneous
elements of $A$ forms  an Ore set with corresponding
localization $Q(A)=A_\mc{C} \cong D[t,t^{-1}; \tau]$. Here $\tau$ is
an automorphism of the division ring $D$, and $D[t,t^{-1}; \tau]$ is the twisted  Laurent
polynomial ring defined by $td= d^\tau t$ for all $d\in D$.
We call $Q(A)$ the \emph{graded quotient ring of
$A$} and \emph{$D$ the function ring of $A$}. \label{fn-ring-defn}
If $D$ is a field, we call $A$ \emph{\bbc}, and note that in this case $D$ will be a finitely
generated field extension of $k$; see, for example, Corollary~\ref{zero-ideal}(1).\label{bbc-defn}
 (One reason for the notation is that, if $A$ were the homogeneous
coordinate ring of an integral projective scheme  $Y$, then $D$ would just
be the function field $k(Y)$.)

We begin with a useful observation about the point scheme data for
factor rings.

\begin{lemma}
\label{inverse image stable} Let $A$ be a noetherian cg $k$-algebra,
generated in degree $1$,  with a  homogeneous ideal $I$.  Let
$(\XX_n, \Phi_n, \Psi_n )$ be the point scheme data for $A$, and let
$(\widehat{\XX}_n, \widehat{\Phi}_n, \widehat{\Psi}_n)$ be the point
scheme data for $A/I$.   Then, for all $n\geq 0$, $\widehat{\XX}_n$ is
canonically identified with a closed subscheme of $\XX_n$ with
$\wh{\Phi}_n=\Phi_n\vert_{\wh{\XX}_n}$ and
$\wh{\Psi}_n=\Psi_n\vert_{\wh{\XX}_n}$, respectively.

Moreover, as sets, $\Phi_n^{-1}(\wh{\XX}_n) = \wh{\XX}_{n+1}$ and
$\Psi_n^{-1}(\wh{\XX}_n) = \wh{\XX}_{n+1}$ for all $n \gg 0$.
\end{lemma}

\begin{proof}
We are given a presentation $A = T(A_1)/J$ so that $\XX_n \subseteq
\mb{P}^{\times n}=\mb{P}(A_1^*)^{\times n}$ is the zero set
$\mc{V}(\wt{J}_n)$  of the multi-linearization $\wt{J}_n$  of
   $J_n$. Now, $A/I \cong T/K$ for some ideal  $K\supseteq J$ and so
   $\wh{\XX}_n$ is defined to be a subscheme of $\mb{P}\bigl((A_1/K_1)^*\bigr)^{\times n}$.
    However,  we can identify  $\mb{P}\bigl((A_1/K_1)^*\bigr)^{\times n}$  with a closed subscheme
   of $\mb{P}(A_1^*)^{\times n}$ via the embedding
   $(A_1/K_1)^* \hookrightarrow A^*_1$, after which   $\wh{\XX}_n$ is the zero set of the larger
   vector space $\wt{K}_n\supseteq \wt{J}_n$. Thus $\wh{\XX}_n$ is indeed closed in
   $\XX_n$.
   (If one considers the corresponding point modules, the embedding
   $\wh{\XX}_n\hookrightarrow \XX_n$ is nothing more than the identification
   of point modules over $A/K$ with the
point modules over $A$ annihilated by $K$.) Since $\Phi_n$ and
$\Psi_n$ are defined as restrictions of projection maps
$\mb{P}^{\times (n+1)} \to \mb{P}^{\times n}$, the maps
$\wh{\Phi}_n$ and $\wh{\Psi}_n$ are just further restrictions of
those maps.

Assume that $I$ is generated, as a right $A$-module, in degrees
$\leq n_0$;  thus  $I_nA_1 = I_{n+1}$ for  $n \geq n_0$.  For $n\geq
n_0$, let $z \in \wh{\XX}_n$ and $y \in \XX_{n+1}$ be closed points
with $\Phi_n(y) = z$.     Thus $z$ corresponds to a truncated point
module $A/L$  where $L$ is a right ideal of $A$ containing $I$ and
$y$ corresponds to a truncated point module $A/M$ with $M_{\leq n} =
L_{\leq n}\supseteq I_{\leq n}$.  As $n\geq n_0$ this implies that
$M \supseteq M_{\leq n} A \supseteq I$ and hence that $y \in
\wh{\XX}_{n+1}$. In other words,  $\Phi_n^{-1}(\wh{\XX}_n) =
\wh{\XX}_{n+1}$ for all $n \geq n_0$. The assertion
$\Psi_n^{-1}(\wh{\XX}_n) = \wh{\XX}_{n+1}$ follows by working with
left modules and using Lemma~\ref{right left duality}(2).
\end{proof}

If $A$ is a strongly noetherian domain which is birationally commutative,
then it is isomorphic in large degree to a twisted homogeneous
coordinate ring $B = B(X, \mc{L}, \sigma)$ (see Section~\ref{subsection-questions}
or  \cite[Theorem 1.2]{RZ})
and  there is a natural correspondence between ideals of $A$ and
$\sigma$-invariant subschemes $Z\subseteq X$
 (see \cite[Lemma~4.4]{AS1}).
If $A$ is merely noetherian then the correspondence, which forms the content of
the next proposition,  is more subtle, basically because the morphisms
$\Psi_n$ and $\Phi_n$ need not be isomorphisms for $n \gg 0$.

We begin by discussing the point modules defined by function fields of factors of $A$.
Let $B=\bigoplus B_n$ be a \bbc\ cg
algebra that is generated in degree one, with graded quotient ring $Q(B)=  K[t,t^{-1};\sigma]$
 for some  field $K$  and $\sigma\in \mathrm{Aut}(K)$. Since $B_1\not=0$ we may
  assume that $t\in B_1$, as is often convenient in these situations.
  Regard
$\mc{G}=K[t;\sigma]=Q(B)_{\geq 0}$ as a right $B_K$-module, where $B$ acts by
right multiplication but $K$ acts from the left. Since $t\in B_1$, $\mc{G}$ is cyclic and hence
 is a $K$-point module for $B$.  Now consider the shift
 $\mc{G}[1]_{\geq 0}$ of $\mc{G}$.
It is immediate that the map
$\mu: \mc{G}\to \mc{G}[1]_{\geq 0}$,  given by
$t^m\ell \mapsto t^{m+1}\ell$ for $\ell\in K$, is an isomorphism of graded right $B$-modules.
However, this is not an isomorphism of $B_K$-modules since the $K$-action is
 twisted: given $a\in K$ and $g=t^{m}\ell\in \mc{G}_{\geq 0}$ for some $\ell\in K$, then
 the left $K$-action is  given by
$$a\circ \mu(g) = a t^{m+1}\ell = t^{m+1}a^{\sigma^{-(m+1)}}\ell =
 \mu\left(a^{\sigma^{-1}}\circ (t^m\ell) \right)  = \mu(a^{\sigma^{-1}}\circ g).
 $$
There is a second way of writing the module  $\mc{G}[1]_{\geq 0}$.
In the notation from page~\pageref{twisted-tensor-defn},
consider the homomorphism $\sigma: K=R\to K=S$ and the resulting
$B_S$-module $K{{}_\sigma\otimes}\,\mc{G}  =
S\otimes_R\mc{G}$.
(We have introduced $R$ and $S$ in order to distinguish between the two copies of $K$.)
As right $B$-modules, $\mc{G} \cong S\otimes_R \mc{G}$
under the map $\lambda:g\mapsto 1\otimes g$.
 However, the $S$-module action is again twisted:
 if $s * r=s r^\sigma $ denotes the action of $r\in R$ on $s\in S$ then
$$
a\circ \lambda(g) = a\otimes g =  (1* a^{\sigma^{-1}})\otimes g =
1\otimes a^{\sigma^{-1}}g =\lambda(a^{\sigma^{-1}} \circ g)
\quad\mathrm{for}\ a\in K=S\ \mathrm{and}\ g\in \mc{G}.
$$
Comparing the two  displayed equations shows that:

\begin{lemma}\label{K-shift}
Keep  the notation from the previous paragraph. Then
$\mc{G}[1]_{\geq 0} \cong K{{}_\sigma\otimes}\,\mc{G}$ as $A_K$-modules. \qed
\end{lemma}

\begin{proposition}
\label{ideals}  Let $A$ be a noetherian cg $k$-algebra that is  generated in degree $1$.
Then there is a bijection between  the  sets:
\begin{enumerate}
\item[(i)] Completely prime, homogeneous ideals $P$ of $A$ (excluding $A_{\geq 1}$)
such that $A/P$ is \bbc.
\item[(ii)] Sequences of (not necessarily closed) points $\{ \eta_n \in \XX_n :   n \gg 0 \}$ such that
$\Phi_n(\eta_{n+1}) = \Psi_n(\eta_{n+1}) = \eta_n$ for all $n \gg0$,
where two such sequences are identified if they agree for all large
$n$.
\end{enumerate}
Under this correspondence the function ring $Q(A/P)_0$ in part (i)
  is isomorphic to the   field
 $k(\eta_n)$  in part (ii).
\end{proposition}

\begin{proof} \  \textbf{Step 1}. Suppose that $P$ is a completely prime
ideal as in (i).  Let $\ol{A}= A/P$ and write $Q(\ol{A}) \cong
K[t,t^{-1}; \sigma]$, where $K$ is a field and $t\in \ol{A}_1$.
  The idea is to think of  $\mc{G}=Q(\ol{A})_{\geq 0}$  \label{generic-module}
  as a ``generic" point module for $\ol{A}$.  More precisely,
  regard $\mc{G}$ as a $K$-point module for  $A$, as  discussed above.
  For $n\geq 1$, the module  $\mc{G}_{\leq n} = \bigoplus_{i =
0}^n Kt^i$ is a truncated $K$-point module of length $n+1$ for $A$
and so corresponds to a morphism $\theta_n: \spec K\to \XX_n$. By
Lemma~\ref{right left duality}(1),   $\Phi_n\theta_{n+1}=\theta_n$
for all $n\gg 0$.

Lemma~\ref{right left duality}(1) also implies that the left
truncation  $\Psi_n'(\mc{G}_{\leq n+1}) =\bigoplus_{i = 0}^n
\mc{G}[1]_i$ corresponds to the morphism $\Psi_n\theta_{n+1}: \spec
K \to \XX_n$ for $n\gg 0$.  By Lemma~\ref{K-shift},
$\Psi_n'(\mc{G}_{\leq n+1})  \cong K {{}_\sigma\otimes}\,
\mc{G}_{\leq n}$ and so $\Psi_n\theta_{n+1} = \theta_n\sigma^\star$,
where $\sigma^\star : \spec K\to \spec K$ is the morphism induced by
$\sigma$. If  $\eta_n=\theta_n(\eta) \in \XX_n$ denotes  the image
of the point $\eta=\spec K$ then  $\Psi_n (\eta_{n+1}) =
\Psi_n\theta_{n+1}(\eta) = \theta_n\sigma^\star(\eta) =
\theta_n(\eta)=\eta_n. $ On the other hand, the previous paragraph
shows that $\eta_{n}=\theta_n(\eta)=\Phi_n\theta_{n+1}(\eta) =
\Phi_n(\eta_{n+1})$ for all $n \gg0$. Thus $\{ \eta_n \}$ is a
sequence of points as in (ii).
\medskip

\noindent \textbf{Step 2}. Conversely, suppose that we are given
points $\eta_n \in \XX_n$ as in (ii), say for $n \geq n_0$.  By
Proposition~\ref{eventually iso}, we can enlarge $n_0$ if necessary
so that,   for all $n \geq n_0$, both $\Psi_n^{-1}$ and
$\Phi_n^{-1}$ are defined and hence local isomorphisms  at $\eta_n$.
\ Define $L = k(\eta_{n_0})$ to be the residue field at
$\eta_{n_0}$; this identification also induces a morphism
$\rho_{n_0}: \spec L \to \XX_{n_0}$ with $\im \rho_{n_0} =
\eta_{n_0}$.  For $n\geq n_0$,  we can inductively define
$\rho_{n+1}=\Phi_n^{-1}(\rho_n)$; thus
 $\rho_{n+1}: \spec L \to \XX_{n+1}$  satisfies
 $\im \rho_{n+1} = \eta_{n+1}$, and  $\rho_n=\Phi_n(\rho_{n+1})$ for all $n \geq n_0$.
 This means that
there is an $L$-point module $M$ for $A$ which corresponds to the   maps $\rho_n$.

Since  $\Phi_n^{-1}$ and $\Psi_n^{-1} $
are defined at $\eta_n$ for $n \geq n_0$,   our hypotheses imply that
$\Psi_n \Phi_n^{-1}$ is a local
isomorphism  at  $\eta_n$ with $\Psi_n
\Phi_n^{-1}(\eta_n) = \eta_n$.
Then there is a  unique scheme automorphism $\wt{\tau}_n: \spec L \to \spec L$
satisfying $\rho_n \wt{\tau}_n = \Psi_n \Phi_n^{-1} \rho_{n} =
\Psi_n \rho_{n+1}$.  We claim that  the morphisms $\wt{\tau}_n$ are
all
 equal to each other. Indeed, for  $n \geq n_0$,  Lemma~\ref{right
left duality}(3) implies  that
\begin{equation}
\label{commuting calculation}
\rho_n \wt{\tau}_n = \Psi_n \rho_{n+1}
= \Psi_n \Phi_{n+1} \rho_{n+2} = \Phi_n \Psi_{n+1} \rho_{n+2} =
\Phi_n \rho_{n+1} \wt{\tau}_{n+1} = \rho_n \wt{\tau}_{n+1}.
\end{equation}
Thus  $\wt{\tau}_n = \wt{\tau}_{n+1}=\wt{\tau}$, say.

By  Lemma~\ref{right left duality}(1), the truncations
  $N=\Phi_{n}'(M)=M_{\leq n+1}/M_{n+1}$ and $N'=\Psi'_{n}(M_{\leq n+1})
  =(M_{\leq  n+1}[1])_{\geq 0}$
 of $M_{\leq n+1}$  correspond to the   scheme morphisms
$\Phi_n\rho_{n+1}=\rho_n $ and $\Psi_n\rho_{n+1}=\rho_n\wt{\tau}$, respectively.
By  the functoriality of $F_n$, as discussed on page~\pageref{twisted-tensor-defn},
  $N'\cong L{{}_\tau \otimes}\, N$,
where   $\tau:L\to L$ is the morphism induced
 from $\wt{\tau}$. Thus
 $M[1]_{\geq 0}\cong   L{{}_\tau \otimes}\,  M$.
Let $P= \operatorname{ann}_A M_s$ for some $0\leq s$. As the
$A$-module and $L$-module actions on $M$ commute, $P$ is also the
annihilator of  $(L{{}_\tau\otimes}\, M)_s = M_{s+1}$.
Hence, by induction, $P=\operatorname{ann}_A
M $. Note that, as $M$ is a point module, $P \neq A_{\geq  1}$.

We will show that $P$ is a completely prime ideal such that $A/P$ is
\bbc, as in (i).  Suppose first that $a\in A_u$ and $b\in A_v$ are
homogeneous elements with $a,b\not\in P$. Thus $M_0\, a\not=0$
which, since $M_u$ is a $1$-dimensional $L$-module, implies that
$M_0 a= M_u$. Similarly, as $P=\mathrm{ann}_A M_u$, we find that
$M_0 a b=M_u b=M_{uv}$. Thus, $ab\not\in P$ and $P$ is a completely
prime ideal of $A$.   Set $\overline{A}=A/P$ and write
   $Q=Q(\overline{A})=Q_0[t,t^{-1};\gamma]$,  for some $\gamma\in \mathrm{Aut}(Q_0)$
 and $t\in \overline{A}_1$.  Given
nonzero homogeneous elements $m \in M_s$, $a \in \overline{A}_u$ and
$b \in \overline{A}_v$ with $u \geq v$ then,  once again,
 $M_{s+u} = Lma = M_{s+u-v}b$ and so
there exists $n \in M_{s + u - v}$  with $n b = m a$.  Then $m
\otimes a b^{-1} = n \otimes 1$  in $M \otimes_A Q$, and so $M \cong
M \otimes_A Q_{\geq 0}$  already has a right $Q_{\geq 0}$-module
structure. As $Q_0$ is a division ring and its right action on $M$
commutes with that of $L$, there is an injection of rings $f : Q_0
\hookrightarrow \End_L M_0=L$. Thus $Q_0$ is   isomorphic to a
subfield $K$ of $L$. In particular, $\overline{A}$ is \bbc.

\medskip
\noindent \textbf{Step 3}. We need to show that the operations from
Steps~1 and 2 are mutual inverses.   To begin, we assume that we
have chosen a sequence $\{ \eta_n \}$ as in (ii), and performed
Step~2 to construct a completely prime ideal $P$ as in $(i)$.  Keep
all of the notation of that step.

As $L$ is a field, the right action of $Q_0$ on $M_0$ coincides with
the left action of the
 field $K\subseteq L$.
 Thus, for  any $0
\neq m \in M_0$ the $A_L$-module surjection $L\otimes_kQ_{\geq 0} \to M$
given by $\ell\otimes q\mapsto \ell m q$ induces
 a graded surjective $A_L$-module homomorphism
$\alpha: L {\otimes}_K Q_{\geq 0} \to M$.
Since  $Q_{\geq 0} = Q_0[t;\gamma]$,
 each graded piece both of  $L\otimes_KQ_{\geq 0}$ and of $M$   is isomorphic to $L$.
Thus $\alpha$ is an isomorphism.

We now perform Step~1 with the ideal $P$ to construct maps
$\theta_n: \spec K \to \XX_n$ corresponding to the $K$-point module
$Q_{\geq 0}$. If $f^\star: \spec L \to \spec K$ is the morphism
associated to the injection of fields $f:K\hookrightarrow L$ then
the isomorphism $L {{}_f \otimes} Q_{\geq 0} \cong M$ implies that
$\rho_n$ factors as $\theta_n f^\star = \rho_n$ for all $n \geq
n_0$.  Thus $\eta_n=\im(\theta_n)$ for all $n\geq n_0$ and so
performing Step~2 followed by Step~1 gives the identity operation.

This also proves the final assertion of the proposition. Indeed, as
the morphism $\rho_{n_0}$ is defined by the identification
$k(\eta_{n_0}) = L$ the factorization  $\theta_{n_0} f^\star  =
\rho_{n_0}$ forces $f$ to be the equality  $K = L$.  Thus Step~2 maps the set of
points $\{\eta_n\}$ to a completely prime ideal $P$ such that $A/P$
has function ring $K=L= k(\eta_n) $ for $n \gg 0$.

\medskip
\noindent \textbf{Step 4}.
 Conversely, suppose we have chosen an ideal $P$ as in (i)
 and performed Step~1, to obtain a sequence of maps $\theta_n: \spec K \to \XX_n$
corresponding to the point module $Q(A/P)_{\geq 0} = K[t; \sigma]$,
together with points $\eta_n = \im(\theta_n) \in \XX_n$ as in (ii).
Then we may perform Step~2 to construct a point module $M$. As in
Step~2, write   $L = k(\eta_{n_0})$ with associated morphisms
$\rho_n:\spec L\to \XX_n$ for $n\geq n_0$. In this case
$\theta_{n_0}$ factors as $\theta_{n_0} = \rho_{n_0} \wt{f}$ for  a
unique morphism of schemes $\wt{f}: \spec K \to \spec L$, with
corresponding ring map $f: L \to K$.  Thus $ \Phi_{n, n_0} \theta_n
= \theta_{n_0} = \rho_{n_0} \wt{f} = \Phi_{n, n_0} \rho_n \wt{f} $
for all $n \geq n_0$. By the choice of $n_0$,    $\Phi_{n, n_0}^{-1}
$ is a local  isomorphism  at   $\eta_{n_0}$  and so  $\theta_n =
\rho_n \wt{f}$ for all $n \geq n_0$.   By  functoriality  this
implies that $K {}_f \otimes M \cong Q(A/P)_{\geq 0}$ as
$A_K$-modules, and hence that
 $\mathrm{ann}_A M =
\ann_A (Q(A/P)_{\geq 0}) = P$.  Thus performing Step~1 followed by
Step~2 is also the identity operation.\end{proof}

One of the most important  cases of Proposition~\ref{ideals}
occurs when $P$ is the zero ideal.

\begin{corollary}\label{zero-ideal}
Let $A$ be a cg \bbc\   noetherian domain which is generated in degree
$1$ and has graded quotient ring $Q(A)=K[t,t^{-1};\sigma]$, for some $t\in A_1$ and
 $\sigma\in \mathrm{Aut}(K)$.
  Then there exists  $n_0\in \mb{N}$ such that:
\begin{enumerate}
\item
 There is a sequence of points $\{\eta_n = \Phi_n(\eta_{n+1})=\Psi_n(\eta_{n+1})\in \XX_n\}$
corresponding to the $K$-point module $Q(A)_{\geq 0}$  such that $k(\eta_n)\cong K$
for all $n\geq n_0$.
\item If $\YY_n$ denotes the closure of $\eta_n$ in $\XX_n$ then $\YY_n$ is a projective
integral scheme such that the morphisms
$\phi_n=\Phi_n{}\vert_{\YY_n}$  and $\psi_n=\Psi_n\vert_{\YY_n}$ are
birational (and hence surjective) maps for all $n\geq n_0$.
\end{enumerate}
\end{corollary}

\begin{proof}  (1)  Apply
Proposition~\ref{ideals}  to the ideal $P=0$. By Step 1 of the proof, the sequence of
points that this constructs does indeed correspond to $Q(A)_{\geq 0}$.

2) By
Proposition~\ref{eventually iso},   $\Phi_n^{-1}$ and
$\Psi_n^{-1}$ will be defined locally  at  $\eta_{n}$ for all
 $n\gg 0$ and
so $\phi_n$ and $\psi_n$ will be birational for such  $n $.
The rest follows from part (1).
\end{proof}

\begin{definition}\label{y-defn}
Keep the hypotheses and notation of Corollary~\ref{zero-ideal}. The
scheme $\YY_n$ will be called  the \emph{\rel
subscheme}\label{relevant-defn} of $\XX_n$ and the data
$(\YY_n,\phi_n,\psi_n)$ will be called the \emph{relevant point
scheme data} for $A$. In Proposition~\ref{component} we will show
that $\YY_n$ is actually an irreducible component of $\XX_n$, after
which we will
 call  $\YY_n$   the \emph{relevant component of $\XX_n$}.
We will usually write $\theta_n:\spec K\to \YY_n$ for the morphism
defined by $K\cong k(\eta_n)$. In the notation of \eqref{phi-mn},
the restrictions of $\Phi_{m,n}$ and $\Psi_{m,n}$ to $\YY_m$ will be
written $\phi_{m,n}$, respectively $\psi_{m,n}$.
\end{definition}

We want to  briefly explain why the schemes $\YY_n$ are important.
So let $A$ be  a noetherian cg domain that is generated in degree
one and \bbc, with $Q(A) = K[t,t^{-1}; \sigma]$. If $A$ is
\emph{strongly} noetherian then, as we explained in
Section~\ref{subsection-proofs},   $A$ is isomorphic  in high degree
  to a twisted homogeneous coordinate ring $B = B(Y,
\mc{L}, \sigma)$ where   the scheme $Y$ equals $\XX_n$ for any $n\gg
0$.

Now suppose that $A$ is merely noetherian.
The aim of this paper is, of course,  to prove that $A$ is a
\nsr\   $R(Y, Z_{\mc{I}},  \mc{L}, \sigma)$, again up to a finite
dimensional vector space, for some scheme $Y.$
Since $K=k(Y)$, the scheme $Y$ will necessarily be integral.
  Unfortunately,  $Y$ will typically not equal any $\XX_n$  or
  even be birational to some $\XX_n$
simply because the $\XX_n$ can have the wrong dimension or have
multiple components  (see Example~\ref{big-xn} and
Remark~\ref{big-xn-remark}). On the other hand   the relevant
subscheme $\YY_n$ is an irreducible component of $\XX_n$
  for all $n \gg  0$ (see Proposition~\ref{component}, below)
  and,  fortunately, it   also contains all the information we need to
  reconstruct $A$. Indeed,  the eventual scheme $Y$ will be obtained
  by blowing down $\YY_n$  for any large $n$.

We will need to  consider how the schemes $\XX_n$ and $\YY_n$  are
affected by the passage to Veronese rings.  Given  a  noetherian cg
domain $A$ generated in degree $1$ and  $p \geq 1$, the {\it
$p^{\mathrm{th}}$-Veronese ring}\label{veronese-defn} $A^{(p)} =
\bigoplus_{n = 0}^\infty A_{np}$ of  $A$ is again a cg noetherian
domain that is generated in degree one (see
\cite[Proposition~5.10]{AZ1}) and  we  write
$(\XX^{(p)}_n,\Phi^{(p)}_n,\Psi^{(p)}_n)$ for its point scheme data.
 If $A$ is also \bbc\, with $Q(A) \cong K[t,t^{-1};
\sigma]$ and  \rel subschemes  $\YY_n \subseteq \XX_n$, then
$A^{(p)}$ will again be \bbc\ (its graded quotient ring is
$K[t^p,t^{- p}; \sigma^p]$) and we let
$(\YY^{(p)},\phi_n^{(p)},\psi_n^{(p)})$ denote the relevant point
scheme data.

\begin{lemma}   \label{veronese}
Let $A$ be a cg \bbc\ noetherian domain that is generated in degree
one and pick $p\geq 1$.
\begin{enumerate}
\item
For all $n \geq 0$ there  exist closed embeddings $\rho_n: \XX_{np}
\to \XX_n^{(p)}$ satisfying $\Phi^{(p)}_{n-1} \rho_n = \rho_{n-1}
\Phi_{np,(n-1)p} $ and $\Psi^{(p)}_{n-1} \rho_n = \rho_{n-1}
\Psi_{np,(n-1)p}$, in the notation of \eqref{phi-mn}.
\item The  morphism
$\rho_n$ restricts to an isomorphism $\YY_{np} \xrightarrow{\sim}
\YY^{(p)}_n$ for all $n \gg 0$.
\end{enumerate}
\end{lemma}

\begin{proof}
(1) Let $C$ be  any commutative $k$-algebra. If $M =
\bigoplus_{i=0}^{np} M_i$ is a truncated  point module of length $np
+1$ over $A_C=A \otimes C$ then   $M^{(p)} = \bigoplus_{i  = 0}^n
M_{ip}$ is a truncated point module of length $n + 1$ over
$A^{(p)}_C$.  Since the map $M \mapsto M^{(p)}$ is compatible with
change of rings, we obtain a morphism of the representing schemes
$\rho_n: \XX_{np} \to \XX^{(p)}_n$.  The given decompositions of
$\Phi^{(p)}_{n-1}\rho_n$ and $ \Psi^{(p)}_{n-1} \rho_n$ then follow
  from Lemma~\ref{right left duality}(1).

It remains to show that  $\rho_n$ is a closed embedding. Write $A =
T(V)/J$, for $V=A_1$; thus $A^{(p)}\cong T(V^{\otimes
p})/\bigoplus_{n = 0}^{\infty} J_{np}$. Under these presentations,
$\XX_{np}$ is the zero set $\mc{V}(J_{np})\subseteq
\mb{P}(V^*)^{\times np}$  while $\XX^{(p)}_n =\mc{V}(J_n^{(p)})
\subseteq \mb{P}(V^{*\otimes p})^{\times n}$.
The map $\rho_n$ is just the
restriction to $\XX_{np} $ of the Segre embedding
$\mb{P}(V^*)^{\times np} = \left(\mb{P}(V^*)^{\times
p}\right)^{\times n} \longrightarrow \mb{P}(V^{*\otimes p})^{\times
n}$. As such, $\rho_n$ is a closed embedding.

(2) If $Q(A) = K[t,t^{-1}, \sigma]$  then $Q(A^{(p)})= K[t^p,t^{-p},
\sigma^p]$. The $K$-point module $Q(A)_{\geq 0}$ corresponds to a
sequence of morphisms $\theta_n: \spec K \to \XX_n$ for $n \geq 0$
and, similarly, $Q(A^{(p)})_{\geq 0}$ is a $K$-point module for
$A^{(p)}$
  corresponding to morphisms $\theta'_n: \spec K \to \XX^{(p)}_n$.
    Since $Q(A)_{\geq 0}^{(p)} = Q(A^{(p)})_{\geq 0}$,
 it follows from part (1)  that $\rho_n
\theta_{np} = \theta'_n$ for each $n$.  Since $\YY_{np}$,
respectively $\YY^{(p)}_n$, is the closure of the image of
$\theta_{np}$, respectively $\theta'_n$, it follows that
$\rho_n(\YY_{np}) = \YY^{(p)}_n$ for each $n$.  As $\rho_n$ is also
a closed embedding, it must restrict to an isomorphism $\YY_{np}
\xrightarrow{\sim} \YY^{(p)}_n$.
\end{proof}

To conclude the section, we will show that the relevant subschemes
of  a \bbc\ ring are irreducible components of the truncated point schemes.
This result will not be used elsewhere in the paper.

\begin{proposition}\label{component}
Let $A$ be a \bbc\ noetherian cg domain, with point scheme data
$(\XX_n, \Phi_n, \Psi_n)$ and relevant subschemes $\YY_n \subseteq
\XX_n$.
  Then $\YY_n$ is an irreducible component of $\XX_n$ for all  $n \gg 0$.
\end{proposition}

\begin{proof} Keep the notation from Corollary~\ref{zero-ideal} and
Definition~\ref{y-defn}; in particular, $\eta_n$ is  the generic
point of $\YY_n$ and  $\Phi_n^{-1}$ and $\Psi_n^{-1}$ define  local
isomorphisms   at
   $\eta_n$ for all $n \geq n_0$.
For each such  $n$, let $S_n$ be the  set of irreducible components
of $\XX_n$ which contain $\YY_n$.

  Suppose first  that $\dim \Phi_n(Z) < \dim Z$ for some $Z \in S_{n+1}$ with $n\geq n_0$.
  Then for every closed point $z \in \Phi_n(Z)$, and in particular for
  $z\in \Phi_n(\YY_{n+1}) = \YY_n$,
the fibre $\Phi^{-1}_n(z)$ has dimension $\geq 1$. This contradicts
the fact that  $\Phi_n^{-1}$ defines a local isomorphism at $\eta_n$
and proves  that $\dim \Phi_n(Z) = \dim Z$ for all $Z\in S_{n+1}$.
Similarly, $\dim \Psi_n(Z) = \dim Z$ for all $Z \in S_{n+1}$.
  For  $n\geq n_0$ let $d_n$ be the maximum dimension of the components in
$S_n$. Since $\Phi_n(Z) \subseteq Z'$ for some $Z' \in S_n$, we have
$d_{n+1} \leq d_n$ for $n \geq n_0$.  Thus $d_n  $ is constant for
$n \gg 0$, say $d_n=d$  for $n \geq n_1 \geq n_0$.

For each $n \geq n_1$ let $T_n \subseteq S_n$ consist of those
components of dimension $d$.  Let $\Omega_n$ be the set of generic
points of the components in $T_n$.  Then $\Phi_n(\Omega_{n+1})
\subseteq \Omega_n$ and similarly  $\Psi_n(\Omega_{n+1}) \subseteq
\Omega_n$ for all $n \geq n_1$.   By
 Corollary~\ref{finite sets}, we may pick $n_2 \geq n_1$ such that
 for $n\geq n_2$  both   $\Phi_n$ and $\Psi_n$ give bijections from
$\Omega_{n+1} $ to $\Omega_n$ and both   $\Phi_n^{-1}$ and
$\Psi_n^{-1}$ define local isomorphisms at each point of $\Omega_n$.
  In particular, $|T_n|$ is a constant value $m$ for $n\geq n_2$.  Now
$\Phi_n \Psi_{n+1} = \Psi_n \Phi_{n+1}$  by Lemma~\ref{right left
duality}(3). Thus, for each $n\geq n_2$, we can label $\Omega_n =
\{\mu_{1,n}, \dots, \mu_{m,n} \}$   in such a way that
$\Phi_n(\mu_{i,{n+1}}) = \mu_{i,n}$ and $\Psi_n(\mu_{i, {n+1}}) =
\mu_{\tau(i), n}$, where $\tau$ is a permutation in the symmetric
group on $\{1, 2, \dots, m \}$, independent of the choice of  $n$.

Set $p =|\tau|$ and write  $A' = A^{(p)}$, with point scheme data
$(\XX'_n, \Phi'_n, \Psi'_n)$ and relevant subschemes  $\YY'_n
\subseteq \XX'_n$.  Let $\rho_n: \XX_{np} \to \XX_n'$
be the closed embeddings of Lemma~\ref{veronese}; thus
$\YY'_n=\rho_n(\YY_{np})$ by that lemma. If $\eta'_n$ denotes the
generic point of each $\YY'_n$, then $\Psi'_n(\eta'_{n+1}) =
\Phi'_n(\eta'_{n+1}) = \eta'_n$ for all $n \gg 0$. On the other
hand, if   $\mu'_n = \rho_n(\mu_{1,np})$,
 then   Lemma~\ref{veronese} shows that
$\Psi'_n(\mu'_{n+1}) = \Phi'_n(\mu'_{n+1}) = \mu'_n$ for all $n \gg
0$.  Applying Proposition~\ref{ideals} to the ring $A'$ shows that the points
$\{\eta'_n\}$ correspond to the zero  ideal while  the
$\{\mu'_n\}$ correspond to some (possibly different) completely prime ideal $P$ of
$A'$.

We need to show that $P=0$. Note that $\eta'_n$ is in the Zariski
closure $Z'_n = \overline{\mu'_n}$ for all $n \geq n_2$.  By
Proposition~\ref{eventually iso}, pick $n_3\geq n_2$ such that
$(\Phi'_n)^{-1}$ defines  a local  isomorphism   both at  $\eta'_n$
and at  $\mu'_n$ for all $n \geq n_3$.  Let $R$ be the local ring $R =
\mc{O}_{Z'_{n_3}, \eta'_{n_3}}$   with field of fractions $L\cong
k(\mu'_{n_3})$ and residue field $K \cong k(\eta'_{n_3})$. As usual,
write $f^\star: \spec L\to \spec R$ and $g^\star:  \spec K\to \spec
R$ for the morphisms induced from the corresponding maps $f: R \to
L$ and $g: R \to K$. Then we can define morphisms $\alpha_{n_3}:
\spec K \to \XX'_n$ with image $\eta'_n$,  as well as $\beta_{n_3}:
\spec L \to \XX'_n$ with image $\mu'_n$, and $\gamma_{n_3}: \spec R
\to \XX'_n$ in such a way that $\alpha_{n_3} = \gamma_{n_3} g^\star$
and $\beta_{n_3} = \gamma_{n_3}f^\star$. Since
$(\Phi_{n,n_3}')^{-1}$ is defined as a morphism on
$\mathrm{Im}(\gamma_{n_3})$, we can inductively  define maps
$\gamma_{n+1}=(\Phi_n')^{-1}\gamma_n : \spec R\to \XX_{n+1}'$ such
that $\Phi_n'\gamma_{n+1} = \gamma_n$ for all $n\geq n_3$. The
morphisms $\alpha_n=\gamma_ng^\star$ and $\beta_n = \gamma_n
f^\star$ then satisfy $\Phi_n'\alpha_{n+1} = \alpha_n$ and $\Phi_n'
\beta_{n+1} = \beta_n$ for all $n\geq n_3$.

The sequences of maps $\{\alpha_n\}, \{\beta_n\}, \{\gamma_n\}$ for
$n\geq n_3$ correspond to point modules $M, N,$ and $Q$ over the
rings $K$, $L$, and $R$, respectively, such that $K {_g \otimes}_R Q
\cong M$ and $L {_f \otimes}_R Q \cong N$. It follows that $\ann_A N
= \ann_A Q \subseteq  \ann_A M$.  By the correspondence of
Proposition~\ref{ideals}, this means that $P \subseteq (0)$, whence
$P = 0$. This forces  $\mu'_n = \eta'_n$ for all $n\gg  0$ which,
since each $\rho_n$ is a closed embedding, implies that $\mu_{1,np}
= \eta_n$ for each such $n$. Thus $\YY_n$ is indeed a component of
$\XX_n$.
\end{proof}


\section{Contracted points}\label{contracted-section}

Let $A$ be a noetherian cg domain that is \bbc.  If the morphisms
$\phi_n, \psi_n :\YY_{n+1}\to \YY_{n}$ are isomorphisms for $n \gg
0$ then it is not hard to understand the structure of $A$ and,
indeed, to prove that $A$ is a twisted homogeneous coordinate ring
in large degree (see \cite[Theorem 4.4, Corollary 4.6]{RZ}).  To
understand the general case, it will be useful to examine the the
points where the maps $\phi_n$ and $\psi_n$ are \emph{not}
isomorphisms, with particular emphasis on the points where they have
infinite fibres, and this is what we   do over the next  four sections.
Eventually, in Theorem~\ref{blowing down}, this information will be used
to construct a surjective image of $\YY_n$, for some $n \gg 0$, on
which the birational map $\sigma=\psi_n\phi_{n}^{-1}$ induces an
isomorphism.

To    achieve  this aim we will
 have to restrict our attention to   surfaces and so,
  from now on,   we will assume:

\begin{assumptions}\label{assumption4}
$A$ is a cg noetherian domain,  generated in degree $1$, which
is \bbc\  with graded quotient ring $Q = K[t,t^{-1}; \sigma]$,
where $K$ is a field with  $\trdeg K/k = 2$.
\end{assumptions}

The points that interest us are described as follows.

\begin{notation}
\label{smn-notation2} Given an algebra  $A$ as in
\eqref{assumption4}, let     $(\YY_n, \phi_n,  \psi_n)$ denote   the
\rel point scheme data from
 Definition~\ref{y-defn}.  By  Corollary~\ref{zero-ideal}, fix
 an integer $n_0$  such that   $k(\YY_n)
\cong K$ and the morphisms $\phi_n, \psi_n$ are birational, for all
$n \geq n_0$. For $n \geq n_0$, let $\Sminus_n$ denote the set of
closed points $y \in \YY_n$ such that the fibre $\phi_n^{-1}(y)$ has
dimension $\geq 1$. Similarly, let $\Splus_n$ be the set of points
$y \in \YY_n$ such that the fibre $\psi_n^{-1}(y)$ has dimension
$\geq 1$.  The points in $\Sminus_n$ and $\Splus_n$  will be called
\emph{contracted points}.  As $\dim \YY_n = 2$
 it  follows, for example   from \cite[Exercise~II.3.22(d)]{Ha}, that
$\Sminus_n, \Splus_n$ are finite sets.\label{contracted-defn}
\end{notation}

\begin{example} \label{simplest-eg}
In order to understand these concepts, consider the example $A =
R(\mb{P}^2, Z, \mc{O}(1), \sigma)$ from Section~\ref{first-eg},
where $Z$ is the single reduced point $c = (1:1:1)$ and $p,q$ are
algebraically independent over the prime subfield of $k$. We will
not prove the assertions we make about $A$ except to note that, at
the level of sets, our description of the closed points of $Y_n$
follows from  \cite[Theorem~6.6]{Ro1}. To see this, note that, by
Corollary~\ref{zero-ideal}(2), each truncated point module $y_n$
corresponding to a point in $Y_n$ is a homomorphic image of a point
module. Those point modules are then described by
\cite[Theorem~6.6]{Ro1}. The proof that we have the correct
description of $Y_n$ as a
 scheme is more complicated and is omitted.

 Write  $r_m=\sigma^{1-m}(c) = (1: p^{1-m}:q^{1-m})$ for $m \in \mb{Z}$.
If $n$ is sufficiently large   then $\YY_n$ is the
scheme obtained as the blowup    $\gamma_n :\YY_n \to \mb{P}^2$ of $\mb{P}^2$ at
$Z_n=\{r_1,\dots,r_{n}\} $.  Write the resulting exceptional curves
  as   $C_{i}^{(n)} = \gamma_n^{-1}(r_i)$ for $1\leq i\leq n$  but, by a slight abuse of notation,
   set $r_m=\gamma^{-1}(r_m)$ for   $m\not\in [1,n]$.
The  maps $\phi_n$ and $ \psi_n$ are easily described.  The morphism
$\phi_n: \YY_{n+1} \to \YY_n$  is simply obtained by blowing up
$r_{n+1}$ to obtain the new exceptional curve $C^{(n+1)}_{n+1}$. On
the other hand, $\psi_n : \YY_{n+1}\to \YY_n$ is obtained by blowing
up $r_0=c$ and then ``shifting by $\sigma$''; thus
 $\psi_n^{-1}(r_0) = C_1^{(n)}$  while $\psi_n^{-1}(C^{(n)}_j) = C_{j+1}^{(n+1)}$
for $1\leq j\leq n$.

We remark that the conventions from the previous paragraph are the ones that are most natural
 for \naive\ blowups (see, for example,  \cite{KRS}).
 However, the notation that is most useful for this paper is  the opposite---thus  $\Sminus_n = \{ r_{n+1} \}$ and $\Splus_n = \{ r_{0} \}$.
\qed\end{example}

The following result gives the most basic properties of the
contracted points.

\begin{lemma}
\label{exceptional points} Let $A$ satisfy \eqref{assumption4} and
keep the notation of \eqref{smn-notation2}. Then there exists
$n_1\geq n_0$ such that:
\begin{enumerate}
\item For all $p\geq n_0$ one has $\psi_p(\Sminus_{p+1}) \subseteq \Sminus_p$
and $\phi_p(\Splus_{p+1}) \subseteq \Splus_p$.
\item For all $q > p \geq n_1$ the birational map
$\psi_{q,p}^{-1}$ is defined at each point $s \in \Sminus_p$, with
$\psi_{q,p}^{-1}(s)\in \Sminus_q$.  Similarly,
$\phi_{q,p}^{-1}(t)$ is defined and belongs to $\Splus_q$ for each
$t\in \Splus_p$.
 \item  For all $p\geq n_1$,      $\Sminus_p \cap \Splus_p = \emptyset$ and
 the cardinalities  $|\Sminus_p|$ and $|\Splus_p|$ are  constant.
  \end{enumerate}
\end{lemma}

\begin{proof}
  In parts (1) and (2),  it suffices by symmetry  to  prove the statements for $\Sminus_p$.

  (1) Recall
that  we have a closed immersion $\iota_p: \XX_p \to
\mathbb{P}^{\times p}$, and the morphisms $\Phi_p, \Psi_p :
\XX_{p+1} \to \XX_p$ are defined to be the restrictions of the
respective projection maps $\pi_{1,p},\pi_{2,p+1} :
\mathbb{P}^{\times(p+1)}\to \mathbb{P}^{\times p}$.  Now let $y\in
\Sminus_{p+1}$, say given by $y=(y_1,\dots, y_{p+1})\in
\mathbb{P}^{\times (p+1)}$.  By definition, $y$ has infinitely many
preimages under $\phi_{p+1}$, say $z_\alpha = (y_1,\dots,
y_{p+1},\alpha)\in \YY_{p+2}\subseteq \mathbb{P}^{\times (p+2)}$,
for some infinite set of $\alpha\in \mb{P}$.  The points $w_\alpha =
\psi_{p+1}(z_\alpha) = (y_2,\dots,y_{p+1},\alpha)$ are then distinct
elements of $\YY_{p+1}$, each of which satisfies $\phi_p(w_\alpha)=
(y_2, y_3, \dots, y_{p+1}) = \psi_p(y) \in \YY_p$.  Hence
$\psi_p(y)\in \Sminus_p$.

(2) The sets $\Sminus_p$ are finite for $p \geq n_0$ and
$\Psi_p(\Sminus_{p+1}) \subseteq \Sminus_p$ by part (1).  By
Corollary~\ref{finite sets}  there exists $n_1 \geq n_0$ such that,
if $p\geq n_1$ then $|\Sminus_p|$ is constant and $\Psi_p^{-1}$ is
locally an isomorphism on  $\Sminus_p$. For  such $p$, the
restriction $\psi_p^{-1} $ will still  be a local isomorphism at
$\Sminus_p$ and so, by induction, the birational map
$\psi_{q,p}^{-1}: \YY_p\dra \YY_{q}$ is defined at   $ \Sminus_p$
for all $q>p\geq n_1$.

(3)  Repeating the argument of part (2) for the maps $\Phi_n$, we
can choose a single $n_1$ so that the cardinalities $|\Sminus_p|$
and $|\Splus_p|$ remain constant for $p\geq n_1$.   Since
$\psi_p^{-1}$  is defined on $\Sminus_p$, the map $\psi_p$ obviously
 cannot have a fibre of  positive dimension
over $\Sminus_p$.  Thus $\Sminus_p \cap \Splus_p =
\emptyset$.
\end{proof}

By definition, the birational map $\phi_p^{-1}: \YY_p \dra
\YY_{p+1}$ is not defined at any point  $y\in \Sminus_p$.
Conversely, we would like to say that for  $p \gg 0$, the domain of
definition of $\phi_p^{-1}$ is exactly  $\YY_p \smallsetminus
\Sminus_p$. Although this holds for Example~\ref{simplest-eg},
 in general it is too strong an assertion.
 However,  something weaker is true; as we next show,  the points where
$\phi_p^{-1}$ is undefined are either contracted points themselves,
or else lie on a curve which   contracts at a later stage.

\begin{lemma}
\label{stein fact} Let $A$ satisfy \eqref{assumption4} and keep the
notation of \eqref{smn-notation2}.  Fix   $n_1 \geq
n_0$. Then there exists $n_2 \geq n_1$
with the following properties.
\begin{enumerate}
\item For some $m\geq n_2$, let $y \in \YY_m$  be a closed point such that
 the map $\phi_m^{-1}$ is not defined at
  $y  $. Then either $y \in \Sminus_m$, or else there exists $m>p\geq n_1$ and a curve
  $C\subset \YY_m$ such that $y \in C $ and   $\phi_{m,p}(C)   \in \Sminus_p$.
  \item   Suppose that $\psi_m^{-1}$ is not defined at a closed point
  $y \in \YY_m$ for some $m \geq n_2$.
Then either $y \in \Splus_m$, or else $y \in C \subset \YY_m$ for
some curve $C$ with $\psi_{m,p}(C)  \in \Splus_n$ for some $m > p
\geq n_1$.
\end{enumerate}
\end{lemma}
\begin{proof}
As usual, we prove only part (1); the proof of part (2) is
analogous.

The main idea of the proof is to study the morphisms $\phi_{m,n}$
through Stein factorizations, as described in
\cite[Corollary~III.11.5]{Ha}, and we start
with the relevant notation. Choose   integers $m > n \geq n_1$ and,
for ease of notation, set $q = n_1$. The \emph{Stein factorization}
$\phi_{m,n} = g_{m,n} \circ h_{m,n}$
 of the morphism $\phi_{m,n}: \YY_m \to \YY_n$  has the following properties:
  $h_{m,n}: \YY_m
\to Z_{m,n}$  is a projective morphism
 to a scheme defined by $Z_{m,n} =
\mathbf{Spec} (\phi_{m,n})_* \mc{O}_{\YY_m}$,  while the morphism
$g_{m,n}: Z_{m,n} \to \YY_n$ is finite. We remark that $h_{m,n}$
necessarily has connected fibres and, since $\phi_{m,n}$ is
birational, the maps $g_{m,n}$ and $h_{m,n}$ must also  be
birational.

By Corollary~\ref{zero-ideal}  $\phi_{m,q}: \YY_m \to \YY_q$ is a
birational surjective morphism for all $m \geq q$ and so
  we get an ascending chain of coherent sheaves
of $\mc{O}_{\YY_q}$-algebras,
\begin{equation}\label{Y-chain}
\mc{O}_{\YY_{q}} \subseteq (\phi_{q +1, q})_* \mc{O}_{\YY_{q +1}}
\subseteq \dots \subseteq (\phi_{m, q})_* \mc{O}_{\YY_m} \subseteq
\dots
\end{equation}
and associated to this a chain of finite morphisms
\begin{equation}
\label{finite maps}
\cdots \overset{e_m}{\longrightarrow} Z_{m,q}
\overset{e_{m-1}}{\longrightarrow} \cdots
 \overset{e_{q +1}}{\longrightarrow} Z_{q +1, q}
 \overset{e_{q}}{\longrightarrow} Z_{q} =  \YY_q
\end{equation}
where  $g_{m,q}= e_{q} \circ e_{q +1} \circ \dots \circ e_{m-1} $ is
the finite morphism half of the Stein factorization of $\phi_{m,q}$.
Choose a finite  open affine cover $\YY_q = \bigcup U_{\alpha}$ of $\YY_q$.
Then  looking at the sections of \eqref{Y-chain} over  any
$U_{\alpha}$ gives an ascending chain of finite ring extensions,
which must eventually stabilize by the finiteness of the integral
closure \cite[Theorem 4.14]{E}.   Thus \eqref{Y-chain} must also
stabilize   and so  $e_m$ must be an isomorphism for all $m\gg q$,
say for all $m\geq n_2$.

For $m > n_2$ we have a commutative square:
\[
\begin{CD}
\YY_{m+1}  @>\phi_m >> \YY_m \\
@V h_{m+1, q} VV @V h_{m,q} VV  \\
Z_{m+1,q} @>e_m >> Z_{m,q}
\end{CD}
\]
By assumption $e_m$ is an isomorphism and so the map $\theta =
h_{m,q} \phi_m = e_m h_{m+1, q} : \YY_{m+1}\to Z_{m,q}$ has
connected fibres, since $h_{m+1, q}$ does.  Let $y \in \YY_m$ be a
closed point where the  birational map $\phi_m^{-1}$ is not defined
and consider $z = h_{m,q}(y) \in Z_{m,q}$.   If the birational map
$\theta^{-1}$ were defined at $z$, then $\phi^{-1}_m= \theta^{-1}
h_{m,q}$ would be defined at $y$, contradicting our assumption.  So
$\theta^{-1}$ is not defined at $z$.

Suppose now that the set-theoretic fibre $\theta^{-1}(z)$ is a
single point $z'$.  By construction,  $\theta= e_m h_{m+1, q}$ is
projective with $\theta_* \mc{O}_{\YY_{m+1}} =\mc{O}_{Z_{m,q}}$.
Thus Lemma~\ref{point fiber}(2) can be applied
to show that the birational map $\theta^{-1}$ must have been defined
at $z$ after all, a contradiction.  As $\theta$ has connected
fibres, we conclude that the fibre  $\theta^{-1}(z)$
must contain an  irreducible curve $C$ such that
$y \in \phi_m(C)$.  Since
\[
\phi_{m+1, q} = \phi_{m,q} \phi_m = g_{m, q} h_{m,q} \phi_m =
g_{m,q} \theta,
\]
we see that $\phi_{m+1, q}(C) = g_{m,q}(z)$ is a closed  point in
$\YY_{q}$. Now let $p$ be the largest integer  with $m \geq p \geq
q$ for which $\phi_{m+1, p}(C)$ is a closed point $w \in \YY_p$.
Then $w\in \Sminus_p$ and $\phi_{m,p}(y) = w$. If $p = m$, then  $w
= y \in \Sminus_m$; otherwise   $m> p \geq q$ and $\phi_m(C) = D$
is a curve, containing $y$, for which   $\phi_{m,p}(D) = w \in
\Sminus_p$.
\end{proof}


\section{The  stable scheme and its  birational map}\label{stable-section}

Throughout this section  $A$ will satisfy
Assumptions~\ref{assumption4} and we keep the notation from
Notation~\ref{smn-notation2}.  In particular, $Q=Q(A)  = K[t,t^{-1};
\sigma]$ for some $\sigma\in \aut_kK$.\label{sigma-notation}   We are
interested in the birational map  that  $\sigma$ induces on $\YY_n$
for $n\gg 0$ and,
 for simplicity, we continue to write this  as $\sigma: \YY_n \dra \YY_n$.
As we will see in this section, the birational map $\sigma$
has   a number of special properties.

These properties are fairly technical but it is at least possible to indicate where
 the results are heading (see
Proposition~\ref{consequences of geometric} and Theorem~\ref{blowing
down} for the formal statements).  The basic idea is that we are
aiming to prove that $\YY_n$ is the blowup $\pi_n:\YY_n\to\wtY_n$ of
a scheme $\wtY_n$ such that (i) $\sigma$ is induced from an
automorphism $\wt{\sigma}$ of $\wtY_n$ and (ii)  the points that are
blown up all lie on Zariski dense orbits of $\wt{\sigma}$. So,
consider the set $V$ of closed points $y\in \YY_n$ at which some
power of $\sigma$ is undefined or not a local isomorphism. Under a
minor extra hypothesis $V$ breaks into a finite union of sets
$V=\bigcup V_i$ that are ``trying to be'' infinite orbits of
$\sigma$. Formally,  within each $V_i$   there are two dense   half
orbits $V_i^+=\{ y_{j,i}=\sigma^j(t_i) : j \geq 0 \}$ and $V_i^-= \{
x_{j,i}=\sigma^{-j}(s_i) : j \geq 0 \}$.  The remaining points in
$V_i$,  which essentially  correspond to points on the exceptional
divisor of $\pi_n$,  are contained in a finite union  of irreducible
curves (and possibly some sporadic closed points)  $C_\ell$  all of
which are mapped  to both $s_i$ and $t_i$ by suitable powers of
$\sigma$; that is, there exist unique integers $a=a(\ell)>0$ and
$b=b(\ell)<0$ such that $\sigma^a(C_{\ell})=t_i$ and
$\sigma^b(C_{\ell})=s_i$.

The results mentioned in the previous paragraph  only hold for sufficiently large values of $n$
and it is useful to  more formally quantify
 ``sufficiently large''   as follows:

\begin{definition} \label{stable-defn}
Define $n_0$ by Corollary~\ref{zero-ideal}.
  Then fix some $n_1 \geq n_0$  such that  the
conclusions of Lemma~\ref{exceptional points} hold.  Then take this
$n_1$ in Lemma~\ref{stein fact} and choose some $n_2 \geq n_1$
satisfying the conclusions of that lemma. Finally, we call any $n >
n_2$ a \emph{stable value of $n$.}  \label{stable-defn2}
\end{definition}

\begin{notation}\label{S-notation}
Fix once and for all a stable value of $n$, and set $\YY = \YY_n$.
Recall from Corollary~\ref{zero-ideal}  that for each $m \geq n$,
the map $\theta_m: \spec K \to \XX_m$ corresponding to the truncated
$K$-point $\bigoplus_{i = 0}^m Q_i$ induces an isomorphism
$k(\eta_m) \to K$, where $\eta_m$ is the generic point of the
relevant component $\YY_m$.   We   identify  $k(\YY) = K$ through
the map $\theta_n$ and,
 as before,  write $\sigma : \YY\dra \YY$ for the birational map induced by the automorphism
 $\sigma\in \mathrm{Aut}(K)$.  Set
$\Splus = \Splus_n$ and $\Sminus = \Sminus_n$ for the  contracted
points, as defined in Notation~\ref{smn-notation2}.
\end{notation}

To simplify notation,  we often write
$\phi_m=\phi$, for any $m \geq n_0$,
 where the subscript will be clear from the domain or
range of $\phi$. Similar comments apply to $\psi$ and so, for
example, for any $m \geq n_0$ the birational map $\psi_{m+s,m}^{-1}$
is now written $\psi^{-s}$. The next result, which shows  how to
interpret $\sigma^m : \YY \dra \YY$
 in terms of the morphisms $\phi$ and $ \psi$,
  will be used frequently and  usually without further comment.

\begin{lemma}\label{sigma-lemma}  {\rm (1)}  As  birational maps $\YY_p
\dra \YY_{p-m-\ell}$ one has $\phi^{\ell} \psi^m = \psi^m
\phi^{\ell}$
 for any $\ell, m \in \mb{Z}$ for which all the maps
involved are birational (that is, as long as $p$,  $p-m$, $p-\ell$
and  $p-m-\ell$ are $\geq n_0$).

{\rm (2)} For any $m \geq 0$, we have $\sigma^m = \psi^m \phi^{-m}$
and $\sigma^{-m} = \phi^m \psi^{-m}$ as birational maps $\YY \dra
\YY.$
\end{lemma}
\begin{proof}
(1) When $\ell,m \geq 0$ this follows  from Lemma~\ref{right left
duality}(3).
 When  $\ell$ or $m$ is negative, multiply through  by the appropriate (birational) inverse map.

(2) Recall from step $1$ of Proposition~\ref{ideals} that $\theta_n
= \Phi_n \theta_{n+1}$ and $\theta_n \sigma^\star  = \Psi_n
\theta_{n+1}$, where $\theta_m : \spec K \to \XX_m$ is defined by
Notation~\ref{S-notation} and $\sigma^\star: \spec K \to \spec K$ is
the map of schemes  associated to $\sigma\in \mathrm{Aut}(K)$.
Restricting to $\YY$ we find that $\psi \phi^{-1} \theta_n = \psi
\theta_{n+1} = \theta_n \sigma^{\star}$,  and it follows that $\psi \phi^{-1}$ is precisely the
birational map $\sigma: \YY \dra \YY$ induced by $\sigma \in
\aut(K)$.  That $\sigma^{-1} = \phi \psi^{-1}$ follows immediately,
and the results for all $m \geq 1$ follow by induction using part
(1).
\end{proof}

The rest of this section is devoted to proving successively stronger
properties of $\Splus$ and $\Sminus$, and of certain special
extremal subsets thereof, for our fixed stable value of $n$. The
properties of $(\Splus,\sigma)$ and $(\Sminus, \sigma^{-1})$ are
analogous and so  we  will  first consider  $(\Splus,\sigma)$,
referring the reader to Remark~\ref{first duality} for
$(\Sminus,\sigma^{-1})$. We begin by  studying   points $\mu \in
\YY$ which are fixed by some power $\sigma^n$ and show that they do
not interfere with $\Splus$.

\begin{lemma}
\label{fixed point and curve} Let $\mu \in \YY$  be
a point such that $\sigma^p$ is defined at $\mu$ and $\sigma^p(\mu)
= \mu$ for some $p > 0$.
\begin{enumerate}
\item If $\mu$ is  a closed point of $\YY$ then  $\mu \not \in \Splus$.
\item If $\mu$ is the generic point of an irreducible curve $C$, then
$C \cap \Splus = \emptyset$.
 \end{enumerate}
\end{lemma}

\begin{proof}
The proof is by contradiction, so we assume that there is a closed
point $z \in \Splus$ where either $z = \mu$ in case (1), or $z \in
C$ in case (2). By Lemma~\ref{exceptional points}(2), $\phi^{-s}$ is
defined at $z$ for any $s \geq 1$ and then $\phi^{-s}(z) \in
\Splus_{s+n}$.  In either case the birational map $\phi^{-s}$ is
also defined at $\mu_n =\mu$   for $s \geq 1$, so we can  define
$\mu_r = \phi^{n-r}(\mu)=\phi^{-1}(\mu_{r-1}) \in \YY_r$ for all $r
n$.

We claim that $\psi^p(\mu_{r+p}) = \mu_r$ for all $r \geq n$. When
$r = n$ this is just the observation that $\psi^p(\mu_{n+p}) =
\psi^p \phi^{-p}(\mu) = \sigma^p(\mu) = \mu$.  We then prove it for
all $r \geq n$ by induction: assuming the result is true for some $r
\geq n$,  then
\[
\psi^p(\mu_{r +1+p}) = \psi^p \phi^{-1}(\mu_{r+p}) = \phi^{-1}
\psi^p(\mu_{r+p}) = \phi^{-1}(\mu_r) = \mu_{r+1}
\]
and so the result holds for $r +1$.

Consider the Veronese ring $A' = A^{(p)}$ of $A$ with  point scheme
data
  $(\XX'_m, \Phi'_m, \Psi'_m)$  and relevant point scheme data $(\YY'_m, \phi'_m,
\psi'_m)$.  By Lemma~\ref{veronese}(2), there is an isomorphism
$\rho_p: \YY_{pm} \to \YY'_m$ for each  $m \gg 0$ and we write
$\mu'_m = \rho_m(\mu_{mp})$ for such $m$. It follows from the previous
paragraph and  Lemma~\ref{veronese}(1) that $\psi'(\mu'_{m+1}) =
\mu'_m = \phi'(\mu'_{m+1})$ for all such $m$.  By
Proposition~\ref{ideals}, there is  then a  completely prime ideal
$P$ of $A'$ corresponding to this sequence of elements such that
$\widehat{A} = A'/P$ is \bbc,   with function ring $L\cong
k(\mu'_m)$ for $m \gg 0$.

Now consider the point scheme data
$(\widehat{\XX}_m,\widehat{\Phi}_m,\widehat{\Psi}_m)$ for
$\widehat{A}$.  By Lemma~\ref{inverse image stable},
$\widehat{\XX}_m$ is a closed subscheme of $\XX'_m$ for each $m$,
and by construction $\mu'_m \in \wh{\XX}_m$ for $m \gg 0$. By
assumption $L$ is a field of transcendence degree $\leq 1$ over $k$,
so \cite[Theorem~0.1(ii)]{AS1} implies that
the Gelfand-Kirillov dimension of  $\widehat{A}$ satisfies $\GKdim \widehat{A} \leq 2$.
By  \cite[Theorem~4.24]{ASZ} $\widehat{A}$ is strongly noetherian and
hence, by \cite[Theorem~E4.4]{AZ2} and Lemma~\ref{right left
duality}(2),   the maps $\widehat{\Phi}_m$ and $\widehat{\Psi}_m$ are
isomorphisms for all $m \gg 0$. Lemma~\ref{inverse image stable}
therefore implies that the fibre $(\Psi')^{-1}(y)$ is a singleton
for all $y \in \widehat{\XX}_m$ and all $m \gg 0$.

Recall that $z=\mu$ in case (1), respectively $z\in C$ in case (2),
and
 set $z_{pm} = \phi^{-pm+n}(z) \in \YY_{pm}$ for all $m \gg 0$,
so $z_{pm}$ is in the closure of $\mu_{pm}$. Then $z'_m =
\rho_m(z_{pm})$ is in the closure of $\mu'_m$, whence $z'_m\in
\widehat{\XX}_m$.   By the argument above,
 the fibre $(\Psi')^{-1}(z'_m)$ is a singleton for all $m \gg
0$. On the other hand, since $z_{pm} \in \Splus_{pm}$, we can find
an irreducible curve $D_{pm+p} \subset \YY_{pm+m}$ such that
$\psi^p(D_{pm+p}) = z_{pm}$.  Then Lemma~\ref{veronese} implies that
$\rho_{m+1}(D_{pm+p})$ is a curve in $\XX'_{m+1}$ that contracts to
$z'_m$
 under $\Psi'$; a  contradiction.
\end{proof}

Next, we reinterpret Lemmas~\ref{exceptional points} and \ref{stein
fact} to see how the contracted points $\Splus$ are related to the
points where $\sigma^{-1}$ is undefined.

\begin{lemma}
\label{Splus properties} {\rm (1)}  Let $y \in \Splus \subseteq
\YY.$ Then $\sigma^{-1}$ is not defined at $y$, but $\sigma^d$ is
defined at $y$ for all $d \geq 0$.
\begin{enumerate}
\item[(2)]  Conversely, let $y \in \YY$ be a closed point where
either $\psi^{-1}$ or $\sigma^{-1}$ is not
defined. Then there exists $e\geq 0$ such that $\sigma^e(y)$ is
defined and belongs to $\Splus$.  \end{enumerate}
\end{lemma}

\begin{proof}
(1) Since $y \in \Splus$  there exists, by definition,  an
irreducible curve $D \subset \YY_{n+1}$ such that $\psi(D) = y$.
Suppose that $\phi(D) = z$ is also a closed point; thus $z \in
\Sminus$.  By Lemma~\ref{exceptional points}(1),
$$w = \phi(y) =\phi\psi(D)=\psi\phi(D)= \psi(z) \in \Splus_{n-1} \cap \Sminus_{n-1},$$
 contradicting Lemma~\ref{exceptional points}(3).
 This contradiction implies that
  $\phi(D) = C$ is a
curve, and $w = \psi(C)=\phi\psi(D)=\phi(y)$ is a closed point in $
\Splus_{n-1}$. By Lemma~\ref{exceptional points}(2), $\phi^{-1}$ is
defined at $w$ and so $\sigma = \phi^{-1} \psi$ is defined at every
point of $C$. Since $\sigma$ contracts $C$ to $y$, $\sigma^{-1}$ is
not defined at $y$. On the other hand, by Lemma~\ref{exceptional
points}(2) again, $\phi^{-d}(y)$ is defined for all $d \geq 0$ and
so $\sigma^d(y) = \psi^d \phi^{-d}(y)$ is defined for all $d \geq
0$.

(2) As $\sigma^{-1} = \phi \psi^{-1}$, if $\sigma^{-1}$ is not
defined at $y$, then the birational map $\psi^{-1}$ cannot be
defined at $y$.  Thus we need only consider the case that
$\psi^{-1}$ is not defined at $y$.  Since $n$ is a stable value,
Lemma~\ref{stein fact}(2) implies that either $y \in \Splus$ or $y
\in C \subset \YY$ for some irreducible curve $C$ with the property
that $z=\psi_{n,p}(C) \in \Splus_p$ for some $n > p \geq n_1$.  In
the former case, we can take $e = 0$; in the latter case,
Lemma~\ref{exceptional points} implies that $\phi_{n,p}^{-1}(z)$ is
defined and in $\Splus$, so taking $e = n-p$ we have $\sigma^e(y) =
\psi^e \phi^{-e}(y) = \phi^{-e} \psi^e(y) = \phi^{-e}(z) \in
\Splus$.
\end{proof}

It can happen that some power $\sigma^m$ of $\sigma$ sends a point  $y \in \Splus$
to another point in $\Splus$ but, as we show next, this cannot happen for $m\gg 0$.
  The transition between these two possibilities will be fundamental and  so we define
  the
\emph{extremal elements} of $\Splus$ to be
\begin{equation}\label{extremal-defn}
\Splusex = \{ y\in \Splus :  \forall\ d\geq 1, \ \
\sigma^d(y)\not\in \Splus \}.
\end{equation}   This set
 will turn out to be the set of distinguished elements
  $\Splusex=\{t_i\}$  mentioned in the
introduction to this  section.  The dual set $\Sminusex=\{s_j\}$
 will be defined after Remark~\ref{first duality}.

\begin{lemma}
\label{propertyD}
If $y \in \Splus$ then $\sigma^e(y) \in
\Splusex$ for a unique integer $e$, necessarily with  $e \geq 0$.
\end{lemma}

\begin{proof}
Let $y_0=y \in \Splus$.  By Lemma~\ref{Splus properties}(1),
$\sigma^d(y_0)$ is defined for all $d \geq 0$.  If $y_0 \not \in
\Splusex$, then there exists $d_0 > 0$ such that $\sigma^{d_0}(y_0)
= y_1 \in \Splus$.  Similarly, $\sigma^d(y_1)$ is defined for all $d
\geq 0$ and so either $y_1 \in \Splusex$ or there is $d_1 > 0$ with
$y_2 = \sigma^{d_1}(y_1) = y_2 \in \Splus$.
If this process continues forever, we will produce an infinite
sequence of points $y_0, y_1, y_2, \dots$ in the finite set
$\Splus$ and  so $y_i = y_j$ for some $j > i$.  But then $\sigma^d(y_i) =
y_j = y_i$ for $d = d_i + d_{i+1} + \dots + d_{j-1} > 0$, contradicting
Lemma~\ref{fixed point and curve}(1).
 Thus the process is finite and $\sigma^e(y_0)
\in \Splusex$ for some $e  \geq 0$.

If $\sigma^{e_1}(y_0) = w \in \Splusex$ and $\sigma^{e_2}(y_0) = z \in
\Splusex$ for some (possibly negative)  integers $e_2 > e_1$, then
Lemma~\ref{Splus properties}(1) implies that  $\sigma^{e_2 - e_1}(w) = z$,
contradicting the definition of $\Splusex$.  So $e$ is unique.
\end{proof}

We next  study the special properties of the points $\Splusex$. The
intuitive idea is that these points mark the boundary of bad
behavior of $\sigma$, in the sense that for $y \in \Splusex$ the
half-orbit $\{ \sigma^d(y) : d \geq 0 \}$  behaves just like the
half-orbit of an automorphism.

\begin{proposition} \label{props of Splusex}  {\rm (1)}
Let $y\in \Splusex$. Then $\sigma^e$ is defined and
a local isomorphism at $\sigma^d(y)$ for all $e,d\geq 0$.

\begin{enumerate}
\item [(2)] Let $y\in \Splus$. The set of points $\{ \sigma^d(y) : d \geq 0 \}$ is Zariski
dense in $\YY$, and each of these points is nonsingular.
\end{enumerate}
\end{proposition}

\begin{proof}
(1)  By Lemma~\ref{Splus properties}(1),  $y_d = \sigma^d(y)$ is defined  for each
$d \geq 0$.  If $\sigma^{-1}$ is undefined at   $y_d$ for some  $d
\geq 1$ then  Lemma~\ref{Splus properties}(2) implies that
$\sigma^e(y_d) = y_{d+e} \in \Splus$ for some $e \geq 0$,
contradicting the definition of $\Splusex$.  Thus
$\sigma^{-1}$ is defined at each $y_d$ with $d \geq 1$ and so, by induction,
$\sigma^{-d}$ is defined at $y_d$ for $d \geq 1$. Thus
$\sigma^{-d}(y_d) = y$ and $\sigma^d$ is a local isomorphism at $y$. In addition,
$\sigma = \sigma^{d+1} \sigma^{-d}$ is defined at $y_d$ for all $d \geq 0$, with
$\sigma(y_d) = y_{d+1}$.  Finally, $\sigma^e$ is a local
isomorphism at every point $y_d$ with $d \geq 0$, since the inverse
map $\sigma^{-e}$ is defined at $y_{d+e}$.

(2) Given   $y \in \Splusex$, the half orbit $P = \{y_d
=  \sigma^d(y) : d \geq 0 \}$ is defined by part (1).  Suppose first
that $y_i=y_j$ for some $j > i \geq 0$. By part (1), $\sigma^{-i}$ is
defined at both $y_i$ and $y_j$, and so $y = \sigma^{-i}(y_i) =
\sigma^{-i}(y_j) = y_{j-i}$.  Then $\sigma^{j-i}(y) = y_{j-i} = y$, contradicting
Lemma~\ref{fixed point and curve}(1).
  So $y_i \neq y_j$ for $j \neq i$ and $P$ is infinite.

Suppose that the Zariski closure $\overline{P}$ of $P$ is a closed
subset of $\YY$ of dimension $1$.  Then $\overline{P}$ is a finite
union $\overline{P}= C_1 \cup \dots \cup C_e\cup D$ of irreducible
curves $C_i$ together with a (possibly empty) finite set of closed
points $D$.    Each $C_i$
  contains infinitely many of the $y_j$, and since both $\sigma$
and $\sigma^{-1}$ are defined at every closed point in $P
\smallsetminus \{y_0 \}$, both $\sigma$ and $\sigma^{-1}$ are
defined at the generic point $\mu_i$   of $C_i$ for each $i$.
  In particular, $\sigma$ and $\sigma^{-1}$ permute
the set $\{ \mu_1, \dots, \mu_e \}$.  Since $D$ is finite, if there
exists $y_d \in D$, then $\sigma^p(y_d) \in C_j$ for some $j$ and
some $p > 0$, which forces $y_d \in \sigma^{-p}(C_j) \subseteq
\bigcup C_i$, a contradiction. Thus $D = \emptyset$ and   $y \in
C_\ell$ for some $\ell$.  But   $\sigma^{e!}(\mu_\ell) = \mu_\ell$
and so Lemma~\ref{fixed point and curve}(2) implies that $C_\ell
\cap \Splus=\emptyset$, a contradiction.  The only  remaining
possibility is that $P$ is dense in~$\YY.$

Finally, since the singular locus of $\YY$ is a proper closed subset
and $P$ is dense, some point of $P$ is nonsingular.  Since $\sigma$
is a local isomorphism at each $y_d$ with $d \geq 0$, the local
rings $\mc{O}_{\YY,y_d}$ are isomorphic for all $d \geq 0$ and hence
every point  $y_d$ is  nonsingular.
\end{proof}

\begin{remark}\label{first duality}
Recall from Lemma~\ref{right left duality}(2) that working with left modules in
place of right modules   interchanges the r\^oles of  $\phi$ and $\psi$, and hence it also
interchanges  the r\^oles of $\sigma$ and $\sigma^{-1}$.
This therefore provides analogues for $(\Sminus,\sigma^{-1})$ of the results just proved for
$(\Splus,\sigma)$; formally, by combining  Lemma~\ref{right left duality}(2)
 with the left-hand versions of  Lemmas~\ref{Splus properties}, \ref{propertyD} and
 Proposition~\ref{props of Splusex},
one obtains the following results:
\begin{enumerate}
\item  Let $y \in \YY$ be a closed point where either $\sigma$ or $\phi^{-1} $ is  undefined.
Then $\sigma^{e}(y)$ is defined and in $\Sminus$ for some $e \leq
0$.
\end{enumerate}
\noindent
By analogy with the definition of  $\Splus$ in \eqref{extremal-defn}, we
 let $\Sminusex = \{ y\in \Sminus : \forall\
n\geq 1,\ \sigma^{-n}(y)\not\in \Sminus\}$ denote the \emph{extremal
elements} of $\Sminus$. \label{extremal-defn2}
\begin{enumerate}
\item[(2)]   If $y \in \Sminus$ then $\sigma$ is not defined at $y$. However
$\sigma^e(y)$ is defined for all $e\leq 0$ and $\sigma^{e}(y) \in
\Sminusex$ for a unique integer $e$, necessarily with $e \leq 0$.
\item[(3)] Let $y \in \Sminusex$.  Then  $\sigma^e$ is
defined and  a local isomorphism  at $\sigma^{d}(y)$, for all $d,e
\leq 0$. The points $\sigma^{d}(y)$   are  nonsingular for   $d \leq
0$.
\end{enumerate}
\end{remark}


\section{\Bg\ algebras}    \label{geometricity-section}
We continue to suppose that  $A$ satisfies Assumptions~\ref{assumption4}   and let
$\YY = \YY_n$ be the relevant component of the truncated point
scheme, for some fixed stable value of $n$. If we want to show that
$A$ is a \nsr, then we first have to contract $\YY$ to a second
variety $\wtY$ where $\sigma$ becomes a (biregular) automorphism.
This  is not possible under the present assumptions on $\YY.$
However, an obviously necessary condition for the existence of
$\wtY$ is for $\sigma$ to be represented by an automorphism of
\emph{some} projective model $Z$ of $K=k(\YY)$. If we make this
extra assumption then $\wtY$ does exist (see Theorem~\ref{blowing
down})  and   in this section we set
 the scene by examining the more elementary consequences of this hypothesis.

We begin with the formal  definition.

\begin{definition}
Let $B$ be a cg algebra which is a \bbc\ domain, with graded
quotient ring $Q(B) \cong K[t, t^{-1}; \sigma]$. Then $B$ is called
\emph{\bg}\ \label{bg-defn} if there is a projective   integral scheme
$Z$ with $k(Z) = K$ and an automorphism $\tau \in \aut(Z)$ inducing
$\sigma$.  If moreover $Z$ is a surface then we say that $B$ is
\emph{\btwog}. \label{btwog-defn}
\end{definition}

If $B$ is a cg domain with $Q(B) \cong K[t, t^{-1}; \sigma]$ where
$K$ is a finitely generated field extension with $\trdeg K/k = 1$,
then $B$ is automatically \bg---take $Z$ to be the unique
nonsingular model of $K$.
 The  surface case is more subtle, but by the following lemma one can at
least assume in the definition that the surface is nonsingular.
\begin{lemma}
\label{props of geometric maps} Let $B$ be \btwog, with $Q(B) = K[t,
t^{-1}; \sigma]$.  Then there exists a nonsingular surface $Z$ with
$k(Z) = K$ and $\tau \in \aut(Z)$ inducing $\sigma$.
\end{lemma}
\begin{proof}   By hypothesis, there is some projective surface
$Z'$ with  $K=k(Z')$ and $\tau' \in \aut(Z')$ inducing $\sigma$. Now
take a resolution of singularities $\pi: Z\to Z'$ following, say,
\cite[Remark~B,  p.~155]{Li})  and check that $\tau'$ lifts to an
automorphism $\tau$ of $Z$.
\end{proof}

\begin{remark}\label{geometric-maps2}
Keep the hypotheses from the lemma. If  $K$ is not the function
field of a ruled surface  then there exists a unique minimal
nonsingular model $Z$ of $K$ and so  $\sigma$ is necessarily induced
by an automorphism of $Z$   (see \cite[Proposition 7.5]{DF}).
On the other hand, if $K$ is the function field of a ruled surface,
then there do exist automorphisms $\sigma: K \to K$ such that
$\sigma$ is not induced by an automorphism of any projective model
$Z$ of $K$.  For example, if $K=k(u,v)$ is purely transcendental,
then the automorphism $\sigma$ defined by $\sigma(u)=u$ and
$\sigma(v)=uv$ is one such example. Further examples are constructed
in \cite[Remark~7.3]{DF}.
\end{remark}

For the rest of this section, we will assume that
$A$ is a \btwog\ cg domain that is generated in degree 1, with graded
quotient ring $Q = K[t,t^{-1}; \sigma]$.
 We continue to use the terminology  from Notation~\ref{S-notation};
 in particular we   fix $\YY = \YY_n$ for some stable value of
$n$ and note  that $k(\YY)=K$ with  induced birational map $\sigma :
\YY  \dra \YY.$ In Section~\ref{stable-section}, we proved various
special properties for the contracted points $\Splusex$, and by
symmetry, for the points in $\Sminusex$.  The main application of
the birationally $2$-geometric hypothesis
is to show that  these two sets are closely  interrelated.

\begin{proposition}\label{consequences of geometric}
 Let $A$ be a \btwog\ cg domain that is generated in degree 1
and assume
that $\Sminusex \not=\emptyset$; say $\Sminusex =\{s_1,\dots,s_N\}$.
Then $\Splusex$ may be uniquely written as $\Splusex  =\{t_1,\dots,
t_N\}$ so that, for each $i$, there exists a unique $d_i\geq 2$ with
$\sigma^{d_i}(s_i) = t_i$ and $\sigma^{-d_i}(t_i) = s_i$.  Moreover,
$\sigma^{d_i}$ is a local isomorphism at~$s_i$.
\end{proposition}

\begin{remark}\label{consequences2} Note that, by the proposition and Remark~\ref{first duality},
$\Sminusex \not=\emptyset$ if and only if  $\Splusex
\not=\emptyset$.   Moreover, by Remark~\ref{first duality} and
Lemmas~\ref{Splus properties}(1,2) and \ref{propertyD},
$\Splusex=\emptyset $ if and only if $\sigma$ is an automorphism of
$\YY.$
\end{remark}

\begin{proof}
Fix $y\in \Splusex$.  By Lemma~\ref{props of Splusex}, the set $P = \{ y_m = \sigma^m(y)
: m \geq 0 \}$ is defined. Moreover, $\sigma$ is defined and a
local isomorphism at each point of $P$, and $P$ is a dense set in
$\YY$ consisting of nonsingular points.  In particular $P$ is
infinite and so $y_i \neq y_j$ for $i \neq j$.
Let $g: \wt{\YY} \to \YY$   be the
normalization of $\YY$ and note that $g^{-1}$ defines an isomorphism
locally at each nonsingular point of $\YY$, in particular at each
point of $P$. Write
 $\wt{P}=g^{-1}(P)=\{\wt{y}_m= g^{-1}(y_m)\}$ and $\wt{\sigma} = g^{-1} \sigma g :
\wt{\YY}\dra \wt{\YY}$ for the induced birational map. Clearly
$\wt{P}$ is dense in $\wt{\YY}$, and $\wt{y}_i \neq
\wt{y}_j$ for $i \neq j$. By Lemma~\ref{props of geometric maps},
there exists a nonsingular model $Z$ of $K$ and a birational map $f
: Z \dra \wt{\YY}$ such that $\tau = f^{-1} \wt{\sigma} f$ is an
automorphism of~$Z$.

Since $f$ is a birational map between normal surfaces, there are at
most finitely many points where $f$ fails to be defined \cite[Lemma
V.5.1]{Ha}, say $G \subseteq Z$,   and finitely many points, say $H
\subseteq \wt{\YY}$, where $f^{-1}$ is undefined.  Since  $H$
can contain only finitely many of the $\wt{y}_i$,  there exists
 $m_0 \geq 0$ such that $f^{-1}$ is defined at $\wt{y}_m$ for $m
\geq m_0$.   Set $z_m = f^{-1}(\wt{y}_m)$ for $m \geq m_0$ and
note that $\{ z_m  : m \geq m_0 \}$ is dense in $Z$.
We compute that
\[
\tau(z_m) = \tau f^{-1}(\wt{y}_m) = f^{-1} g^{-1} \sigma g(\wt{y}_m)
= f^{-1} g^{-1} \sigma (y_m) = f^{-1} g^{-1} (y_{m+1}) = z_{m+1}
\qquad \mathrm{for\ all\  } m \geq m_0.
\]  Thus $\{z_m  :  m \geq m_0 \}$ extends to
a full dense $\tau$-orbit $\{ z_m = \tau^{m-m_0}(z_{m_0}) : m
\in \mb{Z} \}$  and again, $z_i \neq z_j$ for $i
\neq j$.

Recall that $y = y_0 \in \Splusex$, and so $\sigma^{-1}$ is not
defined at $y$ by Lemma~\ref{Splus properties}(1).  On the other
hand, since $G$ is finite we must have $z_m \not \in G$ for all $m
\ll 0$. Hence
\[
\sigma^m(y) = \sigma^{m-m_0}(y_{m_0}) =
 g  f \tau^{m-m_0} f^{-1} g^{-1}(y_{m_0}) =
g  f \tau^{m-m_0}(z_{m_0}) = g f (z_m)
\]
is defined for all $m \ll 0$.  Let $a < 0$ be the largest negative
integer for which $\sigma^a(y)$ is defined, so $a \leq -2$. Since
 $\sigma^{a+1}$ is not defined at $y$, the map  $\sigma$ is
not defined at $\sigma^a(y)$. By Remark~\ref{first duality}(1) there
exists $b\leq 0$ such that $\sigma^b(\sigma^a(y)) = \sigma^{a+b}(y)$
is defined and belongs to $ \Sminus$, and then by Remark~\ref{first
duality}(2) there exists an integer  $c \leq 0$ with
$\sigma^c \sigma^{a+b}(y) = \sigma^{a+b+c}(y) \in \Sminusex$. Set $e
= a + b + c \leq -2$.  If $e'$ is a
second integer with $\sigma^{e'}(y) \in \Sminusex$, say with $e'>e$,
then Remark~\ref{first duality}(3) implies that $\sigma^{e - e'}$ is
defined at $\sigma^{e'}(y)$, with $\sigma^{e - e'}\sigma^{e'}(y) =
\sigma^{e}(y)$, contradicting the definition of $\Sminusex$.  So $e$
is unique.

A symmetric argument   shows that for any $s_i \in \Sminusex$
there   exists a unique integer $d_i$, necessarily  with $d_i \geq 2$, such
that $\sigma^{d_i}(s_i) \in \Splusex$.  Define $t_i =
\sigma^{d_i}(s_i)$; applying the first part of the proof to $t_i$,
we obtain $e_i$ such that $z = \sigma^{e_i+d_i}(s_i) \in \Sminusex$.
The definition of $\Sminusex$ ensures that $e_i+d_i\geq 0$.
If $e_i+d_i > 0$, then  Remark~\ref{first duality}(3) implies that
$\sigma^{-e_i-d_i}$ is defined at $z$ and, of course,
$\sigma^{-e_i-d_i}(z) = s_i$, again contradicting the definition of
$\Sminusex$. Thus $e_i = -d_i$ and $z = s_i$. Finally, beginning with any $t \in
\Splusex$, there exists $e\leq -2$
  with $\sigma^e(t) = s_i$ for some $i$.
Since $\sigma^{d_i}(s_i) = t_i$ a similar argument forces $e = -d_i$
and $t = t_i$.  This defines the  bijection between $\Splusex$ and
$\Sminusex$ and it is immediate that $\sigma^{d_i}$ is
a local isomorphism at $s_i$ since the inverse map $\sigma^{-d_i}$
is defined at $t_i$.
\end{proof}

\begin{notation}\label{C-notation}  This notation will only apply when
$A$ is a \btwog\ cg domain  for which
 $\sigma$ is not an automorphism of $Y$;
equivalently, by Remark~\ref{consequences2}, when
$\Sminusex\not=\emptyset$.  In this case
 we fix the notation from Proposition~\ref{consequences of geometric} and
write $$\Delta =  \max\{ d_i : 1\leq i\leq N\}.$$
  For  $1 \leq i \leq N$ and $j \in \mb{Z}$, define
$$C_{ij} = \{ \text{closed points $z\in Y$ such that $\sigma^{-j}$ is defined
at $z$ with}\ \sigma^{-j}(z) = s_i \}$$  and write
$$
D =  \bigcup_{  {-\Delta < j < -\Delta + d_i}\atop{\phantom{\int^0}
1\leq i\leq N\phantom{\int^0}}} C_{ij} \qquad\mathrm{and}\qquad E =
\bigcup_{{0 < j < d_i}\atop{\phantom{\int^0}1\leq i\leq
N\phantom{\int^0}}} C_{ij}.$$ When $C_{ij}$ happens to be a
singleton we write $C_{ij} = \{c_{ij} \}$.
\end{notation}

\begin{example}  As an
illustration of these concepts, consider Example~\ref{simplest-eg}
for stable $n$.  Recall that, in that example, the $S^{\pm}$ are
singletons with $\Sminusex=\Sminus=\{r_{n+1}\}$ and
$\Splusex=\Splus=\{r_0\}$; thus $N=1$ with  $s_1=c_{10}=r_{n+1}$ and
$t_1=c_{1,n+1}=r_{0}$. The sets $C_{11},\dots,C_{1n}$ are the
exceptional divisors $C_j^{(n)}$ from \eqref{simplest-eg},
except that the ordering is reversed:
$C_{11}=C_n^{(n)},\dots, C_{1n}=C_1^{(n)}.$ All the other $C_{uv}$
are singletons.
 Finally   $\Delta = d_1 = n+1$,
and $\sigma^\Delta$ is defined except at the points
$\{c_{1,-n},c_{1,-n+1}\dots c_{1,-1}\} = \{r_{2n+1},\dots,
r_{n+2}\}$.
\end{example}

The following basic properties of the sets $C_{ij}$ will prove
useful.

\begin{lemma} \label{special locus}  Keep the assumptions and notation from
Notation~\ref{C-notation}.
\begin{enumerate}
\item  If $y \not \in C_{ij}$ for all $i,j$, then $\sigma^m$ is
defined and a local isomorphism at $y$ for all $m \in \mb{Z}$.
\item If $y\in C_{ij}$ for some $i$ and $j$, then $\sigma^d$ is defined
at  $y$ for all $d\leq -j$ and $d\geq -j+d_i$.
\item If $j \not \in [1, d_i-1]$,  then $C_{ij}=\{c_{ij}\}$ is a
singleton, with $c_{ij}=\sigma^j(s_i)$ a nonsingular point.  The map
$\sigma$ is defined and a local isomorphism at $c_{ij}$ for $j < 0$
and $j \geq d_i$.
\item The sets $C_{ij}$ are pairwise disjoint closed subsets of $Y.$
\item   Suppose that $y\in C_{ij}$ for some $i$ and $j$.  If $\sigma^d$ is defined
at $y$ for some $d\in \mb{Z}$, then $\sigma^d(y) \in C_{i,j+d}$.
\end{enumerate}
\end{lemma}

\begin{proof}
(1) Suppose  that $y \in Y$ is a point where some $\sigma^m$ is
either not defined or is defined but fails to be a local
isomorphism.  In the former case, set $z = y$; in the latter case,
set $z = \sigma^m(y)$ and note that $\sigma^{-m}$ is not defined at
$z$.  In either case, some nonzero power of $\sigma$ is undefined at
$z$, so by Lemma~\ref{Splus properties}(2), Lemma~\ref{propertyD}
and Remark~\ref{first duality}(1,2), we have either $\sigma^a(z) \in
\Splusex$ or $\sigma^a(z) \in \Sminusex$ for some $a\in \mb{Z}$.
Applying Proposition~\ref{consequences of geometric}, if necessary,
we conclude that $\sigma^b(z)=\sigma^c(y) \in \Sminusex$ for some
$b,c\in \mb{Z}$.  Thus
 $y \in C_{ij}$ for some~$i,j$.

(2) When $d\leq -j$, Remark~\ref{first duality}(3) implies that
$\sigma^d(y) = \sigma^{d+j}\sigma^{-j}(y)=\sigma^{d+j}(s_i)$ is
defined. On the other hand, Proposition~\ref{consequences of
geometric} implies that $\sigma^{d_i}(s_i)=t_i$ and so, for $d\geq
-j+d_i$, Proposition~\ref{props of Splusex}(1) implies that
$\sigma^d(y)=\sigma^{d+j}(s_i)=\sigma^{d-(d_i-j)}(t_i)$ is defined.

(3)  If $j \not\in [1,d_i-1]$,  then part (2) implies that
$\sigma^j(s_i)$ is defined. Hence, for any $z\in C_{ij}$ one has
$z=\sigma^j\sigma^{-j}(z)=\sigma^j(s_i)$. Thus
$C_{ij}=\{\sigma^j(s_i)\}$ is a singleton.  For $j \geq d_i$,
  Proposition~\ref{props of Splusex}   implies that
$c_{ij}$ is nonsingular and $\sigma$ is a local isomorphism at
$c_{ij}$. Similarly, when $j<0$,
 Remark~\ref{first duality}(3) says that $c_{ij}$ is nonsingular and that
 $\sigma^{-1}$ is a local isomorphism at $c_{i,j+1}$.
  Thus  $\sigma$ is a local isomorphism at $c_{ij}$.

(4,5)
 Suppose that $y \in C_{ij} \cap C_{i'j'}$, where $j \geq j'$.
Then $\sigma^{-j}(y) = s_i$ and $\sigma^{-j'}(y) = s_{i'}$. Since
$s_i \in \Sminusex$,  Remark~\ref{first duality}(3) implies that
$\sigma^{j' - j}(s_{i'})$ is defined and thus equal to  $s_i$. By
the definition of $\Sminusex$ this forces $j = j'$ and $s_i =
s_{i'}$, whence $i = i'$.

Before completing the proof of (4), we prove (5). Suppose that
$\sigma^d$ is defined at $y\in C_{ij}$ but that $\sigma^d(y)\not\in C_{uv}$ for any $u,v$.
Then part (1) implies that $\sigma^{-d-j}\sigma^d(y)=\sigma^{-j}(y)$ is defined. But this element
must equal $s_i$ and so $\sigma^d(y)\in C_{i,j+d}$. This contradiction implies that
  $\sigma^d(y)\in C_{uv}$
to begin with. But
this implies that $\sigma^{-v+d}(y)=s_u$ which, by the previous
paragraph,
  forces $u=i$ and $v= d+j$.

We return to the proof of (4) and suppose that $y$ is in the closure
of some $C_{ij}$.
 By parts (1) and (2), and for any $m\gg 0$, $\sigma^m$ is defined in
 a neighbourhood of $y$. However,  if $z \in C_{ij}$, then parts (2) (3)  and (5) imply that
$\sigma^{m}(z) = c_{i,j + m}$ for such an integer  $m$. Therefore,
by continuity, $\sigma^m(y) = c_{i,j + m}$. But then $\sigma^{-j}(y)
= \sigma^{-j-m}(c_{i,j+m})=s_i$ and so $y \in C_{ij}$ after all.
Thus $C_{ij}$ is closed.
\end{proof}

The significance of the computations in this section is that we
tightly prescribe the closed points where $\sigma^\Delta$ is either
undefined or not a local isomorphism.

\begin{corollary} \label{special locus 2}  Keep the assumptions and notation from
Notation~\ref{C-notation}.
\begin{enumerate}
\item The   points where   $\sigma^\Delta$ is undefined are all
contained in the finite set $D$.
\item The set
of points in $U=Y\smallsetminus D$ where $\sigma^\Delta$ is
not a local isomorphism is  contained in the closed (but possibly
infinite) set $E$.
\end{enumerate}
\end{corollary}
\begin{proof}
(1) Suppose that  $\sigma^\Delta$ is undefined at the  point $y$. By
Lemma~\ref{special locus}(1), $y \in C_{ij}$ for some $i,j$ and then
Lemma~\ref{special locus}(2) implies that $-\Delta < j < -\Delta +
d_i$. Thus $y \in D$. Given $i,j$ with $-\Delta < j< -\Delta+d_i$,
then $j \leq 0$ and so the set $C_{ij}$ is a singleton by
Lemma~\ref{special locus}(3). As there are finitely many such pairs
$\{i,j\}$, the set $D$ is finite.

(2) If $y \in Y \smallsetminus D$ is a point  where
$\sigma^{\Delta}$ is not a local isomorphism, then
Lemma~\ref{special locus}(1) again shows that $y \in C_{ij}$ for
some $i,j$.  Then $\sigma^{-\Delta}$ cannot be defined at $z =
\sigma^{\Delta}(y)$. By  Lemma~\ref{special locus}(2) this forces $0
< j < d_i$ and  $y \in E$.
\end{proof}


\section{Constructing the automorphism}\label{construction-section}

We continue to assume that $A$ is
 a \btwog\ cg domain, generated in degree 1, with graded
quotient ring $Q = K[t,t^{-1}; \sigma]$    and we maintain
the notation from Proposition~\ref{consequences of geometric} and
Notation~\ref{C-notation}.  We are now ready to show that the stable
scheme $\YY$ can be contracted to a second projective scheme $\wtY$
such that the birational map $\sigma$ descends to an automorphism on
$\wtY$. The idea behind the result is rather easy: If $y\not\in
\bigcup_{i,j}C_{ij}$ then Lemma~\ref{special locus}(1)  implies that
$\sigma^m$ is defined at $y$ for all $m\in \mathbb{Z}$. On the other
hand,   Lemma~\ref{special locus}(2,3) implies that $\sigma^m$ is
defined on all of $C_{ij}$ and $\sigma^m(C_{ij})=\{c_{i,j+m}\}$ is a
singleton provided that $m\not\in [-j,-j+d_i]$. Thus, by  patching
$U$ and $\sigma^m(U')$ for appropriate  open subsets $U$ and $U'$ of
$\YY$ and some $m\gg 0$,
 we may hope to construct a new scheme $Y'$
for which the images of the $C_{ij}$ are singletons $\{c_{ij}'\}$. At this point it is easy to define
an action of $\sigma$ on $Y'$ by setting $\sigma(c_{ij}')=c_{i,j+1}'$. That this idea works forms
 the content of the next theorem.

\begin{theorem}
\label{blowing down}
   Let $A$ be a \btwog\ cg domain that is generated in degree 1
   and write  $\YY=\YY_n$ for
some stable value of $n$ with its induced birational action of $\sigma$.
 Then there exists a surjective birational
morphism $\pi: \YY \to \wtY$, where $\wtY=\wtY_n$ is a projective
surface, such that:
\begin{enumerate}
\item $\wt{\sigma} = \pi \sigma \pi^{-1}$ is an automorphism of
$\wtY$,   and
\item  $\pi^{-1}: \wtY \dra \YY$ is defined except at a finite
set $P$  of nonsingular closed points of $\wtY$, each of which lies
 on an dense $\wt{\sigma}$-orbit.
\end{enumerate}
\end{theorem}

\begin{remark}
In the special case of Example~\ref{simplest-eg},  $\wtY \cong
\mb{P}^2$ and the morphism  $\pi : \YY \to \wtY$ is simply the map
$\gamma_n$ from that example that blows up  $\mb{P}^2$ at  the $n$
points $\{r_1,\dots, r_n\}$.
\end{remark}

\begin{proof}
  If $\Sminusex=\emptyset$  then  $\sigma$ is an automorphism by
  Remark~\ref{consequences2} and  so  we can simply  take
$\wtY=\YY$ and let $\pi$ be the identity.  We may therefore assume
that $\Sminusex \not=\emptyset$ and hence, by
Remark~\ref{consequences2} again, that $\Splusex\not=\emptyset$.

Define the integer $\Delta$ and the sets $C_{ij}$, $D$, and $E$  as
in Notation~\ref{C-notation}. Set $U_1 =\YY\smallsetminus D$ and
notice that, by construction, $E \subset U_1$. By
Corollary~\ref{special locus 2}(1,2) we know that $\sigma^\Delta$ is
defined at every point of $U_1$, while $\sigma^\Delta$ is defined
and a local isomorphism on the open subset
$U_1\smallsetminus E$ of $\YY.$ By the definition of the $C_{ij}$, we also know
that   $U_2 = \sigma^{\Delta}(U_1) $ is equal to $  \YY \smallsetminus E$, and so, in
particular, $U_2$ is open in $\YY$ with $\YY=U_1\cup U_2$.
By Lemma~\ref{special locus}(2), $\sigma^{\Delta}(E) = \bigcup \{
c_{ij} : 1\leq i\leq N, \Delta<j < \Delta+d_i\}$ is a finite set of
points.

We now want to construct the surface $\wtY$. Formally, we will do
this by glueing two open sets, but the construction is
 probably clearer if we first indicate what is happening at the level of elements.
 So, set $H=\YY\smallsetminus\{C_{ij} : 1\leq i\leq N,\, j\in \mb{Z}\}$, let $\wt{H}$ be
 an isomorphic  copy of $H$ and for each $i,j$ as above define an abstract point
$\wt{c}_{ij}$. Then the new scheme $\wtY$ will simply be $\wt{H}\cup
\{\wt{c}_{ij} : 1\leq i\leq N,\, j\in \mb{Z}\}$, with the projection
$\pi: \YY\to \wtY$ defined by the given isomorphism $H\to \wt{H}$
together with $\pi(C_{ij})=\wt{c}_{ij}$.

To  make this formal, and especially to define the topology on $\wtY$,
recall  that $\YY$ is covered by the   open sets $U_1$ and $U_2$. Let
$V_1 = U_2$ and $V_2 = U_2$ be two new copies of $U_2$, and set
$V_{12} = U_2 \ \smallsetminus \sigma^{\Delta}(E)$ and $V_{21} = U_1
\smallsetminus E$; thus   $V_{12}$ is an open subset of $V_1$ and
$V_{21}$ is open in $V_2$. The morphism $\sigma^{\Delta}
\vert_{U_1}$ restricts to an isomorphism $f: V_{21} \to V_{12}$. Now
think of $V_1$ and $V_2$ as new quasi-projective schemes and forget
their embeddings in $\YY.$ Then we can use the isomorphism $f$ to
glue $V_2$ and $V_1$ along their open subschemes $V_{21}$ and
$V_{12}$ to obtain a new scheme $\wtY$.  Define a morphism $\pi: \YY
\to \wtY$ by defining $\pi \vert_{U_1} = \sigma^{\Delta}: U_1 \to
U_2 = V_1$ and letting $\pi \vert_{U_2}: U_2 \to U_2 = V_2$ be the
identity map. It is immediate that these definitions agree on $U_1
\cap U_2=V_{21}$, and so $\pi$ is  indeed a well-defined morphism.
It is obviously surjective and birational. By construction,  $\pi$
is a local isomorphism except possibly at points in the set
$E=\YY\smallsetminus U_2$, and so the morphism $\pi^{-1}$ is defined
except possibly at the finite set of points $P = \pi(E) \subset
\wtY$.

As was noted earlier, if $C_{ij}\subseteq E$,  then
Lemma~\ref{special locus}(3,5) implies that
$\wt{c}_{ij}=\pi(C_{ij})\cong \sigma^\Delta(C_{ij})$ is indeed a
closed point of $\wtY$, while for $C_{ij}\not\subseteq E$ it
follows from  Lemma~\ref{special locus}(3,4) that both $C_{ij}$ and
$\wt{c}_{ij} =\pi(C_{ij})$ consist of  a single closed point.
 Since  $H=\YY\smallsetminus\{C_{ij} : 1\leq i\leq N,\, j\in \mb{Z}\}
\subset U_1\cap U_2$, this justifies the comments of the paragraph
before last.   We next want prove that the birational map
$\wt{\sigma} = \pi \sigma  \pi^{-1}: \wtY \dra \wtY$
 is actually an automorphism and once again this is particularly  easy to describe
 at the level of sets: $\wt{\sigma}(\wt{c}_{ij})=\wt{c}_{i,j+1}$, while the action of $\wt{\sigma}$
 on $\wt{H}=\pi(H)$   is   induced  from the fact that $\sigma$ is defined on $H$.

Formally, let $E'=E\cup\{s_i : 1\leq i\leq N\}$ and $U_2'=\YY
\smallsetminus E'\subseteq U_2$. Also, write
$P'=\pi(Z')=P\cup \pi(\{s_i\})$ and $V_2'=\pi(U_2') = \wtY\smallsetminus P'$.
By Lemma~\ref{special locus}(1,3),
$\sigma:U_2'\to U_2$ is defined and therefore an open immersion.
  Since $\pi \vert_{U_2} : U_2 \to V_2$ is an isomorphism,
this implies that $\wt{\sigma}=\pi\sigma\pi^{-1} : V_2'\to V_2$ is
also an immersion.
  Next, set $U_1' = U_1\smallsetminus \bigcup_i
C_{i,-\Delta}$ and $V_1'=\pi(U_1') \subseteq V_1$.
Recall that $V_1$ was defined to be a copy of $U_2$ and write this
identification  as $g: V_1 \to U_2$.  By the definition of $\pi\mid_{U_1}$, we see that
$g \pi = \sigma^{\Delta} $ as morphisms from $U_1$ to $U_2$.  Since
$g$ maps $V_1'$ isomorphically onto $U_2'$, it follows that
$\wt{\sigma}' = g^{-1} \sigma g$ is defined as a morphism from $V_1' $ to $ V_1$
and is even  an immersion. As birational maps,
 $g^{-1} \sigma g = \pi \sigma^{-\Delta} \sigma \sigma^{\Delta}
\pi^{-1} = \pi \sigma \pi^{-1}$.  Thus $\wt{\sigma}$ and
$\wt{\sigma}'$ represent the same birational map $\wtY \dra \wtY$,
and so they patch together to give a morphism $\wt{\sigma}: V_1'
\cup V_2' \to V_1 \cup V_2$.

 From the definitions of $D$ and $E$ in
Notation~\ref{C-notation} it is clear that  $\YY=U_1'\cup U_2'$.
Applying $\pi$, this  gives $V_i'\cup V_2'=\wtY$.  Thus
$\wt{\sigma}: \wtY \to \wtY$ is a globally defined birational
morphism which, by the previous paragraph, is a local isomorphism at
every point.   Thus the inverse birational map is defined at every
point of $\wtY$ and so   $\wt{\sigma}$ must be an isomorphism.  It
is also clear that, at the level of points, $\wt{\sigma}$ has the
form described in the paragraph before last; in particular
$\wt{\sigma}(\wt{c}_{ij})=\wt{c}_{i,j+1}$ for each $i,j$.

It remains to prove that $\wtY$ is projective, which follows from
the next, more general proposition.
\end{proof}

\begin{lemma}\label{immersion lemma}
  Let $W$ be an integral scheme of finite type.
   Let $\mc{M}_1, \mc{M}_2$ be invertible sheaves on $W$ with
$M_i \subseteq \HB^0(X, \mc{M}_i)$ such that each $M_i$ generates
$\mc{M}_i$.  Suppose that the map $W \to \mb{P}(M_1^*)$ is an
immersion.  If $M$ denotes the image of $M_1 \otimes M_2$ in
$\HB^0(X, \mc{M}_1 \otimes \mc{M}_2)$, then the map $W \to \mb{P}(M^*)$
is an immersion. \qed
\end{lemma}

 \begin{proposition}\label{projective}
Let $\YY$ be a projective $k$-scheme and $\wtY$ be an integral
$k$-scheme of finite type with a birational surjective morphism
$\pi: \YY \to \wtY$. Assume that
 $\wtY$ has an automorphism $\tau$ and a finite set $P$ of nonsingular closed points such  that
 \begin{enumerate}
 \item[(i)] $\pi^{-1}$ is defined on $U=\wtY\smallsetminus P$, and
 \item[(ii)]  every point of $P$ lies  on an infinite orbit of $\tau$.
\end{enumerate}
Then $\wtY$ is a  projective scheme.
\end{proposition}

\begin{proof}
We begin by  proving that $\wtY$ is separated and proper, for which
we use  the valuative criteria from \cite[Section~II.4]{Ha}.  By
\cite[Exercise~II.4.11(c)]{Ha}, we need only consider discrete
valuation rings in the criteria.
 So,  fix a discrete valuation ring $R$ with field of fractions $L$ and
write $\iota: B=\spec L \hookrightarrow C=\spec R$  for the
induced morphism. Let $c$ denote the closed point of $C$ and write
$b$ for the generic point of $C$ (and $B$).  We
also have the structure map $f:\wtY \to \spec k$ and
assume that there exist morphisms $\alpha:C\to \spec k$, $\beta : B\to \wtY$ such that
$\alpha \iota = f \beta$.

We first prove that $\wtY$ is separated, so assume that there exist
two morphisms $\gamma_i: C\to \wtY$ both of which make
(\ref{separate-equ})(I)  commute.
\begin{equation}\label{separate-equ}
(I) \quad
     \xymatrix@C=6.2em@R=5.2em{ B \ar[r]^{\beta} \ar[d]_{\iota} &
 \wtY \ar[d]^{f}   \\
     C \ar[r]_{\alpha}
     \ar@<-2pt>[ru]^{\gamma_1}\ar@<+2pt>[ru]_{\gamma_2} &
    \spec k}
 \qquad {}\qquad   \qquad (II) \quad
     \xymatrix@C=6.2em@R=5.2em{ B \ar[r]^{\pi^{-1} \tau^m \beta} \ar[d]_{\iota} &
 \YY \ar[d]^{f \tau^{-m} \pi}   \\
     C \ar[r]_{\alpha}
     \ar@<-2pt>[ru]^{\pi^{-1} \tau^m \gamma_1}\ar@<+2pt>[ru]_{\pi^{-1} \tau^m \gamma_2} &
    \spec k}
    \end{equation}
Write $y_1=\gamma_1(c)$ and $y_2 = \gamma_2(c)$. By assumption (ii)
we choose $m \gg 0$ such that $\tau^m(y_1)$ and $\tau^m(y_2)$ lie in
$U = \wtY\smallsetminus P$; this forces $\tau^m \beta(b) \in U$ as
well.  Then $\pi^{-1} \tau^m \gamma_i: C \to \YY$ is a well-defined
morphism for $i = 1, 2$, and similarly $\pi^{-1} \tau^m \beta: B \to
\YY$ is well-defined.  Moreover, the diagram
(\ref{separate-equ})(II) commutes.  Since $\YY$ is separated,  \cite[Theorem~II.4.3]{Ha}
implies that $\pi^{-1} \tau^n \gamma_1=\pi^{-1} \tau^n \gamma_2$ and
hence $\gamma_1=\gamma_2$. By \cite[Theorem~II.4.3]{Ha}, again, this
implies that $\wtY$ is separated.

We now  prove properness. By hypothesis,  there exists $m\gg 0$ such that
$\wtY$ is covered by the open sets $U\cong \pi^{-1}(U)$ and
$\tau^{m}(U)$, so  $\wtY$ is noetherian. By  the Valuative Criterion
for Properness \cite[Theorem~II.4.7]{Ha} it therefore suffices to
find a morphism $\gamma: C\to \wtY$ making \eqref{separate-equ4}(I)
 commute.
\begin{equation}\label{separate-equ4}
(I)\quad  \xymatrix@C=6.2em@R=5.2em{ B \ar[r]^{\beta} \ar[d]_{\iota}
&
 {\wtY} \ar[d]^{f} \\
     C \ar[r]_{\alpha}  \ar@{-->}     [ru]^{\gamma}
     &    \spec k}
 \qquad {}\qquad   \qquad (II) \quad
   \xymatrix@C=6.2em@R=5.2em{ B \ar[r]^{\pi^{-1}\beta} \ar[d]_{\iota} &
 {\YY} \ar[d]^{f\pi} \\
     C \ar[r]_{\alpha}  \ar@{-->}[ru]^{\gamma'} &
    \spec k}
\end{equation}
First, if $x=\beta(b)$ is a closed point of $\wtY$ then clearly
putting $\gamma(c) = \gamma(b) = y$ gives a well-defined map
$\gamma$ making \eqref{separate-equ4}(I) commute.  Alternatively, if
$x$ is not a closed point, then $x \not\in P$ and so $\pi^{-1}$ is
defined at $x$.  As $\YY$ is proper, there then exists $\gamma' :
C\to \YY$ making the diagram \eqref{separate-equ4}(II) commute. Thus
$\gamma=\pi\gamma'$ makes  \eqref{separate-equ4}(I) commute and
$\wtY$ is proper.

Before proving that $\wtY$ is projective, we prove the following easy facts.
 The constant sheaf of rational functions on $\YY$  and on $\wtY$  will be written
 $\mc{K}$.

\begin{sublemma}\label{pushforward1}  Keep the notation from the proposition.
Suppose that $\mc{M} \subset \mc{K}$ is
an invertible sheaf  on $Y$
and set $\mc{F}=\pi_*\mc{M}\subset \mc{K}$.  Then:
\begin{enumerate}
 \item The reflexive hull $ \mc{L} =  \mc{F}^{**}=
 \calHom_{\mc{O}_{\wtY}}( \calHom_{\mc{O}_{\wtY}}
 (\mc{F}, \mc{O}_{\wtY}), \mc{O}_{\wtY})$\label{reflexivehull-defn}
 of $\mc{F}$ is invertible.
\item  $\mc{L}$ is the unique reflexive (or invertible) sheaf
on $\wtY$ satisfying $\mc{L}\vert_{U} = \mc{F}\vert_U$.
 \item  $\mc{L}^{\otimes n}$ is the reflexive hull of $\mc{F}^{n}$, where the
 multiplication takes place inside $\mc{K}$.
 \end{enumerate}
 \end{sublemma}

 \begin{proof}
 (1)    Since    $\mc{F}\vert_{U}
= \mc{M}\vert_{\pi^{-1}(U)}$ under the identification of $\mc{O}_U$
with $\mc{O}_{\pi^{-1}(U)}$,   the sheaf $\mc{L}$ is locally free on $U$.
On the other hand, if $p\in P$, then $\mc{O}_{\wtY,p}$ is regular
and hence factorial   by  \cite[Theorem~19.19]{E}
 and so  any reflexive ideal of  $\mc{O}_{\wtY,p}$ is cyclic. Thus $\mc{L}$
 is locally free at each point of  $\wtY$.

 (2) We need only prove this locally at a closed point $y\in P$.
 If $\mc{Q}$ is a second reflexive sheaf with $\mc{Q}\vert_{U}=\mc{F}\vert_{U}$,
 then  certainly  $(\mc{Q}_y)_{\mf{p}}=(\mc{F}_y)_{\mf{p}}=(\mc{L}_y)_{\mf{p}}$
 for all height one prime   ideals $\mf{p}$
   of the regular local ring $\mc{O}_{\wtY,y}$. But
      $\mc{Q}_y=\bigcap\{(\mc{Q}_y)_{\mf{p}} : \mf{p} \mathrm{\ height\ } 1\}$. Thus
         $\mc{Q}_y=\mc{L}_y$ and hence $\mc{Q}=\mc{L}$.

   (3)  Since $\mc{L}^{\otimes n}\vert_{U}=\mc{F}^{n}\vert_{U}=(\pi_*(\mc{M}^{\otimes n}))\vert_{U}$,
   this follows from part (2) applied to $\mc{M}^{\otimes n}$.
 \end{proof}

  We return to the proof of the proposition and
 prove that $\wtY$ is projective. Let $\mc{M}\subset\mc{K}$
be any very ample invertible sheaf on $\YY$ and set $\mc{F} = \pi_*
\mc{M}$  with reflexive hull   $\mc{L}=\mc{F}^{**}$. By the sublemma
$\mc{L}$ is locally free.
  Since $\mc{M}$ is generated by its  sections $M=\mathrm{H}^0(\YY,\mc{M})$ and
 $M =\mathrm{H}^0(\wtY,\mc{F})\subseteq \mathrm{H}^0(\wtY,\mc{L})$,
 the sections  $\mathrm{H}^0(\wtY, \mc{L})$ do at least generate
  the sheaf  $\mc{L}\vert_U=\mc{F}\vert_U\cong \mc{M}\vert_{\pi^{-1}(U)}$ as an $\mc{O}_U$-module.
  Since  $P=\wtY\smallsetminus U$ is finite and $\wtY$   is proper,
    it therefore follows from   \cite[Corollary~1.14]{Fj1}  that,  for any
$n \gg 0$,  the sheaf $\mc{L}^{\otimes n}$ is   generated by its global  sections
   $V_1=\mathrm{H}^0(\wtY, \mc{L}^{\otimes n})$.
   For $r\in \mb{Z}$, write $\mc{N}_r=(\tau^{-r})^*(\mc{L}^{\otimes n})$.

  Let   $\theta: \wtY \to  \mb{P}(V_1^*)=\mb{P}^d$  denote the map determined by  $V_1$ and write
   $V_2= \mathrm{H}^0(\YY, \mc{M}^{\otimes n})$. Then  Sublemma~\ref{pushforward1}(3)
   implies that
   $V_2=   \mathrm{H}^0(\wtY,\, \pi_* \mc{M}^{\otimes n})\subseteq
    \mathrm{H}^0(\YY, \mc{L}^{\otimes n})=V_1.$
Since  $\mc{M}^{\otimes n}$ is very ample, $V_2$ defines a closed  immersion
of $\YY$ into $\mb{P}(V_2^*)$ and hence an immersion
  $\theta':\pi^{-1}(U)\hookrightarrow \mb{P}(V_2^*)$.  Since $\theta'$ factors through $\theta$,
  this ensures that $\theta$   defines an immersion of   $U$ into $\mb{P}(V_1^*)$.

By  hypothesis (ii) of the proposition,    $\wtY = U\cup \tau^{m}(U)$ for  $m\gg 0$
  and, for some such $m$,  we  consider  the  map $\alpha: \wtY\to  \mb{P}(V_3^*)$
  determined by   $V_3=\mathrm{H}^0(\wtY,\,\mc{N}_0\otimes \mc{N}_m)$.
Since the  $\mc{N}_r$ are generated by their sections, we can apply
 Lemma~\ref{immersion lemma}  with $W=U$ to show that the   restriction of $\alpha$ to $U$
 is  an immersion.  For exactly  the same reasons  the map  $\theta_1: \wtY \to
\mb{P}^d$ determined by    $\mathrm{H}^0\big(\wtY,\, \mc{N}_m)$
immerses $\tau^m(U)$ into $\mb{P}^d$ and  so the   restriction of  $\alpha$ to $\tau^{m}(U)$  is   also an immersion.

 Since $\alpha$ is defined on all of $\wtY = U \cup \tau^m(U)$, we can  define
  $C$ to be  the closure in $\mb{P}(V_3^*)$ of $\alpha(\wtY)$. Since $\wtY$
is proper and irreducible, $C=\alpha(\wtY) $ and obviously   the
corresponding map $\beta: \wtY \to C$  is birational.
 The conclusion of the previous  paragraph implies that
$\beta$ is a local isomorphism
at every  point of $\wtY$.  Therefore,  for any  point $c\in C$ we can
define an inverse map to $\beta$ in a neighbourhood of $c$ and so the
inverse birational map $\beta^{-1}$  is defined at every
point of $C$.  Hence $\beta: \wtY \to C$ is an isomorphism.  In
other words, $\alpha: \wtY \to \mb{P}(V_3^*)$ is a closed immersion
and  $\wtY$ is projective.
\end{proof}

 \begin{remark}\label{projective-remark}
For future  reference we note that the last part of the proof of
Proposition~\ref{projective} proves
 the following fact: Suppose that $\mc{M}$ is a very ample invertible sheaf on $\YY$
 and let $\mc{B}=\bigl(\pi_*\mc{M})^{**}.$ Then,  for all  $q\gg 0$, the sheaf
 $\mc{B} \otimes (\tau^q)^*\mc{B}$
 is an ample invertible sheaf on $\wtY$.
\end{remark}

Although the  construction of the   scheme $\wtY_n$   in
Theorem~\ref{blowing down} depends upon  $n$, this is unimportant as
there is a   natural isomorphism $\wtY_n\cong \wtY_{n+1}$ for any
stable value of $n$. Before proving this we need  a technical lemma.

\begin{lemma}\label{ptmodule-lemma}
Let $y \in \YY = \YY_n$, for a fixed stable value of $n$.  Then
there exists $m_0\geq 0$ such that, for all $m \geq m_0$ and $r \geq
0$,  $\sigma^m$ is defined at $y$ and the birational map $\phi^{-r}$
is defined at $\sigma^m(y)$.
\end{lemma}
\begin{proof}   Let $w\in \YY$ and suppose that $\phi^{-r}$ is not defined at $w$ for
some   $r \geq 1$. If we choose the minimal such
$r$  then $\phi^{-1}$ is undefined at $\wt{w} = \phi^{-r+1}(w) \in
\YY_{n + r -1}$. Since $n + r -1$ is again a stable value,
  Remark~\ref{first duality}(1) implies
   that $\sigma^e(\wt{w}) \in \Sminus_{n+r-1}$ for some $e
\leq 0$.  Thus $\sigma^{e +r -1} = \psi^{r-1} \sigma^e \phi^{-r+1}$
is defined at $w$ and, by Lemma~\ref{exceptional points}(1),
$\sigma^{e+r-1}(w) \in \Sminus$. Using Notation~\ref{C-notation} and
Remark~\ref{first  duality}(2), it follows that $w \in C_{i,j}$ for some $i,j$.

Return to the proof of the lemma and suppose that $y \not \in C_{ij}$ for any $i,j$.  Then
Lemma~\ref{special locus}(1)  implies that  $y_m = \sigma^m(y)$ is defined for all
$m$.  If $\phi^{-r}$ is undefined at $w=y_m$ for some $m \geq
0$ and $r \geq 0$, then the previous paragraph
shows that $y_m \in C_{ij}$ for  some $i,j$.  This forces $y \in C_{i,
j-m}$, a contradiction.  Thus in this case we may take $m_0 = 0$.

Otherwise, $y \in C_{ij}$ some $i,j$, and so $\sigma^{e}(y) = t_i
\in \Splusex$ for $e = d_i -j$.  Set $z = t_i$ and notice that, by  Lemma~\ref{exceptional points}(2),
$\phi^{-r}$ is defined at $z$ for all $r \geq 0$.  For  $d \geq 0$, the map
$\sigma^d$ is defined and a local isomorphism at $z$ by
Proposition~\ref{props of Splusex}.  Hence $\phi^{-r} = \psi^d
\phi^{-r - d} \sigma^{-d}$ is also defined at $z_d = \sigma^d(z)$
for all $r \geq 0$.  Finally, $\sigma^{e +d}(y) = z_d$ for all $d
\geq 0$ and so the conclusion of the lemma holds with  $m_0 = \max(0, e)$.
 \end{proof}

\begin{lemma}\label{mainthm-lemma1}
The scheme   $\wtY_n$ is independent of the stable value of   $n$
 in the sense that, for any such $n$, there exists an isomorphism $\alpha_n:\wtY_{n+1}\to \wtY_{n}$
such that $\pi_n \phi_n = \alpha_n \pi_{n+1}$.
\end{lemma}

\begin{proof}
Let $y\in \wtY_{n}$.  By Theorem~\ref{blowing down}(2),
$\pi_{n}^{-1}$ is defined at $\wt{\sigma}_n^e(y)$ for all $e\gg 0$,
say for all $e\geq e_0$.  Set $z=\wt{\sigma}_n^{e_0}(y)$. Then by
Lemma~\ref{ptmodule-lemma}, we may choose $e_1 \geq 0$ such that
$\phi^{-1}_n$ is defined at $\wt{z}_e = \sigma_n^e
\pi_n^{-1}(z)=\pi_n^{-1}\wt{\sigma}_n^e(z) $, for all $e \geq e_1$.
So, for all $w$ in some neighborhood of $y$ we may define
$$\alpha_{n,m}^{-1}(w) =
 \wt{\sigma}_{n+1}^{-m}\circ\pi_{n+1}\circ\phi_n^{-1}\circ\pi_n^{-1}\circ\wt{\sigma}_n^{m}(w),
\qquad\mathrm{for}\quad m\geq m_0 = e_0 + e_1.$$ By
Lemma~\ref{sigma-lemma},   we have the equality
$\sigma_{n+1}\phi_n^{-1}=\phi_n^{-1}\sigma_n$ of birational maps.
Also, $\wt{\sigma}_n\pi_n = \pi_n\sigma_n$ as birational maps by the
definition of $\wt{\sigma}_n$.  Thus $\alpha_{n,m}^{-1} = \wt{\sigma}_{n+1}^{-m}
\pi_{n+1} \phi_n^{-1} \pi_n^{-1} \wt{\sigma}_n^{m} = \pi_{n+1}
\phi_n^{-1} \pi_n^{-1}$ as birational maps, independent of the
choice of $m \geq m_0$.  Thus the birational map $\alpha_n^{-1} = \pi_{n+1}
\phi_n^{-1} \pi_n^{-1}$ is defined in a neighborhood of $y$ which,  since
$y$ was arbitrary, implies that $\alpha_n^{-1}$ is globally defined. A similar
(but easier) argument shows that the birational map $\alpha_n =
\pi_{n} \phi_n \pi_{n+1}^{-1}$ is globally defined, and so
$\alpha_n$ is an isomorphism.
\end{proof}

We end the section by giving a more algebraic interpretation of the
scheme $\wtY$; the set of closed points in $\wtY$ is canonically
isomorphic to a set $\wt{P}_{\mathrm{qgr}}$ of $k$-point modules
in $\rqgr A$.
 This is both suggestive and curious. It is curious
since, as is shown by \cite[Theorem~0.1]{KRS},
$\wt{P}_{\mathrm{qgr}}$ is definitely \emph{not} represented by a
scheme of finite type when $A$ is a \nsr. However, it
does suggest that there should be some category of modules over $A$
which is represented by $\wtY$.

Formally,  for $m\geq n$,  a closed point $y\in \YY_m$ can be
identified with a truncated $k$-point module $M_y$  of length $m+1$.
Mimicking the definition in Section~\ref{sect-pts}, let
$\wt{\mc{P}}_m$ denote the set of isomorphism classes of all such
truncated point modules and let
$\wt{\mc{P}}=\displaystyle\lim_{\leftarrow}\wt{\mc{P}}_m$ denote the
set of isomorphism classes of point modules corresponding these
truncated modules. Finally, let $\wt{\mc{P}}_{\mathrm{qgr}}$ denote
the image of $\wt{\mc{P}}$ in $\rqgr A$.
Although we will not prove it here, the set
$\widetilde{\mc{P}}_{\mathrm{qgr}}$ is  the image in qgr of the set
of all point modules from $\rgr A$
and hence also the set of all point modules in $\rqgr A$ in the 
sense of  the introduction.

\begin{corollary}\label{ptmodule-cor}  Pick a stable value of $n$ and set $\wtY=\wtY_n$.
There is a bijective correspondence $\rho$ from the  set of closed
points in the scheme $\wtY$ to the set of objects in
$\wt{\mc{P}}_{\mathrm{qgr}}$ with the following property: Given
$\wt{y}\in \wtY$, write $\wt{y}=\pi(y)$ for some $y\in \YY$, say
corresponding to the truncated point module $M=M_y\in \wt{P}_n$.
Then $\rho(\wt{y})$ is the image $\overline{N}$ in $\rqgr A$ of
\emph{any} point module $N \in \wt{\mc{P}}$ satisfying $M \cong
N_{\leq n}$.
\end{corollary}

\begin{proof}   Let $y \in \YY = \YY_n$, for a fixed stable value of $n$
and recall the definitions from Notation~\ref{C-notation} and Theorem~\ref{blowing down}.
We first claim that  the map $\pi': \wt{P}_n \to \wt{\mc{P}}_{\mathrm{qgr}}$
given by $y\mapsto \overline{N}$ is well-defined.

To prove this,  pick   $m \geq 0$ by  Lemma~\ref{ptmodule-lemma} so that
   the map  $\phi^{-r}$ is  defined at $\sigma^m(y)$ for all $r\geq 0$.
 We now  translate this
into module theory.
 Let $M_y \in \wt{\mc{P}}_n$ be the truncated point
module corresponding to the point $y\in \YY_n$.  By
Corollary~\ref{zero-ideal}(2),
 the maps $\phi_d$
are surjective for $d \geq n$,  and so there is a sequence of points
$y_d \in \YY_d$ with $y_n = y$ and $\phi(y_{d+1}) = y_d$ for all $d
\geq 0$. Corresponding to this sequence  is a point module $N$ such
that $N_{\leq n} \cong M_y$.   The  truncation
shift $L = N[m]_{\geq 0}$ is also a  truncated point module but now,
by (\ref{rtrunc}) and (\ref{ltrunc}),
  $L_{\leq n}$ corresponds to the
point $\sigma^m(y)\in \YY$.  Since $\phi^{-r}$ is defined at
$\sigma^m(y)$ for all $r \geq 0$, $L$ must be the unique point
module $L'$ with $L'_{\leq n} \cong L_{\leq n}$.  In other words,
$N[m]_{\geq 0}$ and hence $N_{\geq m}$ are uniquely determined by
$M_y$ and $y$. Since  $N$ and $N_{\geq m}$ have the same image in
$\wt{\mc{P}}_{\mathrm{qgr}}$, this says that  $\pi'$ is well
defined.

In order to  show  that the correspondence $\rho$ is well-defined
and bijective,  it suffices to prove the following: given closed
points $x_1\not=x_2 \in \YY$ then  $\pi(x_1) = \pi(x_2)$ if and only
if $\pi'(x_1) =\pi'(x_2)$.
 By the construction of $\pi$ in Theorem~\ref{blowing down},  $\pi(x_1) =\pi(x_2)$ if
and only if $x_1, x_2 \in C_{ij}$ for some $i,j$. But in this case,
Lemma~\ref{special locus}(2,3,5) implies that $\sigma^m(x_1) =
\sigma^m(x_2) = c_{i,j +m}$ for $m \gg 0$. Now suppose that $N(i)$
are point modules with $N(i)_{\leq n}=M_{x_i}$ for $i = 1,2$. Then,
just as in the previous paragraph, it follows from (\ref{rtrunc}) and
(\ref{ltrunc}) that $L(i)=N(i)[m]_{\geq 0}$   is a truncated point module corresponding
to $\sigma^m(x_i) $, which  is independent of $i$. In other words,
$N(1)$ and $ N(2)$ have isomorphic tails, $\ol{N(1)}\cong
\ol{N(2)}$, and $\rho$ is well-defined. Conversely, if
$\pi'(x_1)=\pi'(x_2)$, then it follows from the argument of the previous
paragraph that $\sigma^m(x_1)=\sigma^m(x_2)$ for all $m\gg 0$. Thus
$x_1=x_2\in C_{ij}$ for some $i,j$ and so $\pi(x_1)=\pi(x_2)$.
\end{proof}


\section{Na{\"\i}ve blowing up}\label{naive-section}

We now have the data required for the statement of the
 Main Theorem~\ref{intro-mainthm} and in the next three
sections we show how to use this to describe the relevant cg domains as \nsr s.
In this section we give the definitions and some of the basic properties
of  these algebras.
 Despite their name,    \naive\  blowups  have  properties that are very different from the classical
 case---just look at Corollary~\ref{intro-maincor}.
They  have been studied  in detail in \cite{KRS}; but only  for the case where
  one is   blowing up a single closed point, whereas  for applications in this paper we
  need to be able to \naive ly blow up a zero dimensional  subscheme
 of the given scheme and this complicates several of the results from \cite{KRS}.  A
 detailed examination of these more general rings is given in the companion paper
  \cite{RS}, to which the
 reader is referred for the details behind the assertions given here.

The following assumptions will be in force throughout the section.

\begin{assumptions}\label{assumptions7} Fix an integral projective scheme $X$ of dimension
$d\geq 2$ with an automorphism $\tau$ and
  field of rational functions   $K=k(X)$.
The corresponding  sheaf of rational functions will be written
 $\mc{K}$.\label{mcK-defn}
 Let $\LL\subset\mc{K}$ be an invertible sheaf   and let
 $\mc{I}$ be an   ideal sheaf such that the associated scheme
 $Z_{\mc{I}}=\mc{V}(\mc{I})$ is zero dimensional.
   We allow $\mc{I}=\mc{O}_X$ but, to avoid trivialities,
 we assume that if $s\in Z_{\mc{I}}$  then $s$
 lies on an infinite $\tau$-orbit in $X$.
\end{assumptions}

    Given an $\mc{O}_X$-module $\mc{M}$,   write $\tau^*(\mc{M})=\mc{M}^\tau$. If
   $\mc{M}\subseteq \mc{K}$, set $\mc{M}_0=\mc{O}_X$ and
   $\mc{M}_n=\mc{M}\mc{M}^\tau\cdots\mc{M}^{\tau^{n-1}}$ for $n>0$, under the  multiplication
induced   from $\mc{K}$. Of course,
    $\mc{M}_n\cong\mc{M}\otimes \mc{M}^\tau\otimes \cdots\otimes \mc{M}^{\tau^{n-1}}$
    when $\mc{M}$ is invertible \label{Ln-defn}
but  we will also need  use the construction for
  $\mc{M}= \LL\mc{I}$ where $\mc{M}\mc{M}^\tau$ is typically not isomorphic
  to $\mc{M}\otimes \mc{M}^\tau$.

  Given this data, and regarding  $\mc{K}[t,t^{-1};\tau]$ as a   sheaf of $\mc{O}_X$-algebras
in the natural way,   we can form the \emph{bimodule algebra}
 \begin{equation}\label{bimodule-defn}
 \mc{R}(X,Z_{\mc{I}},\LL,\tau) \ = \
  \bigoplus_{n\geq 0}\mc{R}_nt^n\ \subseteq \ \mc{K}[t,t^{-1};\tau],
 \qquad\mathrm{where}\quad
\mc{R}_n=\mc{L}_n \mc{I}_n\cong\mc{L}_n\otimes\mc{I}_n\quad\mathrm{for}\ n \geq 0.
\end{equation}
The definition of  $\mc{R}=\mc{R}(X,Z_{\mc{I}},\LL,\tau)$ in
\cite{KRS} or \cite{RS} is slightly different---although
equivalent---since it is given in terms of $\mc{O}_X$-bimodules. One
advantage of the present definition is that the multiplication in
$\mc{R}$ is induced from that in  $\mc{K}[t,t^{-1};\tau]$
and is given by the natural  map
 \begin{equation}\label{bimodule-defn2}
 \mc{R}_mt^m\otimes \mc{R}_nt^n \longrightarrow  \mc{R}_mt^m  \mc{R}_nt^n
=  \mc{R}_m\mc{R}_n^{\tau^m}t^{m+n} = \mc{R}_{m+n}t^{m+n}.
\end{equation}
We remark that, as is noted in \cite[p.~252]{AV},  $ \mc{R}(X,Z_{\mc{I}},\LL,\tau)$
 is not a sheaf of algebras
 in the usual sense, since  the presence of $\sigma$ means one has to play a game of
 musical chairs with the open sets.
 The \emph{\nsr\ of $X$ at $Z_{\mc{I}}$ }  is then defined to be the  ring of global
sections
 \begin{equation}\label{naive-defn}
  R(X,Z_{\mc{I}},\LL,\tau) \ = \  \HB^0(X,\mc{R})
=\bigoplus_{n\geq 0} R_n = \bigoplus_{n\geq 0} \overline{R}_n t^n
 \qquad\mathrm{where}\quad
\overline{R}_n = \HB^0(X,\,\mc{I}_n\otimes \mc{L}_n)\quad\mathrm{for}\ n \geq 0.
\end{equation}

As a matter of convention we write a graded subspace of $R$ or $Q(R)$ as
 $I=\bigoplus I_n=\bigoplus \overline{I}_nt^n$,  for some $\overline{I}_n\subset k(X)$,
  but we do not use the overline for subsets of $\mc{R}$.

Although  we have implicitly assumed that  $\mc{I}\not=\mc{O}_X$ in these definitions, they
 make perfect sense when $\mc{I}=\mc{O}_X$.   In that case
one obtains the   bimodule algebra $\mc{B}(X,  \mc{L},\sigma)
=\RR(X,\emptyset, \LL,\tau)$ whose  ring of global sections is  the familiar
\emph{twisted homogeneous coordinate ring}\label{twisted-defn}
$B(X,\LL,\tau) = R(X,\emptyset, \LL,\tau) =\bigoplus_{n\geq 0}\mathrm{H}^0(X,\,\mc{L}_n)t^n$
 from \cite{AV}.  We note in passing
that  $ R(X,Z_{\mc{I}},\LL,\tau)\subseteq B(X,\LL,\tau) $ for any $Z_{\mc{I}}$.

 By definition,  a \emph{graded right $\mc{R}$-module}\label{scriptR-mod}
  is a quasi-coherent $\OO_X$-module
 $\mc{M}=\bigoplus_{i\in\mb{Z}}\mc{M}_n$ together with a right $\OO_X$-module map
 $\mu:\mc{M}\otimes\mc{R}\to \mc{M}$ satisfying the usual axioms.
  The module $\mc{M}$ is called \emph{coherent}\label{coherent-defn}
   (as an  $\mc{R}$-module) if   each
$\mc{M}_n $ is a coherent $\OO_X$-module, with   $\mc{M}_n = 0$
 for $n \ll 0$, such that
the natural map $\mu_n:\mc{M}_n\otimes \mc{R}_1^{\tau^n}\to \mc{M}_{n+1}$
is surjective for   $n\gg 0$.   (By  \cite[Lemma~2.11]{RS}, this is equivalent
 to the somewhat more general definition from \cite{VdB1}.
  The superscript of $\tau^n$ that appears here  does not appear
 in that lemma   since the two papers use slightly different notation; in \cite{RS}
$\mc{R}_1$ is a bimodule into which one absorbs  the action of $\tau$.)
The bimodule algebra $\RR$ is \emph{(right) noetherian} if all (right) ideals of $\RR$ are coherent.

The usual definitions for graded modules over cg algebras also
apply to bimodule algebras. In particular,
the graded right $\mc{R}$-modules form an abelian category $\rGr \mc{R}$, with
homomorphisms graded  of degree zero. The subcategory of coherent modules
  is denoted $\rgr \mc{R}$ and  is abelian when $\mc{R} $  is
  noetherian \cite[Proposition~2.10]{KRS}.
 Let  $\rtors \mc{R}$ denote  the full subcategory of $\rgr \mc{R}$ consisting of
coherent $\mc{O}_X$-modules, and write $\rqgr \mc{R}$ for  the
quotient category $\rgr \mc{R}/\negthinspace\rtors \mc{R}$.
\label{script-Rmod-defn}  The category  $ \rqgr \mc{R}$ is called \emph{the
\naive\ blowup } of $X$ at $Z_{\mc{I}}$ (see, also, Theorem~\ref{VdB
main theorem}).\label{naive-defn2}

\Nsr s   only behave  well if the sequence
  $ \{\mc{R}_n=\mc{L}_n\II_n \}$ is ample in the following sense.
Assume that  $\mc{M}\subset \mc{K}$. Then  the sequence $\{\mc{M}_n : n\geq 0\}$ is   \emph{ample},
equivalently   $\mc{M}$ is  \emph{$\tau$-ample},\label {ample def} if
for every  coherent $\mc{O}_X$-module   $\mc{F} $
 the sheaf $\mc{F} \otimes_{\mc{O}_X} \mc{M}_n$
is generated by its global sections   and
$\coH^i(X,\mc{F} \otimes \mc{M}_n) = 0$ for all $i > 0$ and all $n \gg 0$.
If $\mc{M}$ is an invertible sheaf, then $\tau$-ampleness is a subtle but
 well-understood concept (see \cite{Ke1}), but much less is known
 when $\mc{M}$ is not invertible.

A basic method for relating $\mc{R}$ to $R$ is given by  the following result;
see  \cite[Theorem~2.12]{KRS} or,    for a more general version,  \cite[Theorem~5.2]{VdB1}:

 \begin{theorem} {\rm (Van den Bergh)}
\label{VdB main theorem}  Keep   Assumptions~\ref{assumptions7} and let
 $\mc{R} = \mc{R}(X,Z_{\mc{I}},\LL,\tau)$.
Assume that  $\mc{R}$ is   noetherian  and that
  $\{ \mc{R}_n  : n \in \mb{N}\}$ is a
  ample sequence of $\mc{O}_X$-modules. Then
the \nsr\ $R=R(X,Z_{\mc{I}},\LL,\tau)$
 is  noetherian, and
there is an equivalence of categories $ \rqgr \mc{R} \to  \rqgr R$ given by
$\mathrm{H}^0(X,\,-)$.   \qed
\end{theorem}

In order to use this theorem to pass from $R$ to $\mc{R}$
  we need to know when $\mc{R}$ is noetherian and when
 $\{\mc{R}_n\}$ is an ample sequence. The answer is provided by the next two results.

\begin{theorem}\label{ampleness}  Keep  Assumptions~\ref{assumptions7}
and assume that $\mc{L}$ is $\tau$-ample and
that each point $z\in Z_{\mc{I}}$ lies on a Zariski dense $\tau$-orbit.   Then the
sequences   $\{\mc{R}_n= \mc{L}_n \otimes \mc{I}_n \}$
 and $\{\mc{R}_n^{\tau^{-n}}\} $ are  ample.
\end{theorem}

\begin{proof}    The ampleness of   $\{\mc{R}_n\}$ is  proved in
 \cite[Theorem~3.1]{RS}.  In order to prove
  the second claim write $\RR_n^{\tau^{-n}} = \mc{S}\otimes \mc{S}^\alpha
 \otimes\cdots\otimes \mc{S}^{\alpha^{n-1}}$ for
 $\mc{S}=\mc{L}^{\tau^{-1}}\otimes \mc{I}^{\tau^{-1}}$ and $\alpha=\tau^{-1}$. By
  \cite[Corollary~5.1]{Ke1}, $\mc{L}^{\alpha}$ is $\alpha$-ample and
  clearly the other hypotheses of the theorem also hold for $\alpha$ and $\mc{S}$.
So   \cite[Theorem~3.1]{RS} also shows that $\mc{S}$ is
$\alpha$-ample; equivalently, that $\{\mc{R}_n^{\tau^{-n}} \} $ is
ample.
    \end{proof}

The noetherian property  for $\RR$ and $R$
 requires stronger hypotheses.  For $z \in X$, the
$\tau$-orbit $\{\tau^i(z) : i\in\mb{Z}\}$ is called \emph{critically
dense}\label{critical-defn} if it is infinite and each infinite
subset of this orbit is (Zariski) dense in $X$.  We say that   $Z_{\mc{I}}$  is
 \emph{saturating}\label{saturating}
  if either $Z_{\mc{I}}=\emptyset$ or,
   for each point $z\in Z_{\mc{I}}$, the $\tau$-orbit $\{\tau^i(z) : i\in\mb{Z}\}$ is
   critically dense.
 \begin{proposition}
\label{bimod alg noeth}
Keep   Assumptions~\ref{assumptions7}. Then
\begin{enumerate}
\item
The bimodule algebra $\mc{R} = \mc{R}(X, Z_{\mc{I}}, \mc{L}, \tau)$
is  noetherian if and only if $Z_{\mc{I}}$ is a saturating subset
of $X$.
\item  Assume that $\mc{L}$ is $\tau$-ample  and  that either $Z_{\mc{I}}=\emptyset$ or that
each point $z\in Z_{\mc{I}}$ lies on a Zariski dense $\tau$-orbit.
 Then  $R(X, Z_{\mc{I}}, \mc{L}, \tau)$ is
noetherian if and only if $Z_{\mc{I}}$ is saturating.
 \end{enumerate}
\end{proposition}

\begin{proof}  (1)  See \cite[Proposition~2.12]{RS}.

(2) If $Z_{\mc{I}}$ is saturating then this follows from part (1) combined with
Theorems~\ref{ampleness} and \ref{VdB main theorem}.  If
 $Z_{\mc{I}}$ is  not saturating then it follows from  \cite[Proposition~3.16]{RS}.
 \end{proof}

\section{Surjectivity in large degree}
\label{sect-surj}

The aim of this section is to show that if a cg algebra  is ``close
to'' being a \nsr, then it actually equals one (see
Theorem~\ref{surj in large degree}  for the precise statement). This
is a  generalization  of   \cite[Theorem~4.1]{AS1}, which proved the
analogous result for twisted homogeneous coordinate rings, and the
present proof closely follows the argument from \cite{AS1}.

The   following  assumptions and notation will remain in force through to the end of
the proof of Theorem~\ref{surj in large degree}.

\begin{assumptions}\label{assumptions8} Fix an integral  projective
scheme $\XXX$   of dimension $\geq 2$ with sheaf of rational
functions $\mc{K}$,  an automorphism $\tau$ and a $\tau$-ample
invertible sheaf $\LL\subset \mc{K}$.  Let $\II\subseteq
\mc{O}_{\XXX}$ be an ideal sheaf   defining a zero dimensional
subscheme $Z_{\II}$ of $\XXX$ and write $\QQ_n = \mc{L}_n
\mc{I}_n\cong \mc{L}_n\otimes \II_n$ for each $n\geq 0$. Let
$\PPP_1$ be  a finite set of closed points of $\XXX$, with $  \supp
Z_{\II}\subseteq \PPP_1$, and assume  that each
  $x \in \PPP_1$ lies on a (Zariski) dense $\tau$-orbit of $\XXX$.
   \end{assumptions}

 For $n\geq 1$ write  $\PPP_n =\bigcup_{j=0}^{n-1} \tau^{-j}(\PPP_1) $;
thus,  $\PPP_n \supseteq \bigcup_{j=0}^{n-1} \tau^{-j}(\supp
Z_{\II}) = \supp \mc{O}_\XXX/\mc{I}_n$.  We continue to use the notation from the previous section;
in particular a \nsr\ will be written
in the form  $R=R(\XXX , Z_{\mc{I}}, \mc{L},\tau)
=\bigoplus\overline{R}_nt^n\subseteq K[t,t^{-1};\tau] $, with
$\overline{R}_n=\mathrm{H}^0(\XXX, \mc{R}_n)$, as described in
\eqref{naive-defn}.

We can now state the first form of the surjectivity result. For
rings generated in degree one, this  will be further   refined in
Corollary~\ref{surj-cor}.

\begin{theorem}
\label{surj in large degree} Keep the data from
Assumptions~\ref{assumptions8}. Let
  $T =\bigoplus_{n\geq 0} \ol{T}_n t^n$
  be a cg noetherian subring of the \nsr\ $R=R(\XXX , Z_{\mc{I}}, \mc{L},\tau)$ such that
\begin{enumerate}
\item for $n\gg 0$, $ \ol{T}_n$ generates the sheaf $\QQ_n$, and
\item  for $n\gg 0$,   the birational map $\theta_n: \XXX \dra \mb{P}(\ol{T}_n^{\,*})$
defined by the sections $\ol{T}_n\subseteq \ol{R}_n\subseteq
\mathrm{H}^0(\XXX,\mc{L}_n)$ restricts to give  an immersion
 $(\XXX \smallsetminus \PPP_n)\hookrightarrow  \mb{P}(\ol{T}^{\,*}_n)$.
\end{enumerate}
Then  $T=R$ up to a finite dimensional vector space and $Z_{\mc{I}}$ is saturating.
\end{theorem}

Note that, by (1),  the sections $\ol{T}_n$ do generate
$\mc{L}_n(\XXX \smallsetminus \PPP_n)$ and so part (2) makes sense.
The proof of Theorem~\ref{surj in large degree} will be obtained
through a series of lemmas and so the assumptions and notation from
the theorem  will remain in force throughout that proof. The first
lemma is a standard application of $\tau$-ample sheaves.

\begin{lemma}
\label{ampleness conseqs}
 Let $V \subseteq \HB^0(\XXX, \FF)$ be a vector space
that generates the coherent sheaf $\FF$ and pick $r\in \mathbb{Z}$.
Then:
\begin{enumerate}
\item   For $p\gg 0$, the natural homomorphism
$\alpha: V\otimes_k \text{H}^0(\XXX ,\, \QQ_p^{\tau^r})\to
\text{H}^0(\XXX ,\, \mathcal F\otimes \QQ_p^{\tau^r})$ is
surjective.

\item For  $p \gg 0$,  the natural   homomorphism
$\beta: V^{\tau^p}\otimes_k \text{H}^0(\XXX ,\, \QQ_p)\to
\text{H}^0(\XXX ,\, \mathcal F^{\tau^p}\otimes \QQ_p)$ is
surjective.
\end{enumerate}
\end{lemma}

\begin{proof}
(1)  Set   $V' = V^{\tau^{-r}} \subseteq \HB^0(\XXX, \,
\FF^{\tau^{-r}})$.   Clearly
$\alpha$ is surjective if and only if $\alpha': V' \otimes
\HB^0(\XXX ,\, \QQ_p)\to \HB^0(\XXX ,\, \FF^{\tau^{-r}} \otimes
\QQ_p)$ is surjective. Thus we may assume that $r = 0$.

There is a short exact sequence $ 0\to\HH \to V\otimes_k \CO_{\XXX }
\to \mathcal F\to 0 $ for the appropriate  sheaf $\mc{H}$.
 Tensoring this   sequence with $\QQ_p$  gives an
exact sequence
\[
0\too \mc{G} \ \buildrel{\alpha}\over{\too}\ \HH \otimes \QQ_p \too
V\otimes_k \QQ_p
 \too  \FF \otimes \QQ_p\too 0.
\]
The sheaf   $\mc{G}$ is supported at the finite set of points
$\PPP_p$ where $\QQ_p$ need not be locally free, and so $\HB^j(\XXX
,\mc{G})=0$ for $j>0$. Therefore, if
$\mc{A}_p=\operatorname{Coker}(\alpha)$, then $\HB^j(\XXX, \mc{A}_p
) =\HB^j(\XXX ,  \HH \otimes \QQ_p)$ for all $j>0$
 and all $p \geq 0$.  As $\{\QQ_p\}$ is ample, this implies that $\HB^j(\XXX ,  \mc{A}_p )
=0 $, for all $j>0$ and $p\gg 0$.  Thus,  taking global sections of
the last displayed equation
 gives the required surjection
$V\otimes_k \HB^0(\XXX ,\,\QQ_p)
 \to \HB^0(\XXX,\, \FF \otimes \QQ_p)$.

(2) By applying $\tau^{-p}$ to the given equation,  it suffices to
prove  surjectivity of  the map
 $ V\otimes \text{H}^0(\XXX ,\, \QQ_p^{\tau^{-p}})\to
  \text{H}^0(\XXX ,\, \mathcal F\otimes \QQ_p^{\tau^{-p}})$.
Now use the proof of part (1),   combined with the second
assertion from Theorem~\ref{ampleness}.
\end{proof}

\begin{lemma}
\label{tensor versus times} The natural map $\HB^0(\XXX, \, \QQ_m
\otimes \QQ_n^{\tau^m}) \to \HB^0(\XXX, \, \QQ_{n+m})$ is a
surjection for any $m,n \geq 0$.
\end{lemma}
\begin{proof}
Since $\QQ_{m+n} = \QQ_m \QQ_n^{\tau^m}$, the natural surjective map
from the tensor product to the product fits into a short exact sequence
$0 \to \mc{H} \to \QQ_m \otimes \QQ_n^{\tau^m} \to \QQ_{m+n} \to 0$
for some sheaf $\mc{H}$.  Since $\QQ_m$ is locally free except at a finite
set of points,  $\mc{H}$ is supported on a finite set of
points and thus $\HB^1(\XXX, \, \mc{H}) = 0$.  Taking global
sections of the short exact sequence  gives the required surjection.
\end{proof}

\begin{lemma}
\label{finite extension} $R$ is  finitely generated as a left and
right  $T$-module.  Moreover, $R$ is noetherian, and $Z_{\II}$ is a
saturating subscheme of $\XXX $.
\end{lemma}

\begin{proof}
Fix $n \gg 0$ for which $\ol{T}_n$ generates $\QQ_n$. Then
Lemma~\ref{ampleness conseqs}(1), with $\mc{F} = \QQ_n$ and
$V=\ol{T}_n\subseteq \HB^0(\XXX, \, \QQ_n)$,   together with
Lemma~\ref{tensor versus times},  implies that
the multiplication map
\[
\ol{T}_n\otimes   \ol{R}_p^{\,\tau^n}=
\ol{T}_n\otimes \HB^0(\XXX ,\,\QQ_p^{\tau^n})\; \too \HB^0(\XXX, \,
\QQ_n \otimes \QQ_p^{\tau^n}) \too \HB^0(\XXX ,\,\QQ_{n+p})=
\ol{R}_{n+p}
\]
is surjective for all $p\gg 0$; say for all $p\geq p_0$. Therefore,
$R$ is  finitely generated by $R_0\oplus\cdots \oplus R_{p_0+n-1}$
as a left $T$-module. An analogous argument, using
Lemma~\ref{ampleness conseqs}(2), shows that $R_T$ is finitely
generated.

As $T$ is   noetherian, this certainly implies that $R$ is
noetherian.  By Assumption~\ref{assumptions8} the hypotheses of
 Proposition~\ref{bimod alg noeth}(2)  are therefore satisfied and, by that  result,  $Z_{\II}$
is saturating.
\end{proof}

\begin{lemma}
\label{reduce to veronese} It suffices to prove Theorem~\ref{surj in
large degree} for some pair of Veronese rings $T^{(d)}\subseteq
R^{(d)}$.
\end{lemma}

\begin{proof}
The Veronese ring $R^{(d)} \cong R(\XXX , Z_{\mc{I}_d},  \mc{L}_d,
\tau^d)$
 is again a \nsr\ with data satisfying
Assumptions~\ref{assumptions8} and so the hypotheses of
Theorem~\ref{surj in large degree} do carry over to the Veronese
rings $T^{(d)}\subseteq R^{(d)}$.

Assume that the theorem holds for $T^{(d)}\subseteq R^{(d)}$; thus
$T^{(d)}_s=R^{(d)}_s$ for all $s\gg0$, and hence  $\ol{T}_{sd}=
\ol{R}_{sd}$ for all such $s$.  By hypothesis, we can find $n_0$
such that $  \ol{T}_n$
generates $\QQ_n$ for all $n \geq n_0$. Pick
$n$ satisfying $n_0\leq n\leq n_0+d$.   Then Lemma~\ref{ampleness
conseqs}(1), with $p=sd$, $r=n$  and $V = \ol{T}_n \subseteq
\HB^0(\XXX, \, \QQ_n)$,
 together with Lemma~\ref{tensor versus times},  implies that the multiplication map
$$ \ol{T}_n \otimes  \ol{T}_{sd}^{\,\tau^n}\;=\;
  \ol{T}_n \otimes  \ol{R}_{sd}^{\,\tau^n} \;\too\; \ol{R}_{n+sd}$$
is surjective for $s\gg 0$, say for all $s\geq s_0$. Since there are only finitely many choices of
$n$, we may assume that $s_0$ is
chosen independently  of the choice of $n$.
Hence, $ \ol{T}_{u}=\ol{R}_{u}$ for all $u\geq n_0+s_0d$.
\end{proof}

The  key step in the proof of Theorem~\ref{surj in large degree}  is to relate ideals of $R$ to
$\tau$-invariant ideal sheaves on $\XXX$. This is given by the next  result.

\begin{proposition}
\label{ideals of B} Let $M=\bigoplus_{n \geq 0} \overline{M}_nt^n$
be a two-sided ideal of the \nsr\ $R$. Then  there exists a
$\tau$-invariant   sheaf of ideals $\JJ\subseteq \mc{O}_{\XXX }$
such that $  \HB^0(\XXX ,\JJ \QQ_n)=\overline{M}_n$ for all $n\gg
0$.
\end{proposition}

\begin{proof}  Let
$\RR=\RR(\XXX,Z_{\mc{I}},\mc{L},\tau)=\bigoplus_{n\geq 0}\RR_nt^n$
be the bimodule algebra associated to $R$, in the sense of
\eqref{bimodule-defn}.  By Theorem~\ref{ampleness} the sequence
$\{\RR_n\}$ is ample and by Proposition~\ref{bimod alg noeth} $\RR$
is noetherian, so Theorem~\ref{VdB main theorem}  applies. By that
result, there exists $n_0\geq 0$ and a coherent
  submodule $\GG=\bigoplus_{n\geq 0} \GG_nt^n \subseteq \RR$ such that
$\ol{M}_n=\HB^0(\XXX , \GG_n)$   for each $n\geq n_0$.
By the definition of coherent $\RR$-modules in Section~\ref{naive-section}   and increasing
$n_0$ if necessary, we may assume that $\GG_{n+1} =\GG_{n}
\QQ_{1}^{\tau^{n}}$ for all $n\geq n_0$. Since  $\RR_m = \RR_1
\cdots \RR_1^{\tau^{m-1}}$, induction implies that $\GG_{n+m}
=\GG_{n} \QQ_{m}^{\tau^{n}}$ for $m>0$ and $n\geq n_0$.
Equivalently, if $\HH_n =\GG_n \LL_n^{-1}\cong
\GG_n\otimes\LL_n^{-1}$, then
  $\HH_n \subseteq \OO_{\XXX} $ is
a sheaf of ideals satisfying
\begin{equation}\label{ideals1}
 \HH_{n+m} = \HH_n\mc{I}_m^{\tau^n} \subseteq \HH_n \subseteq \HH_{n_0}\qquad \mathrm{ for
\ each}\ n\geq n_0\ \mathrm{ and}\  m \geq 0 .
\end{equation}
 As $\GG_{m + n_0} = \GG_{n_0}\RR_{m}^{\tau^{n_0}}$  for all $m \geq 0$, the ampleness
 of $\{\RR_m\}$ implies that $\GG_n$ is generated by its sections $\overline{M}_n$ for all $n\gg 0$.

  Now consider the fact that $M$ is also a left ideal of $R$.   The hypotheses on $R$ and $\mc{R}$
are left-right symmetric and so the equivalence of categories,
Theorem~\ref{VdB main theorem}, also applies to left ideals.  Thus
there also  exists a coherent  left submodule
 $\mc{G}'=\bigoplus_{m\geq 0} \mc{G}'_mt^m\subseteq \mc{R}$
such that $\mathrm{H}^0(Y,\mc{G}_n')=\overline{M}_n$  for all $n\gg
0$. By ampleness,  again, $\mc{G}_n'$ is generated by its sections
$\overline{M}_n$ for all $n\gg 0$. Thus
$\mc{G}'_n=\mc{G}_n$ for all $n\gg 0$, which we may assume happens
for all $n\geq n_0$. As $\mc{G}' $ is a left submodule of $\mc{R}$
this implies that   $\mc{R}_m\mc{G}_n^{\tau^m}
=\mc{R}_m(\mc{G}_n')^{\tau^m} \subseteq \mc{G}'_{n+m}= \mc{G}_{n+m}$
for all $m\geq 0$ and $n\geq n_0$. Equivalently
\begin{equation}\label{ideals2}
\II_m \HH_{n}^{\tau^m}  \subseteq \HH_{n+m}  \qquad \mathrm{ for \
each}\ n\geq n_0\ \mathrm{ and}\  m \geq 0 .
\end{equation}

Set $\PPP_\infty = \{ \tau^i(\PPP_1) : i \in \mb{Z} \}$; thus
$\PPP_\infty
  \supseteq \supp\mc{O}_{\XXX}/\mc{I}_r$ for all $r\geq 0$.
  If $\mc{N}\subseteq \mc{O}_{\XXX}$ is a sheaf of ideals, write
 $\mc{S}(\NN)$ for  the largest sheaf of
ideals  with the property that $\supp \mc{S}(\NN)/\NN \subseteq
\PPP_\infty$.  Clearly $ \mc{S}(\II\NN)= \mc{S}(\NN)$ and so it follows from
\eqref{ideals1} and \eqref{ideals2}   that
$$ \mc{S}(\HH_{n_0})^{\tau} =  \mc{S}(\HH_{n_0}^{\tau})
=\mc{S}(\II\HH_{n_0}^\tau ) \subseteq   \mc{S}(H_{n_0+1}) \subseteq
 \mc{S}(\HH_{n_0}).$$  Since $\XXX $ is  noetherian, this forces
$ \mc{S}(\HH_{n_0})^{\tau} =  \mc{S}(\HH_{n_0})$. Thus $\JJ =
 \mc{S}(\HH_{n_0})$ is a $\tau$-invariant   sheaf of ideals and we will show that it satisfies the conclusions of the proposition.

We therefore need  to show that $\GG_n= \JJ \mathcal{R}_n$ in high
degree. Tensoring  this statement by $\LL_n^{-1}$, we need to show
that $\HH_n= \JJ \mc{I}_n$ for $n \gg 0$.   Since $\XXX $ is
noetherian, there exists an ideal sheaf $\AA$ such that
$\HH_{n_0}\supseteq \JJ \II_{n_0} \AA$, and  such that
$\OO_X/\mc{A}$ is supported at a finite number of points of
$\PPP_\infty$. For any $m\geq 0$, \eqref{ideals2} gives
\[
\HH_{n_0+m} \supseteq \II_m \HH_{n_0}^{\tau^m} \supseteq \II_m
\II_{n_0}^{\tau^m} \AA^{\tau^m} \JJ^{\tau^m} = \AA^{\tau^m}
\II_{n_0+m} \JJ.
\]
On the other hand,  \eqref{ideals1} implies that
\[
\HH_{n_0+m} = \HH_{n_0} \II_m^{\tau^{n_0}} \supseteq
\II_m^{\tau^{n_0}} \II_{n_0} \AA \JJ = \AA \II_{n_0+m} \JJ.
\]
Since $\OO_{\XXX} /\AA$ is supported at a finite number of points,
each lying on an  infinite $\tau$-orbit,
  $\AA^{\tau^m}+\AA=\OO_\XXX $ for $m \gg 0$.
 Combining the last two displayed equations therefore  shows that
$\HH_{n_0+m}\supseteq \JJ \II_{n_0+m},$ for all $m\gg 0$.

As $\mc{J}$ is $\tau$-invariant, $\supp \mc{O}_{\XXX}/\mc{J}$ is a
 proper $\tau$-invariant closed subspace of $X$. Since each $x \in \PPP_\infty$  lies on a dense
 $\tau$-orbit, this forces $W=\supp(\mc{O}_{\XXX}/\mc{J}) \cap \supp(\mc{O}_{\XXX}/\mc{I}_n)$
 to be empty.  But $\mc{J}\cap\mc{I}_n/\mc{J}\mc{I}_n$ is supported
on $W$.
 Thus $\mc{J}\cap\mc{I}_n=\mc{J}\mc{I}_n$ for all $n$ and so the previous paragraph implies that
 $\HH_n \supseteq \mc{J}\cap\mc{I}_n$ for all $n\gg n_0$.
 On the other hand, $\HH_n\subseteq \HH_{n_0}\subseteq \mc{S}(\mc{H}_{n_0})=\mc{J}$ by
 \eqref{ideals1}  and $\mc{H}_n\subseteq \mc{I}_n$ simply because
 $\mc{G}_n\subseteq \mc{R}_n=\mc{L}_n\mc{I}_n$. Thus  $\HH_n \subseteq  \mc{J}\cap
 \mc{I}_n =  \mc{J}\mc{I}_n$
 and hence   $\HH_n = \mc{J}\mc{I}_n$ for all $n\gg n_0$.
Equivalently,  $\GG_n = \JJ \mathcal{R}_n$  and $\overline{M}_n=
\HB^0(\XXX ,\JJ \QQ_n)$, for all $n\gg 0$.
\end{proof}

\begin{remark}\label{ideals of B2} The proof of Proposition~\ref{ideals of B} also
proves the following result: Let $\GG =\bigoplus \GG_nt^n$ be an
ideal of the bimodule algebra $\RR$. Then there exists   a
$\tau$-invariant  sheaf of ideals $\JJ\subseteq \mc{O}_{\XXX }$ such
that $ \GG_n=\mc{J}\RR_n$ for all $n\gg 0$.
\end{remark}

\noindent {\bf Proof of Theorem~\ref{surj in large degree}.}
First, by Lemma~\ref{reduce to veronese} we may replace $R$ and $T$
by
  Veronese rings and thereby assume that conditions (1) and (2)
of the theorem actually hold for all $n\geq 1$. By
Assumptions~\ref{assumptions8},   $R$ and hence $T$  are contained
in $k(\XXX)[t,t^{-1};\tau]$ and by condition (1),  $T_1\not=0$; say
$0\not=t_1=at\in T_1$ for some $a\in k(X)$. Therefore,
$k(\XXX)[t,t^{-1};\tau] =k(\XXX)[t_1,t_1^{-1};\tau_1]$ and one even
has $\tau=\tau_1$ since they differ by conjugation by the element
$a\in k(\XXX)$.
  Since $T$ is noetherian, it   has a graded Goldie quotient ring which, by universality,
  will  have the form  $Q(T)=F[t_1,t_1^{-1};\tau]
\subseteq Q(R)\subseteq k(\XXX)[t_1,t_1^{-1};\tau]$  for some field
$F\subseteq k(\XXX)$ invariant under $\tau$.

 We claim that $F=k(\XXX)$. To see this, set
 $V=\XXX\smallsetminus \PPP_1$ and note that, by  conditions (1) and (2),
   $\ol{T}_1$ generates $\mc{L}\vert_{V}$ and induces an immersion
$\chi: V \hookrightarrow \mb{P}(\ol{T}_1^{\, *})$.  Fix a basis $\{s_j : 0\leq j\leq r\}$ of
$\ol{T}_1$, write $V_0=\{x\in V : (s_0)_x\notin
\mathfrak{m}_x\mc{L}_x\}$ and let $U_0=\{x_0\not=0\}$ be the
standard open affine subset of $\mb{P}(\ol{T}_1^{\, *})$ with
$\mc{O}(U_0)=k[x_1x_0^{-1},\dots,x_rx_0^{-1}].$ Then $\chi(V_0)$ is
open in its closure inside $U_0$.    As in the proof of
\cite[Theorem~II.7.1]{Ha}, the morphism  $\chi \vert_{_{\scriptstyle V_0}}$
  is then given by the map of rings  $\alpha: k[x_1x_0^{-1},\dots,x_rx_0^{-1}] \to
\Gamma(V_0,\mc{O}_{V_0})$ defined by $x_jx^{-1}_0\mapsto s_j
s^{-1}_0$, and it follows that $\Gamma(V_0,\mc{O}_{V_0})$ is a
localization of the image of $\alpha$.
On the other hand, as $X$ is an integral
scheme, it has the same field of rational functions as both $V$ and
$V_0$   and so $k(X)$ is
generated as a field by the $ s_js^{-1}_0\in F$. Thus $k(X)=F$.

It follows from the previous paragraph that   $Q(T) = Q(R) = k(X)[t_1,t_1^{-1};\tau]$ and so $R/T$
is Goldie torsion as a right (or left)  $T$-module.
 By  Lemma~\ref{finite extension}, $R/T$ is therefore a finitely
generated, torsion right $T$-module and so the left annihilator
$J_\ell=\lann_T(R/T)$ is a non-zero ideal of $T$. Similarly,
$J_r=\rann_T(R/T)$ and $M=J_r J_\ell $ are  non-zero ideals of $T$.
 But $M$ is also an ideal of $R$. So, if $M=\bigoplus M_n$ then
   Proposition~\ref{ideals of B} produces  a  non-zero
$\tau$-invariant ideal sheaf $\JJ\subseteq \mc{O}_\XXX $ such that
\begin{equation}\label{factor05}
 \HB^0(\XXX , \JJ \QQ_n) t^n =  M_n
 \qquad\mathrm{ for\  all}\quad n \gg 0.
 \end{equation}

 We next  consider the factor bimodule algebra  $\QQ/\mc{J}\QQ =
\bigoplus_{n\geq 0} \QQ_nt^n/\mc{J}\QQ t^n$. Let $W \subseteq \XXX $
be the $\tau$-invariant subscheme defined by $\JJ$. As every point
of $\supp Z_{\II_n}\subseteq \PPP_n$ lies on a dense $\tau$-orbit,
and $W$ is a proper $\tau$-invariant   subscheme of $\XXX $, it
follows that $Z_{\II_n}\cap W = \emptyset$ and hence that
$\QQ_n/\mc{J}\QQ_n = \QQ_n\otimes_{\mc{O}_{\XXX}} \mc{O}_{W}=
\mc{L}_n\otimes_{\mc{O}_{\XXX}}\mc{O}_{W}$
 for all $n\geq 0$.  In other words
 $\QQ/\mc{J}\QQ \cong \bigoplus_{n\geq 0} (\mc{L}_n\vert_{W})\ol{t}^{\,n}$,
 where $\ol{t}$ is now a formal symbol denoting the image of $t$.
 If we write $\ol{\mc{L}}=\mc{L}\vert_W$, then
  $\mc{L}_n\vert_W \cong \ol{\mc{L}}
 \otimes_{\mc{O}_W}\ol{\mc{L}}^{\,\tau}\otimes\cdots\otimes\ol{\mc{L}}^{\,\tau^{n-1}}=\ol{\mc{L}}_n$
and so    $\QQ/\mc{J}\QQ$  is just the bimodule algebra $\mc{B}(W, \LL \vert_{W}, \tau \vert_{W}) $,
 as  defined in Section~\ref{naive-section}.  Taking global sections we obtain
 $$\mathrm{H}^0(\XXX, \QQ/\mc{J}\QQ)\ = \
  \bigoplus_{n \geq 0} \HB^0(\XXX ,\, \mc{L}_n\otimes_{\mc{O}_\XXX}\mc{O}_W)\ol{t}^{\,n}
\ = \  \bigoplus_{n \geq 0} \HB^0(W ,\,  \mc{L}_n\vert_{W})\ol{t}^{\,n}
\ = \ B(W, \LL \vert_{W}, \tau \vert_{W}) ,$$
where $B=B(W, \LL \vert_{W}, \tau \vert_{W}) $ is the
 ordinary twisted homogeneous coordinate ring.

 Now take global sections of the short exact sequence
 of $\OO_\XXX $-modules $ 0 \to \JJ \RR \to \RR \to\RR/\mc{J}\QQ \to 0$ to obtain
 the  exact sequence of finitely generated
 right $R$-modules
$ 0\to N  \to R \to B,$
where $N=\mathrm{H}^0(\XXX ,\JJ\RR)$. Set
\begin{equation}\label{factor}
T' \ = \ \frac{\displaystyle T+N}{\displaystyle N}
  \ = \  \bigoplus_{n\geq 0}\ol{T}'_n\ol{t}^{\,n}\  \subseteq\  R'  \ = \
  \frac{\displaystyle R}{\displaystyle N}  \ \subseteq\
B  \ = \ B(W, \LL \vert_{W}, \tau \vert_{W}).
\end{equation}

We claim that the hypotheses of the  theorem pass to the pair of
rings $T' \subseteq B$ or,  more precisely, that $T'\subseteq B$
satisfy the analogous
 hypotheses of \cite[Theorem~4.1]{AS1}.
To see this, observe first that  $\ol{T}'_n$ is the restriction to
$W$ of the global  sections contained in $\ol{T}_n\subseteq
\overline{R}_n=\mathrm{H}^0(\XXX ,\QQ_n)$. Therefore, by
condition~(1) of the theorem, $\ol{T}'_n$ generates the sheaf
$\QQ_n\vert_W =
\mc{L}_n\vert_W=\ol{\mc{L}}_n$ for all $n \gg 0$. Moreover,
the rational  map $\theta'_n: W \dra \mb{P}({\ol{T}'_n}^*)$ is  the
restriction of the rational map   $\theta_n : \XXX  \dra \mb{P}(\ol{T}_n^{\,*})$.
Once again, as $W$ is $\tau$-invariant, $W \cap
\PPP_n = \emptyset$ for all $n$, and so, by condition~(2),
  $\theta'_n$ is an everywhere-defined  immersion.
  Since $W$ is proper, its image is closed and so
$\theta_n'$ is a closed immersion.   Since $\mc{L}$ is $\tau$-ample,
it follows from
  \cite[Proposition~2.3(2) and Corollary~5.1]{Ke1}
 that $\mc{L} \vert_{W}$ is both $\tau
\vert_{W}$-ample and $(\tau\vert_{W})^{-1}$-ample as an invertible
sheaf over $W$. Moreover, by \cite[Theorem~1.2(3)]{Ke1}, $B$ has
finite Gelfand-Kirillov dimension. The hypotheses of
\cite[Theorem~4.1]{AS1} are therefore satisfied and, by that result,
$B/T'$ is finite dimensional. But \eqref{factor05} implies that
  $N_n=\mathrm{H}^0(\XXX ,\JJ\RR)_n\subseteq T_n$ for all $n\gg 0$ and so
 \eqref{factor} implies that $ R/T $ is also finite-dimensional.

Finally, $Z_{\mc{I}}$ is saturating by Lemma~\ref{finite extension}.
 \qed

\medskip
In applications of   Theorem~\ref{surj in
large degree}  it is inconvenient to have to check conditions (1)
and (2)  of that result for infinitely many choices of  $n$.
Fortunately, for rings generated in degree one,  just one choice
will do. In this final result we do not assume that Assumptions~\ref{assumptions8} are satisfied.

 \begin{corollary}\label{surj-cor}
Let $\XXX $ be an integral projective scheme, with $\dim \XXX\geq
2$, sheaf of rational functions $\mc{K}$  and with a $\tau$-ample
invertible sheaf $\mc{L}\subset \mc{K}$   for some $\tau \in
\aut(\XXX )$.
  Let $T  = \bigoplus \ol{T}_m t^m$
be a noetherian cg subalgebra of
   $B(\XXX , \mc{L}, \tau) = \bigoplus \HB^0(\XXX , \mc{L}_m)
t^m \subseteq k(\XXX)[t,t^{-1};\tau]$   that  is generated in degree
$1$.   Let $\PPP$ be a finite set of closed points in $\XXX $
such that each $x\in \PPP $ lies on a dense $\tau$-orbit. Finally,
suppose that there exists $n\geq 1$ such that:
\begin{enumerate}
\item[(i)]   $\ol{T}_n$ generates the sheaf $\mc{L}_n$, except, perhaps, at points  of $\PPP$.
\item[(ii)] The  birational map $\XXX
 \dra \mb{P}(\ol{T}_n^{\,*})$ defined by $\ol{T}_n$ restricts to give
 an immersion $(\XXX \smallsetminus \PPP)
 \hookrightarrow \mb{P}(\ol{T}_n^{\,*})$.
\end{enumerate}
Set $\mc{R}_1 = \ol{T}_1 \mc{O}_{\XXX}  \subseteq \mc{L}$
and let $\mc{I}= \mc{R}_1 \mc{L}^{-1}$.
Then $Z_{\mc{I}}  $ is a saturating zero-dimensional
subscheme of $\XXX$
with $\supp  Z_{\mc{I}}\subseteq \PPP$. Moreover $T\subseteq R=R(\XXX , Z_{\mc{I}},
\mc{L}, \tau)$ with $\dim_k T/R < \infty$.
\end{corollary}

\begin{proof}  By definition, $\ol{T}_1$ generates the
sheaf $\mc{R}_1 = \mc{I} \mc{L}$.  Since $T$ is generated in degree
$1$, the space
$\ol{T}_m=\ol{T}_1\ol{T}_1^{\,\tau}\cdots\ol{T}_1^{\,\tau^{m-1}}$
generates $\QQ_m = \QQ_1 \QQ_1^{\tau} \dots \QQ_1^{\tau^{m-1}}
=\II_m\mc{L}_m$ for all $m\geq 1$.  Thus we have the following
inclusion of subrings of $k(\XXX)[t,t^{-1};\tau]$:
\begin{equation}\label{maincor-equ}
T  = \bigoplus_{m \geq 0} \ol{T}_m t^m \subseteq \bigoplus_{m \geq
0} \HB^0(\XXX , \QQ_m) t^m = R(\XXX , Z_{\mc{I}}, \mc{L}, \tau).
\end{equation}

It is convenient to replace the integer $n$ by $1$ in the
statement of the corollary, which we will accomplish by passing to
 Veronese rings.
Now  the hypotheses of the corollary clearly pass
across to the rings $T^{(n)}\subseteq   B(\XXX , \LL_n,
\tau^n)$   (in which case take $n=1$ in the statement of the result).
Assume for the moment  that the corollary holds for this pair of
rings; thus $\supp Z_{\mc{I}_n}\subseteq P$
 and  $T^{(n)}$ is equal in large degree to $R(\XXX,
Z_{\mc{I}_n}, \mc{L}_n, \tau^n)$ or, in other words,
$\overline{T}_{jn}=\mathrm{H}^0(\XXX, \QQ_{jn})$ for all $j\gg 0$.
This ensures  that Assumptions~\ref{assumptions8} also hold for these Veronese rings,
 and we can apply Lemma~\ref{ampleness conseqs} with
 $V=\ol{T}_a\subseteq \HB^0(\XXX, \mc{F})$ where $\mc{F} = \QQ_a$ and $r = a$ for
 any $1\leq a\leq n$.  Together with Lemma~\ref{tensor versus times},
 this  shows that the multiplication map
 $$\ol{T}_a\otimes \ol{T}_{jn}^{\,\tau^a}
 =\ol{T}_a\otimes   \mathrm{H}^0(\XXX, \QQ_{jn}^{\tau^a})
\too   \HB^0(\XXX, \QQ_a \otimes \QQ_{jn}^{\tau^a}) \too
\mathrm{H}^0(\XXX, \QQ_{jn+a})$$ is surjective for all $j\gg 0$.
Thus $\ol{T}_m = \mathrm{H}^0(\XXX, \QQ_m)$ for $m\gg 0$.  Moreover,
since we are assuming that $Z_{\mc{I}_n}$ is
saturating and supported along $P$,   clearly the same holds for
$Z_{\mc{I}}$.  Therefore, if the corollary holds for  $T^{(n)}\subseteq   B(\XXX , \LL_n,
\tau^n)$, it  also holds for $T\subseteq B(\XXX , \LL, \tau)$.

Thus we have reduced to the case $n = 1$ and we aim to show that the
hypotheses (1) and (2) of Theorem~\ref{surj in large degree} are
satisfied for all $m \geq 1$.   We have already seen
 that $\ol{T}_m$ generates $\mc{R}_m$ for all $m \geq 1$, so
hypothesis (1) of the theorem holds.  Set $\PPP_1=\PPP$ and, for
 $m\geq 1$, define $\PPP_m = \bigcup_{i = 0}^{m-1}
\tau^{-i}(\PPP_1)\supseteq \supp \mc{O}_{\XXX}/\mc{I}_m$.  Note that
$\ol{T}_1$ generates $\mc{L}$ along $\XXX \smallsetminus P_1$.
For any $j \in \mb{Z}$, condition (ii) implies that the space
$\ol{T}(j)=\ol{T}_1^{\,\tau^j}$ defines an immersion
 $\XXX \smallsetminus \PPP_n \hookrightarrow
 \XXX \smallsetminus  \tau^{-j}(\PPP_1)\hookrightarrow \mb{P}(\ol{T}(j)^{\,*})$,
and so  an easy induction  using  Lemma~\ref{immersion lemma} shows that
$\ol{T}_m$ defines an immersion
 $  \bigl(\XXX \smallsetminus \PPP_m\bigr)\hookrightarrow \mb{P}(\ol{T}_m^{\,*})$
 for all $m\geq 1$.  Thus hypothesis (2), and hence all the hypotheses,
  of Theorem~\ref{surj in large
 degree} are verified  and the corollary follows.
\end{proof}


\section{The   main theorem}\label{main-section}

We can now put everything together to prove the main theorem of the
paper, which is an elaboration of Theorem~\ref{intro-mainthm}.

\begin{theorem}\label{mainthm}
 Let $A=\bigoplus_{i\geq 0}A_i$   be a cg noetherian  domain  that is generated in degree one
 and  \btwog.  Then
 \begin{enumerate}
 \item Up to a finite dimensional vector space, $A$ is isomorphic to a \nsr\
 $R(\wtY,Z_{\mc{I}},\mc{L},\wt{\sigma})$, where $ Z_{\mc{I}}$   is either  
 empty or a saturating,
  zero dimensional subscheme of $\wtY$.
   \item The data $\{\wtY, Z_{\mc{I}},\mc{L},\wt{\sigma}\}$ can be made more precise
 as follows:
 \begin{enumerate}
 \item[(i)] For any stable value of $n$, the scheme $ \wtY=\wtY_n$ is the blow down
  $\pi_n(\YY_n)$ of  the relevant component $ \YY_n$ of the truncated point scheme, as
 constructed by Theorem~\ref{blowing down}.
 \item[(ii)] The  automorphism $\wt{\sigma}=\pi_n\sigma\pi^{-1}_n$ of $\wtY$ is
  induced from the birational map $\sigma_n=\psi_n\phi^{-1}_n$ of~$\YY_n$.
 \item[(iii)] If one writes $A=\bigoplus \ol{A}_nt^n\subset Q(A)=k(\wtY)[t,t^{-1},\sigma]$
  for some $0\not=t\in A_1$, then  $\mc{L}=(\ol{A}_1\mc{O}_{\wtY})^{**}$ with
 $\ol{A}_1\mc{O}_{\wtY}=\mc{I}\mc{L}$.
  \end{enumerate}
 \end{enumerate}
 \end{theorem}

Before proving the theorem we need some subsidiary results and we
begin by setting up the notation that will remain in force
throughout the section.
 Fix a stable value of $n$, let $\YY_n$ denote the relevant component of $\XX_n$
  and drop the  subscript from $\YY=\YY_n$, etc.
 Let $\pi:\YY\to \wtY$ be the birational surjective morphism defined by Theorem~\ref{blowing down}.
 Recall from that result that  $ \wt{\sigma} = \pi
\sigma \pi^{-1}$ is  a biregular automorphism of $\wtY$, and that the set
of points where $\pi^{-1}$ is not defined is a finite set of
nonsingular points $P$, each lying on a dense
 $\wt{\sigma}$-orbit.
 Set $V=\wtY\smallsetminus P$ and $U=\pi^{-1}(V)\subseteq
\YY$; thus $U\cong V$. As in the statement of the  theorem, fix $t\in A_1\smallsetminus\{0\} $
and write $A_i=\ol{A}_it^i$ for $i\geq 0$; thus each $\ol{A}_i\subset K=Q(A)_0$ and
$A=\bigoplus
\ol{A}_nt^n\subset Q(A)=K[t,t^{-1} ; \sigma]$. We  recall that
$K\cong k(\YY)=k(\wtY)$ by Corollary~\ref{zero-ideal} and we  always write the
corresponding sheaf of rational functions (on $\YY$ or $\wtY$)
as~$\mc{K}$.

 We begin with some elementary observations about  invertible sheaves.

\begin{lemma}\label{pushforward}
Let $\mc{A} = \sum_{i=1}^a s_i\mc{O}_\YY\subset \mc{K}$ be
a  sheaf
 generated by   sections $s_j \in \mathrm{H}^0(\YY,\mc{A})$
 and set $\mc{B}=\pi_*\mc{A}\subset \mc{K}$ with reflexive hull  $\mc{P}=\mc{B}^{**}$.
 If $\mc{A}(u)=\sum_{i=1}^a s_i^{\sigma^{u-1}} \mc{O}_{\YY}$ is
  invertible  for   $1\leq u\leq v$, then
 $\mc{B}(u) =\pi_*\mc{A}(u)$ has reflexive hull $\mc{P}(u)=\mc{P}^{\wt{\sigma}^{u-1}}$
 for $1\leq u\leq v$. Moreover, $\mc{P}(1)\cdots \mc{P}(v)$ is the reflexive hull of
 $\pi_*\big( \mc{A}(1)\cdots \mc{A}(v)\big).$
\end{lemma}

\begin{proof}
As $\wt{\sigma}$ and $\sigma$ induce the same automorphism of   $K$, we have
   $s_i^{\wt{\sigma}^v}=s_i^{\sigma^v}$ for all $i,v$.  Let  $\tau=\wt{\sigma}^{u-1}$
   for some $1\leq u\leq v$.
   Since $\mc{P}\vert_{V} =\sum s_i\mc{O}_V$
we have    $\mc{P}^{\tau}\vert_{{\tau^{-1}(V)} } =
    \sum s_i^{\tau}  \mc{O}_{\tau^{-1}(V)}$.
Set $V'  =V\cap \tau^{-1}(V)$ and note that $\wtY\smallsetminus V'$
consists of a finite set of closed nonsingular  points of $\YY$.  Since
   $$\mc{B}(u)\vert_{V'} =  \pi_*\left( \mc{A}(u)\vert_{\pi^{-1}(V')}\right)
   =  \sum_{i=1}^a s_i^{\tau} \mc{O}_{V'} =
   \mc{P}^{\tau}\vert_{V'}, $$
Sublemma~\ref{pushforward1}(2) implies that $ \mc{P}^{\tau}=\mc{B}(u)^{**}$.
The proof of the   last assertion is similar.
\end{proof}

\begin{lemma}\label{sigmatilde-ample}{\rm (1)} The scheme $\wtY$ has a
$\wt{\sigma}$-ample invertible sheaf.

\begin{enumerate}
\item[(2)] Let $\mc{P}$ be an invertible sheaf on $\wtY$,
and write $\mc{P}_\ell = \mc{P} \otimes \mc{P}^{\wt{\sigma}}
\otimes \dots \otimes \mc{P}^{\wt{\sigma}^{\ell-1}}$ for $\ell \geq 1$.
  If $  \mc{P}_n\otimes \mc{P}_n^{\wt{\sigma}^q}$ is ample for some $n \geq 1$
and all  $q \gg  0$  then $\mc{P}$ is $\wt{\sigma}$-ample.
\end{enumerate}
\end{lemma}

 \begin{proof}  (1)   As $A$ is noetherian, \cite[Theorem~0.1]{SteZh}
implies that $A$ has sub-exponential growth.  Since $Q(A) = K[t,
t^{-1}; \sigma]$, \cite[Proposition~3.5 and Lemma~2.3(2)]{RZ}
together with \cite[Theorem~1.2(2)]{Ke1} imply that $\wtY$ has a
$\wt{\sigma}$-ample invertible sheaf.

(2)  By part (1) and \cite[Theorem~1.3]{Ke1}, $\wt{\sigma}$ is quasi-unipotent in the sense
of \cite{Ke1}.   By hypothesis, $\mc{F}=\mc{P}_n\otimes \mc{P}_n^{\wt{\sigma}^{mn}}$
 is ample for some $m\geq 1$.  Hence, so is
$\mc{G}=\mc{F}\otimes \mc{F}^{\wt{\sigma}^n}\otimes\cdots\otimes \mc{F}^{\wt{\sigma}^{nm-n}}.$
But expanding terms, one finds that  $\mc{G}=\mc{P}_{2nm}$.
By \cite[Theorem~1.3]{Ke1}, $\mc{P}$ is therefore $\wt{\sigma}$-ample.
\end{proof}

\noindent {\bf Proof of Theorem~\ref{mainthm}.}
We continue to use  the notation set up immediately after the statement of the theorem.
The aim of the proof is to apply  Corollary~\ref{surj-cor}   and this
 will be easy  to do   once we understand the relationship between
sheaves of $\mc{O}_\YY$-modules generated by the $\ol{A}_j$ and those induced from the
embedding of $Y$ into $\mb{P}(A_1^*)$. Thus to begin,
 we want to be very precise about the maps of
schemes involved, after which we will translate the resulting identifications into
ones involving sheaves.

Fix a basis $\{x_j = a_jt :  0\leq j\leq r\}$ of $A_1$,
where each $a_j \in \ol{A}_1 \subset K$.  Write
$\mathbb{P}=\mathbb{P}(A_1^*)$  and let
$\iota:\YY\hookrightarrow \XX\hookrightarrow \mb{P}^{\times
n}$ denote  the restriction to $\YY$ of the natural inclusion described
in Section~\ref{sect-pts}.
For $1\leq m \leq n$,  let $\rho_m :
\mathbb{P}^{\times n}\to \mathbb{P}$ be the projection onto the
$m^{\text{th}}$ copy of $\mathbb{P}$. As usual, we regard
$Q(A)_{\geq 0}$ as a ``generic'' right $K$-point module, where $K$
acts by left multiplication.  By Corollary~\ref{zero-ideal}, the
truncated $K$-point module $\bigoplus_{i = 0}^n Q(A)_i$ corresponds
to a map $\theta: \eta = \spec K \to \YY$
and hence, for $1\leq m\leq n$, induces a map
\begin{equation}\label{mainthm-maps}
\chi_m :   \eta \  \buildrel{\theta} \over\longrightarrow \
 \YY\   \buildrel{\iota} \over\longrightarrow \  \mathbb{P}^{\times n}
\  \buildrel{\rho_m} \over\longrightarrow \  \mathbb{P}.
  \end{equation}
   As  in Notation~\ref{S-notation},
we identify $K=k(Y)$ through  $\theta$.

We need the precise form of $\theta$, so we recall its
construction  from  \cite{ATV1}.    Write
$ \underline{a} = (a_0:a_1: \dots: a_r)$,  regarded as a $K$-valued
point of $\mathbb{P}$.   The decomposition  $\bigoplus_{i = 0}^n Q(A)_i =
\bigoplus_{i = 0}^n K t^i$   fixes the  generator $t^i$
of the free $K$-module $Q(A)_i$   and the right action of
$A_1$   is then given by the formula $t^i x_j=a_j^{\sigma^i} t^{i+1}$ for all $i,j\geq 0$.
  Thus, by the proof of
 \cite[Proposition~3.9]{ATV1} (see \cite[Equations~3.10 and 3.11]{ATV1} in particular), the
 map $\iota \theta$ associates to  $\eta$ the $K$-valued point
 $(\underline{a},\, \underline{a}^\sigma,\dots,\, \underline{a}^{\sigma^{n-1}})$ of
$  \mathbb{P}^{\times n}$.
 Equivalently,  for $1\leq m \leq
n$,   $\chi_m$ associates $\eta$ to the $K$-valued point $\underline{a}^{\sigma^{m-1}}$
of $\mb{P}$.   Alternatively, regard the
 $\{x_j\}$   as coordinate functions on $\mathbb{P}$ and  set
$U_j=\{p\in \mathbb{P} : x_j(p)\not=0\}$ for $0\leq j\leq r$.
Then, for $0\leq j\leq r$ and $1\leq \ell\leq n$,
the map of rings $\mc{O}(U_j)=k[x_1x_j^{-1},\dots, x_rx_j^{-1}] \to K$ induced by
$\chi_\ell$ is simply
 defined by  $x_ix_j^{-1} \mapsto s_i(\ell)s_j(\ell)^{-1}$
for $s_m(\ell)=    a_m^{\sigma^{\ell-1}}.$

We now use this information to construct  invertible sheaves on  $\YY$,
 which we do by following the proof of  \cite[Theorem~II.7.1]{Ha}.
  Fix $1\leq m \leq n$ and let $P(m)=K$,
regarded as a locally free sheaf on $\eta=\spec K$ generated by the
  sections  $\{ s_i(m)  : 0\leq i\leq
r\}$. By the proof of  \cite[Theorem~II.7.1(b)]{Ha},  this data defines a morphism
 $\chi_m': \eta\to \mathbb{P}$ given  by $x_ix_j^{-1}\mapsto s_i(m)s_j(m)^{-1}$
 as a homomorphism  $\mc{O}(U_j)\to K$;
  in other words it defines the morphism $\chi_m$. Conversely, by
  \cite[Theorem~II.7.1(a)]{Ha},  we can identify
  $P(m)=\chi_m^*\mathcal{O}_{\mathbb{P}}(1)$ via  the explicit formula
  $s_i(m)=\chi_m^*(x_i)$ for $0\leq i\leq r$.

 Now    apply  \cite[Theorem~II.7.1(a)]{Ha}  to the map   $\widetilde{\chi}_m^{\phantom{*}} =
\rho_m^{\phantom{*}}\iota : Y\to \mb{P}$ for some $1\leq m\leq n$.
   Since    $\mc{O}_{Y}$ is embedded into
    $ \mc{K}$ via   $\theta$,   we   have
$  \wt{\chi}^*_m(x_i)=\chi^*_m(x_i)=s_i(m)$ for each $i$.
 Therefore, \cite[Theorem~II.7.1(a)]{Ha} shows that   the elements
 $\{ s_i(m) : 0\leq i\leq r\}$
 generate a locally free sheaf of $\mc{O}_\YY$-modules
 $\mathcal{N}(m) \cong \wt{\chi}_m^*(\mc{O}_{\mb{P}}(1))
  =\iota^*\rho^*_m (\mc{O}_{\mb{P}}(1))$.

To summarize,  as   $\{a_i^{\sigma^{m-1}} : 0\leq i\leq r\}$ is a $k$-basis of
  $\overline{A}_1^{\,\sigma^{m-1}}$, we have shown:
\begin{equation}\label{mainthm-eq1}
  \overline{A}_1^{\,\sigma^{m-1}} \subseteq  \mathrm{H}^0(\YY,\,
  \mc{N}(m)) \subseteq K
\qquad\mathrm{and}\qquad \mc{N}(m)  =
\overline{A}_1^{\,\sigma^{\,m-1}} \mc{O}_\YY \subseteq \mc{K}
\qquad\mathrm{for}\ 1\leq m\leq n.
\end{equation}

 Now consider the invertible
sheaf of $\mb{P}^{\times n}$-modules  $ \mc{O}(1,\dots,1) =  \rho_{1}^* \mc{O}_{\mb{P}}(1)
\otimes \rho_{2}^* \mc{O}_{\mb{P}}(1) \otimes \cdots\otimes
\rho_{n}^* \mc{O}_{\mb{P}}(1) ,$ where the tensor products are over
$\mc{O}_{\mb{P}^{\times n}}$.
  Since pull-backs commute with tensor products, and the tensor product of the
 $\mathcal{N}(m)$ is isomorphic to  their product  inside $\mc{K}$,  it follows that
 $ \mc{M}=  \mathcal{N}(1) \mathcal{N}(2)\cdots \mathcal{N}(n)$ satisfies
$$\mc{M} \ \cong \   \iota^*\rho^*_1 \mc{O}_{\mb{P}}(1)\otimes\cdots\otimes
  \iota^*\rho^*_n \mc{O}_{\mb{P}}(1)
 \ = \  \iota^* \mc{O}(1,\dots,1) .$$
Therefore,    by \eqref{mainthm-eq1},
\begin{equation}\label{mainthm-eq2}
\overline{A}_n   =
  \overline{A}_1\overline{A}_1^{\,\sigma}\cdots \overline{A}_1^{\,\sigma^{n-1}}
 \subseteq  \  \mathrm{H}^0\bigl(\YY,\,    \mathcal{N}(1)\bigr)
 \cdots   \mathrm{H}^0\bigl(\YY,\,    \mathcal{N}(n)\bigr) \
\subseteq \  \mathrm{H}^0\left(\YY,\,  \mathcal{M}\right), \
\quad\mathrm{where}\quad \   \mc{M} = \ol{A}_n\mc{O}_{\YY}.
 \end{equation}
Since
  $\mathrm{H}^0\bigl(\mb{P}^{\times n},\,  \mc{O}(1,\dots,1)\bigr)\ = \
\mathrm{H}^0\bigl(\mb{P}^{\times n}, \,
\rho^*_1\mc{O}_{\mb{P}}(1)\bigr) \otimes_k \cdots \otimes_k
\mathrm{H}^0\bigl(\mb{P}^{\times n},\,
\rho^*_n\mc{O}_{\mb{P}}(1)\bigr),$
it follows that $\overline{A}_n$ is just the image of the  natural  map
$$
   \iota^* : \mathrm{H}^0(\mb{P}^{\times n},\,  \mc{O}(1,1,\dots,1))
 \to \mathrm{H}^0(\YY, \iota^*( \mc{O}(1,1,\dots,1)) \
=
   \mathrm{H}^0(\YY,\,   \mathcal{M}).
$$

On the other hand,   $\mc{O}(1,1,\dots,1)$
 is a very ample invertible sheaf over $\mb{P}^{\times n}$ and so its set of  global sections
  $ N=\mathrm{H}^0\bigl(\mb{P}^{\times n},\,  \mc{O}(1,\dots,1)\bigr)$
  defines a closed  immersion   $\xi: \mb{P}^{\times n} \hookrightarrow
  \mb{P}(N^*)$.
    By restriction,  $\xi$ defines a closed immersion of    $\YY$ into
    projective space;  to be precise,   $\overline{N}=\iota^*(N)$
    defines a closed immersion of   $\YY$ into   $\mb{P}(\overline{N}^{\,*})\hookrightarrow
    \mb{P}(N^*)$.     Comparing this with  the conclusion of the previous paragraph
     we deduce:
 \begin{equation}\label{mainthm-eq4}
  \mc{M}\ \mathrm{is\  invertible\ and \    the\ sections\  }  \overline{A}_n  \subseteq
 \mathrm{H}^0(\YY,\, \mc{M})    {\mathrm{ \  define\ a\ closed \  immersion\  }}  \YY\hookrightarrow
\mb{P}(\overline{A}_n^{\,*}).
 \end{equation}

We now push everything forward to $\wtY=\wtY_n$ via the map $\pi :
\YY\to \wtY$ constructed in Theorem~\ref{blowing down}. Recall that
 $V=\wtY\smallsetminus P$ and that  $\pi^{-1} : V\to U=\pi^{-1}(V)$ is an isomorphism.
For each $m$, write $\mc{F}(m)=\pi_*\mc{N}(m)$, with reflexive hull
$\mc{L}(m)=\mc{F}(m)^{**}$.
 By Sublemma~\ref{pushforward1}(1),   each $\mc{L}(m)$ is an  invertible sheaf of
 $\mc{O}_\wtY$-modules   contained in
 $ \mc{K}$.  Moreover,  Lemma~\ref{pushforward}  and \eqref{mainthm-eq1}
imply that  $\mc{L}(m)=\mc{L}(1)^{\wt{\sigma}^{m-1}}$ for $1\leq
m\leq n$ and  that $\mc{L}(1)\mc{L}(2)\cdots \mc{L}(n)=( \pi_*\mc{M})^{**}$.   In other
words,
$ ( \pi_*\mc{M})^{**}=\mc{L}_n=\mc{L}\mc{L}^{\wt{\sigma}}\cdots\mc{L}^{\wt{\sigma}^{n-1}}$
 for $\mc{L}=\mc{L}(1)\subset \mc{K}$.

 By \eqref{mainthm-eq4},  $\mc{M}$ is a very ample invertible sheaf
on $\YY$  while $\mc{L}_n=(\pi_*\mc{M})^{**}$ by definition. Thus,
Remark~\ref{projective-remark} implies that $\mc{L}_n \otimes
(\mc{L}_n)^{\wt{\sigma}^q}$ is   ample   for all  $q \gg 0$. Hence, by
Lemma~\ref{sigmatilde-ample}(2) and the previous paragraph, $\mc{L}$ is  a
$\wt{\sigma}$-ample invertible sheaf.

We are now ready to check the hypotheses of
Corollary~\ref{surj-cor}.
By  its construction in Theorem~\ref{blowing down},
 $\wt{\sigma}$ induces the  automorphism $\sigma$ of $K$. Since
$\mc{L}\subset \mc{K}$ this implies that $B(\wtY,\mc{L},\wt{\sigma})
=\bigoplus_{n\geq 0} \mathrm{H}^0(\wtY,\, \mc{L}_n)t^n \subset
K[t;\sigma] .$   As $A$ is generated in degree one by
$A_1=\overline{A}_1t\subseteq \mathrm{H}^0(\wtY,\, \mc{L})t$,  it
follows  that $A$ is a subring of $ B(\wtY,\,\mc{L},\wt{\sigma})$,
and hence   $\ol{A}_n t^n \subseteq \mathrm{H}^0(\wtY,\,
\mc{L}_n)t^n.$   Since   $\mc{L}$ is  a $\wt{\sigma}$-ample invertible sheaf by the previous paragraph,
and  $P$   consists of a finite set of points lying on dense $\wt{\sigma}$-orbits by
Theorem~\ref{blowing down},   the  hypotheses from  the first three sentences of
Corollary~\ref{surj-cor} are satisfied in  the present setup.
 By    \eqref{mainthm-eq2},
  $\mc{L}_n\vert_V=\pi_*\mc{M}\vert_V= \mc{M}\vert_U$  is generated by $\overline{A}_n$
  and so    condition  (1) of the corollary holds, while   condition (2) of the
  corollary  follows from \eqref{mainthm-eq4}. We can
   therefore apply the corollary to find  a saturating zero-dimensional
   subscheme $Z_{\mc{I}}\subset \wtY$
 such that $A=R(\wtY,Z_{\mc{I}}, \mc{L},\wt{\sigma})$, up to a finite dimensional vector space.
 This proves part (1) of the theorem.

 Finally,   parts (2.i) and (2.ii)   of the theorem both follow from
 Theorem~\ref{blowing down} together with our construction of $\wtY$.
   Corollary~\ref{surj-cor} also ensures that $\mc{I}$ is defined by $ \mc{L}\mc{I} =
\overline{A}_1 \mc{O}_{\wtY}$, while \eqref{mainthm-eq1} shows that
$\mc{L}=(\pi_*\mc{M})^{**} = (\ol{A}_1\mc{O}_{\wtY})^{**}$. Thus
part (2.iii) of the theorem  holds. \qed

As an application of the theorem, we get the following criterion for
a noetherian cg domain to be a twisted homogeneous coordinate ring.
 Up to changes in phraseology  (the set $\Splusex$ does not appear in \cite{RZ})
 this result is also  a special case of \cite[Theorem~4.4]{RZ}.

\begin{corollary}\label{Sinfty is zero}
 Let $A$ be a
noetherian cg domain that is generated in degree one and  \btwog.
Suppose that the set of extremal elements
 $\Splusex$, as defined in \eqref{extremal-defn},  is empty. Then,
 up to a finite dimensional vector space, $A\cong B(\wtY,\mc{L},\sigma)$, where
 $\wtY = \YY_n$   for any   stable value of $n$.
 \end{corollary}

\begin{proof}
The first paragraph of the proof of  Theorem~\ref{blowing  down}
shows that $\wtY=\YY = \YY_n.$ Therefore, the set $P$ constructed in
the proof of Theorem~\ref{mainthm} is also empty and so
Corollary~\ref{surj-cor}---or indeed \cite[Theorem~4.1]{AS1}---shows
that $A\cong R(\wtY,\emptyset,\mc{L},\sigma)=B(\YY,\mc{L},\sigma)$
up to a finite dimensional vector space.
 \end{proof}


\section{Examples}\label{examples-section}

In  this section, we give two  examples which illustrate the earlier
results in the paper.   The first example shows that the truncated point schemes $\XX_n$
and the relevant components $\YY_n$ are, in general,  distinct; indeed the
 dimension can differ by  arbitrarily large amounts.
Surprisingly, this happens when one takes the \nsr\ at a high power of a sheaf of
 maximal ideals.

\begin{example}\label{big-xn}
  Let $\XXX$ be an integral projective scheme  with $\dim \XXX \geq 2$ and $
\sigma\in \mathrm{Aut}\,\XXX$ for which  there exists a very ample
and $\sigma$-ample   invertible sheaf $\mc{L}$  and a point $c \in
\XXX$ such that the orbit $\langle \sigma\rangle \cdot c$ is
critically dense. Write $\mc{I}=\mc{I}_c$ for the ideal sheaf of the
point $c$. Then, for any integer $t\geq 1$ there exist $p,q\in
\mathbb{N}$ with the following properties:
\begin{enumerate}
\item The  ring $R=R(\XXX,Z_{\mc{I}^p},\mc{L}^q,\sigma)$ is noetherian and
generated in degree one.
\item
For any $n \geq 1$ the truncated point scheme $\XX_{n}$ for $R$
satisfies $\dim \XX_{n}\geq t$.
\item   For any $n \geq 1$ there exists a point $w\in \XX_n$ such that the fibre $\phi_n^{-1}(w)$
has dimension   $\geq t$.
\end{enumerate}
\end{example}

\begin{remarks}\label{big-xn-remark}
 (1) Note that, if $\XXX$ is a surface in this example and $n$ is a stable value
 then  $\dim \YY_n=2$, by Corollary~\ref{zero-ideal}. Therefore, if $t>2$, then  $\XX_n\not=\YY_n$
for all $n \gg 0$.

(2) An explicit example can be obtained by using the data from Section~\ref{first-eg};
thus, modulo a slight change of notation, take $\XXX=\mb{P}^2$ with  $\mc{L}=\mc{O}(1)$,
 $c=(1:1:1)$    and   $\sigma$  defined by $(\lambda_0: \lambda_1:\lambda_2)\mapsto
(\lambda_0: \alpha \lambda_1: \beta\lambda_2)$, where
$\alpha,\beta\in k^*$  are algebraically independent over the prime
subfield of $k$.

(3) This result is in striking contrast to the  commutative case, where
blowing up at a maximal ideal $\mc{I}$
and at some power $\mc{I}^p$ of $\mc{O}_{\XXX}$ gives the same  scheme
 \cite[Exercise~II.7.11]{Ha}.
\end{remarks}

\begin{proof} In the proof, tensor products will be taken over
$\mc{O}_{\XXX}$. As $\dim \XXX \geq 2$ we can pick $p$ such that
length$\,\mc{I}^p/\mc{I}^{p+1}\geq t+1$. Set $\mc{J}=\mc{I}^p$.
  We  want to pick $q$ so that there are  no complications arising from  cohomology.
  Explicitly, by   \cite[Proposition~3.20 and Corollary~3.21]{RS}
 we can choose $q\geq 1$ such that:
\begin{enumerate}
  \item[(i)] the  ring $R(\XXX,Z_{\mc{J}},\mc{L}^q,\sigma)$ is   generated in
degree one, and
  \item[(ii)] $\mathrm{H}^1(\XXX,\, \mc{I}^{\sigma^u}\otimes \mc{J}_a
\otimes\mc{L}_a^q) = 0$,
  for all $u\in \mb{Z}$ and  $a\geq 1$.
  \end{enumerate}
For these values of $p$ and $q$ set
$R=R(\XXX,Z_{\mc{J}},\mc{L}^q,\sigma)$ and  $   \mc{R}=\mc{R}(\XXX,Z_{\mc
{J}},\mc{L}^q,\sigma)$.
  The ring   $R$ is noetherian by
  Proposition~\ref{bimod alg noeth}(2) and so (1) holds.
Let  $\mc{H}=\mc{I}^{\sigma^n}$ for  some $n \geq 1$.
The point of the proof will be to find a $t$-dimensional family of
truncated point modules of length $n+2$ among the quotients of the
$R$-module $M=\mathrm{H}^0(\XXX,\, \mc{R}/\mc{H}\mc{R})$.

Taking cohomology of the exact sequence
  $0\to \mc{H} \mc{R}\to \mc{R}\to   \mc{R}/\mc{H}\mc{R}\to 0$ gives  the
exact sequence  of $R$-modules
  $$0\longrightarrow \mathrm{H}^0(\XXX, \, \mc{H}  \mc{R}) \longrightarrow R \buildrel
{\chi}\over\longrightarrow
  M \longrightarrow
   \bigoplus_{a=0}^\infty \mathrm{H}^1(\XXX,\, \mc{H}  \mc{J}_a\otimes
\mc{L}_a^q).$$  Suppose first that $a\geq 1$; thus
   $ \mathrm{H}^1(\XXX,\, \mc{H}\otimes \mc{J}_a\otimes \mc{L}_a^q)=0$,
by the choice of $q$. Since the kernel $\mc{D}$ of the natural
surjection $\mc{H}\otimes \mc{J}_a\otimes \mc{L}_a^q
\twoheadrightarrow \mc{H} \mc{J}_a\otimes \mc{L}_a^q$ is supported
at the point $\sigma^{-n}(c)$,
certainly
 $\mathrm{H}^2(\XXX,\,\mc{D})=0$.
Combining these observations shows that $ \mathrm{H}^1(\XXX,\, \mc{H}
\mc{J}_a\otimes \mc{L}_a^q)=0$ as well and so
  $\chi_a:R_a\to M_a$ is surjective for all $a\geq 1$.
  When $a=0$, $R_0=M_0=k$ and
$\mathrm{H}^0(\XXX,\, \mc{H} \mc{R}_0)=\mathrm{H}^0(\XXX,\, \mc{H})=0$.
Thus $\chi_0$ is  an isomorphism, $\chi$ is surjective   and $M$ is
cyclic.

 For any $r\geq 0$,
$$M_r \ = \ \mathrm{H}^0\left(\XXX,\, \frac{ \mc{J}_r\otimes \mc{L}^q_r}
{\mc{H}   \mc{J}_r\otimes \mc{L}^q_r} \right) \  \cong \
\mathrm{H}^0\left(\XXX,\, \frac{ \mc{J}_r} {\mc{H}   \mc{J}_r } \right)
.$$
If $r\leq n$ then $\mc{H}=\mc{I}^{\sigma^n}$ and
$\mc{J}_r=\mc{J}\cdots \mc{J} ^{\sigma^{r-1}}$ are comaximal and so
$M_r\cong  \mathrm{H}^0(\XXX,\,   \mc{O}_X/\mc{H})\cong k .$ However,
$$M_{n+1}\  \cong\
  \mathrm{H}^0\left(\XXX,\, \frac{\mc{J}\cdots
  \mc{J} ^{\sigma^{n}}}{\mc{H} \mc{J}\cdots \mc{J} ^{\sigma^{n}}} \right)
\ \cong \   \mathrm{H}^0\left(\XXX,\, \frac{\mc{H}^{p}}{\mc{H}^{p+1}}
\right).$$ This has dimension $T \geq t+1$  by the choice of $p$.
For each $(T-1)$-dimensional subspace $V\subset M_{n+1}$, let
$M(V)=M/(VR+ M_{n+2}R) \cong R/Q(V)$, say. By the above analysis,
this is a cyclic $R$-module with $\dim_k M(V)_j=1$ for $0\leq j\leq
n+1$; in other words $M(V)$ is a truncated point module of length
$n+2$. Now, truncated point modules $R/Q$ and $R/Q'$ are isomorphic
if and only if $Q=Q'$. Since  $Q(V_1)\not=Q(V_2)$ if $V_1\not=V_2$, we
have therefore constructed a family of truncated point modules of
length $n+2$ parametrized by $\mb{P}(M_{n+1}^*)$.  Thus $\dim \XX_{n+1} \geq
t$ for all $n \geq 0$, proving (2). Finally, each of the modules
$M(V)$ has the same truncation $M/M_{n+1}R$ and so, by
Lemma~\ref{right left duality}(1),   the fibre of
$\phi_n$ over the point $w \in \XX_n$ corresponding to $M/M_{n+1}R$
also has dimension at least $t$.
  \end{proof}

In order to illustrate the questions raised in
Section~\ref{subsection-questions}, we end by giving
an example  of a (non-noetherian) cg domain
$A$ that is  \bbc\ but not \btwog. As mentioned there,
we conjecture that examples of this type are always non-noetherian.

\begin{example}\label{non-noeth}
Let $K = k(u,v)$ be a rational function field over $k$,
and define $\sigma: K \to K$ by  $u
\mapsto u$, $v \mapsto uv$.  Let $Q = K[ z, z^{-1}; \sigma]$ and
$A = k \langle uz, vz, z \rangle \subset Q$.
Then $A$ is a non-noetherian domain that is \bbc\ but not \btwog.
\end{example}

\begin{proof}
It is a routine exercise, which we leave to the reader,  to show
 that  $\GKdim A=4$ (mimic the proof of   \cite[Example~8.2.16]{MR}).
 Thus $A$ is an Ore domain by \cite[Corollary~8.1.21]{MR}
  and since $u,v \in A_1 A_1^{-1}$ it is clear that $Q$ is the graded
ring of fractions of $A$.

In order to prove that  $A$ is not noetherian, write   $f_n = (vz) z^n =vz^{n+1}$
for  $n \geq 0$   and let $I^{(n)} = \sum_{i=0}^n Af_i$.
Note that any element $a\in A_{ r}$ has the form
 $a=qz^{r}$, for some $q\in k[u,v]$. Hence, if $r>0$ then
 $af_{n-r}=qz^{r}vz^{n-r+1} = qu^rvz^{n+1}\in k[u,v]uf_n$.
Thus,  $f_{n} \not \in I^{(n-1)}$ for $n\geq 1$
 and  $I^{(0)} \subsetneq I^{(1)} \subsetneq
\dots $ gives a proper ascending chain of left ideals of $A$.

Clearly $A$ is \bbc.  There are several ways to see that $A$ is not
\btwog. First of all, \cite[Remark 7.3]{DF} shows  that
there is no projective model $Y$ of $K = k(u,v)$ with an
automorphism inducing the given $\sigma \in \aut(K)$, so $A$ is not
\btwog\ by definition. Alternatively  if $A$ were \btwog,
 then  \cite[Theorem~1.6]{RZ}  would imply that
  $\GKdim(A)=\GKdim(B)$ for some twisted homogeneous
coordinate ring $B(Y, \mc{L}, \sigma)$ with $Q(B) = K[t, t^{-1};
\sigma]$, where   $Y$ would necessarily be a surface and $\mc{L}$ would be $\sigma$-ample.
  But, by  \cite[Theorem~1.7]{AV}, no such twisted homogeneous coordinate ring has
GK-dimension $4$.
\end{proof}

\section*{Acknowledgments}
We thank James Zhang, Mike Artin, and Johan de Jong for helpful
conversations.

\section*{Index of Notation}\label{index}
 \begin{multicols}{2}
{\small  \baselineskip 14pt

   bimodule algebra $\mc{R}(\XXX,Z_{\mc{I}},\mc{L},\tau)$, \hfill \pageref{bimodule-defn}

\btwog,\hfill\pageref{btwog-defn}

\bbc,   \hfill\pageref{bbc-defn}

$\chi_n=\phi_n$ or $\psi_n$, \hfill \pageref{chi-defn}

 contracted points $\Sminus_n$ and $\Splus_n$
 \hfill \pageref{contracted-defn}

critically dense set, \hfill \pageref{critical-defn}

 extremal elements $\Splusex,\Sminusex$, \hfill \pageref{extremal-defn}, \pageref{extremal-defn2}

   function ring $Q(A)_0$, \hfill \pageref{fn-ring-defn}

   $\mc{K}$, the sheaf of rational functions on $Y$,\hfill\pageref{mcK-defn}

  $\LL_n$, \hfill\pageref{Ln-defn}

\nsr\  $R(\XXX,Z_{\mc{I}}, \mc{L},\tau)$,  \hfill   \pageref{naive-defn}

\naive\ blowup, $\rqgr \mc{R}(\XXX,Z_{\mc{I}}, \mc{L},\tau)$,  \hfill   \pageref{naive-defn}

$\Phi_n, \Psi_n :\XX_{n+1}\to \XX_n$,  \hfill \pageref{phi-defn}

$\Phi_{m,n},\Psi_{m,n} : \XX_m\to \XX_n$,  \hfill\pageref{phi-mn}

$\phi_{m,n},\psi_{m,n} : \YY_m\to \YY_n$,  \hfill\pageref{y-defn}

$\mb{P}=\mb{P}(A_1^*)$,\hfill\pageref{pee-defn}

$\PP(R), \PP^\ell(R)$, set of point modules, \hfill \pageref{point-defn}

$\PP_r(R), \PP_r^\ell(R)$, set of  truncated point modules,
\hfill \pageref{trunc-defn}

$Q(A)=D[z,z^{-1}; \tau]$, graded quotient ring  \hfill \pageref{fn-ring-defn}

 qgr, quotient category, \hfill \pageref{qgrade-defn},\pageref{script-Rmod-defn}

{\rel point scheme data} $(\YY_n, \phi_n, \psi_n)$, \hfill\pageref{y-defn}

 $\Splus=\Splus_n,\Sminus =\Sminus_n, \YY=\YY_n$ for fixed $n\gg 0$, \hfill\pageref{S-notation}

 $\sigma=\sigma_n=\psi_n\phi^{-1}_n: \YY \to \YY$, \hfill\pageref{sigma-notation}

$\sigma$-ample sheaf, \hfill\pageref{ample def}

saturating set, \hfill\pageref{saturating}

shift of grading $M\to M[1]$,\hfill\pageref{shift-defn}

stable value of $n$,   \hfill\pageref{stable-defn2}

$\alpha^\star : \spec S\to \spec R$, induced map of schemes,\hfill\pageref{star-defn}

truncations $\Phi'_n(M), \Psi'_n(M)$ of $M$,\hfill\pageref{rtrunc}

twisted coordinate ring  $B(\XXX,\LL,\tau)$,\hfill\pageref{twisted-defn}

$(\XX_n, \Phi_n,\Psi_n)$,  point scheme data,  \hfill\pageref{X-data-defn}

$(\XX^{\ell}_n, \Phi^{\ell}_n, \Psi^{\ell}_n)$,   left point scheme data, \hfill\pageref{left-defn}}

\end{multicols}

\providecommand{\bysame}{\leavevmode\hbox to3em{\hrulefill}\thinspace}
\providecommand{\MR}{\relax\ifhmode\unskip\space\fi MR }
\providecommand{\MRhref}[2]{
  \href{http://www.ams.org/mathscinet-getitem?mr=#1}{#2}
} \providecommand{\href}[2]{#2}

\end{document}